\documentclass[10pt]{article}
\usepackage{amsfonts,color}
\usepackage{amssymb}
\usepackage{footnote}
\usepackage{mathtools}
\usepackage{amsthm,wasysym}
\usepackage[english]{babel}
\usepackage{bbm}
\usepackage{colonequals}
\usepackage{setspace}
\usepackage{titlesec}
\usepackage{float}

\usepackage{graphicx}
\usepackage{amsfonts}
\usepackage{hyperref}
\usepackage{subcaption}
\usepackage{amsmath}
\usepackage{nicefrac}
\usepackage{gensymb}
\usepackage{enumerate}
\usepackage[nameinlink,capitalize]{cleveref}
\usepackage{chngcntr}

\usepackage{cite}
\usepackage[nottoc,numbib]{tocbibind}

\usepackage{pgf,tikz,pgfplots}
\pgfplotsset{compat=1.14}
\usepackage{mathrsfs}
\usetikzlibrary{arrows}
\usetikzlibrary{arrows.meta}

\DeclareMathOperator{\at}{\bigg\vert}
\newcommand{\vb}[1]{\mathbf{#1}}
\newcommand{\bm}[1]{\boldsymbol{#1}}

\DeclareMathOperator{\range}{\mathrm{range}}

\DeclareMathOperator{\sgn}{\mathrm{sgn}}
\DeclareMathOperator{\sym}{\mathrm{sym}}
\DeclareMathOperator{\spa}{\mathrm{span}}
\DeclareMathOperator{\vol}{\mathrm{vol}}
\DeclareMathOperator{\dev}{\mathrm{dev}}
\DeclareMathOperator{\skw}{\mathrm{skew}}
\DeclareMathOperator{\Anti}{\mathrm{Anti}}
\DeclareMathOperator{\tr}{\mathrm{tr}}

\DeclareMathOperator{\cof}{\mathrm{cof}}
\newcommand{\jump}[1]{\ensuremath{[\![#1]\!]} }
\newcommand{\one}{\bm{\mathbbm{1}}}
\newcommand{\con}[2]{\langle {#1} , \, {#2} \rangle}
\newcommand{\norm}[1]{\| {#1} \|}

\newcommand{\dd}{\mathrm{d}}
\newcommand{\D}{\mathrm{D}}
\DeclareMathOperator{\di}{\mathrm{div}}
\DeclareMathOperator{\Di}{\mathrm{Div}}

\DeclareMathOperator{\curl}{\mathrm{curl}}
\DeclareMathOperator{\Curl}{\mathrm{Curl}}

\newcommand{\SO}{\mathrm{SO}}
\newcommand{\so}{\mathfrak{so}}
\newcommand{\Sym}{\mathrm{Sym}}
\newcommand{\Dev}{\mathfrak{sl}}
\newcommand{\Diag}{\mathrm{Diag}}

\newcommand{\Ned}{\mathcal{N}_{I}}

\newcommand{\Nedtwo}{\mathcal{N}_{II}}
\newcommand{\Lag}{\mathcal{L}}

\newcommand{\Po}{\mathit{P}}
\newcommand{\Q}{\mathit{Q}}
\newcommand{\Le}{{\mathit{L}^2}}

\newcommand{\Hone}{\mathit{H}^1}

\newcommand{\HsD}[1]{\mathit{H}^\mathrm{sym}(\mathrm{Div}{#1})}
\newcommand{\Hd}[1]{\mathit{H}(\mathrm{div}{#1})}

\newcommand{\Hc}[1]{\mathit{H}(\mathrm{curl}{#1})}
\newcommand{\HC}[1]{\mathit{H}(\mathrm{Curl}{#1})}

\newcommand{\HdD}[1]{\mathit{H}(\di\mathrm{Div}{#1})}
\newcommand{\HsC}[1]{\mathit{H}(\mathrm{sym}\,\mathrm{Curl}{#1})}
\newcommand{\HdsC}[1]{\mathit{H}(\mathrm{dev}\,\mathrm{sym}\,\mathrm{Curl}{#1})}
\newcommand{\HsCd}[1]{\mathit{H}^{\mathrm{dev}}(\mathrm{sym}\,\mathrm{Curl}{#1})}

\newcommand{\elem}{T}
\newcommand{\body}{V}
\newcommand{\surf}{A}
\newcommand{\curv}{s}

\newcommand{\R}{\mathbb{R}}
\newcommand{\U}{\mathit{U}}
\newcommand{\C}{\mathit{C}}

\renewcommand{\H}{\mathit{H}}
\renewcommand{\S}{\mathcal{S}}
\newcommand{\Y}{\mathcal{Y}}

\newcommand{\Dp}{\mathcal{D}}
\newcommand{\M}{\mathcal{M}}

\newcommand{\tem}{\mathcal{T}}
\newcommand{\phys}{\mathcal{P}}
\newcommand{\ver}{\mathcal{V}}
\newcommand{\edge}{\mathcal{E}}
\newcommand{\face}{\mathcal{F}}
\newcommand{\cell}{\mathcal{C}}

\newcommand{\J}{\mathbb{J}}

\newcommand{\lame}{\lambda_{\mathrm{e}}}
\newcommand{\lammi}{\lambda_{\mathrm{micro}}}
\newcommand{\lamma}{\lambda_{\mathrm{macro}}}
\newcommand{\mue}{\mu_{\mathrm{e}}}
\newcommand{\muc}{\mu_{\mathrm{c}}}
\newcommand{\mumi}{\mu_{\mathrm{micro}}}
\newcommand{\muma}{\mu_{\mathrm{macro}}}
\newcommand{\Lc}{L_\mathrm{c}}
\newcommand{\Ce}{\mathbb{C}_{\mathrm{e}}}
\newcommand{\Cc}{\mathbb{C}_{\mathrm{c}}}
\newcommand{\Cmic}{\mathbb{C}_{\mathrm{micro}}}

\newcommand{\wrt}{\text{w.r.t.}}

\newcommand{\ud}{\vb{u}}
\newcommand{\Pm}{\bm{P}}

\allowdisplaybreaks[1]

\makeindex

\newtheoremstyle{break}
{\topsep}{\topsep}%
{\itshape}{}%
{\bfseries}{}%
{\newline}{}%
\theoremstyle{break}

\newtheorem{theorem}{Theorem}
\newtheorem{lemma}{Lemma}

\newtheorem{remark}{Remark}
\newtheorem{observation}{Observation}
\newtheorem{definition}{Definition}

\counterwithin{theorem}{section}
\counterwithin{lemma}{section}
\counterwithin{remark}{section}
\counterwithin{observation}{section}
\counterwithin{definition}{section}

\counterwithin{figure}{section}
\counterwithin{equation}{section}

\makeatletter
\let\@fnsymbol\@arabic
\makeatother

\crefname{Problem}{Problem.}{Problem.}

\title{Novel $\HsC{}$-conforming finite elements for the relaxed micromorphic sequence}

\author{\normalsize{Adam Sky}\thanks{Corresponding author: Adam Sky, Institute of Computational Engineering and Sciences, Department of Engineering, Faculty of Science, Technology and Medicine, University of Luxembourg, 6, Avenue de la Fonte, L-4362 Esch-sur-Alzette, Luxembourg, email: adam.sky@uni.lu}
	, \quad
	\normalsize{Michael Neunteufel}\thanks{Michael Neunteufel, Institute of Analysis and Scientific Computing, Technische Universit\"at Wien, Wiedner Hauptstr. 8-10 , 1040 Wien, Austria, email: michael.neunteufel@tuwien.ac.at}
	, \quad
    \normalsize{Peter Lewintan}\thanks{Peter Lewintan, Chair for Nonlinear 
		Analysis and Modelling, Faculty of Mathematics, Universit\"{a}t Duisburg-Essen,
		Thea-Leymann Str. 9, 45127 Essen, Germany, email: peter.lewintan@uni-due.de}
    , \quad
    \normalsize{Andreas Zilian}\thanks{Andreas Zilian, Institute of Computational Engineering and Sciences, Department of Engineering, Faculty of Science, Technology and Medicine, University of Luxembourg, 6, Avenue de la Fonte, L-4362 Esch-sur-Alzette, Luxembourg, email: andreas.zilian@uni.lu}
	\quad
    \\
	\normalsize{and} \quad
	\normalsize{Patrizio Neff}\thanks{Patrizio Neff, Chair for Nonlinear 
		Analysis and Modelling, Faculty of Mathematics, Universit\"{a}t Duisburg-Essen,
		Thea-Leymann Str. 9, 45127 Essen, Germany, email: patrizio.neff@uni-due.de}
}

\begin{document}

\maketitle

\begin{abstract}
In this work we construct novel $\HsC{}$-conforming finite elements for the recently introduced relaxed micromorphic sequence, which can be considered as the completion of the $\di\Di$-sequence with respect to the $\HsC{}$-space. 
The elements respect $\HC{}$-regularity and their lowest order versions converge optimally for $[\HsC{}\setminus \HC{}]$-fields.
This work introduces a detailed construction, proofs of linear independence and conformity of the basis, and numerical examples.
Further, we demonstrate an application to the computation of metamaterials with the relaxed micromorphic model.
\\
\vspace*{0.25cm}

\end{abstract}


\section{Introduction}
The deviatoric $\HsCd{,\body}$ space has been introduced in \cite{PaulyDiv} as part of the $\di \Di$-complex in the context of biharmonic equations, for which discrete sequences can be found in \cite{CRMECA_2023__351_S1_A8_0,Hu2021,Chen2022,HuDivDIV,dipietro2023discrete}. Moreover, in \cite{arnold_complexes_2021} the authors show that many Hilbert space complexes \cite{PaulyDeRham,PaulyEl,PaulyOberwolfach,ChenStokes,Angoshtari,Cap,Neilan,pauly_elasticity_2022,hu2023nonlinear} are in fact related to each other through various operators. In contrast to the $\di \Di$-complex, in \cite{Lewintan2021} the full $\HsC{,\body}$-space has been defined in the context of the relaxed micromorphic model \cite{Neff2014,SKY2022115298,KNEES2023126806,knees2023global}, where the model has been shown to maintain well-posedness in this larger space due to the generalised incompatible Korn-type inequalities \cite{Lewintan2021,LewintanInc,LewintanInc2,NEFF20151267,Neff2012,Gmeineder1,Gmeineder2}. 
The relaxed micromorphic model is expressed as a two-field minimisation problem via the energy functional 
\begin{align}
		I(\ud, \Pm) = \dfrac{1}{2} \int_\body &\langle\Ce \sym(\D \ud - \Pm) , \,  \sym(\D \ud - \Pm) \rangle + \langle \Cmic \sym\Pm , \,  \sym \Pm \rangle \notag \\
		& + \langle \Cc\skw(\D \ud - \Pm) , \,  \skw(\D \ud - \Pm) \rangle + \muma\Lc^2 \langle \mathbb{L} \sym \Curl \Pm , \,  \sym \Curl \Pm \rangle \, \dd \body \notag \\ 
		& \qquad - \int_\body \langle \ud , \, \vb{f} \rangle + \langle \Pm , \, \bm{M} \rangle \, \dd \body \to  \min \quad \wrt \quad \{\ud, \Pm\}
		\, ,
  \label{eq:rmm}
	\end{align}
 where the $\sym \Curl$-operator for second order tensors is defined via
    \begin{align}
    	2\sym \Curl \Pm &=  \begin{bmatrix}
    	2(P_{13,y} - P_{12,z}) & P_{23,y} - P_{22,z} + P_{11,z} - P_{13,x} & P_{12,x} - P_{11,y} + P_{33,y} - P_{32,z} \\
    	P_{23,y} - P_{22,z} + P_{11,z} - P_{13,x} & 2(P_{21,z} - P_{23,x}) & P_{31,z} - P_{33,x} + P_{22,x} - P_{21,y} \\
    	 P_{12,x} - P_{11,y} + P_{33,y} - P_{32,z} & P_{31,z} - P_{33,x} + P_{22,x} - P_{21,y} & 2(P_{32,x} - P_{31,y})
    \end{bmatrix} \, ,
    \end{align}
    and the domain $\body \subset \R^{3}$ denotes a general bounded, open, non-emtpy set.
	The displacement field and the microdistortion field are functions of the domain
	\begin{align}
		&\ud : \overline{\body} \subset \R^3 \to \R^3 \, , && \Pm : \overline{\body} \subset \R^3 \to \R^{3 \times 3} \, .
	\end{align}
	The tensors $\Ce, \Cmic \in \R^{3\times 3 \times 3 \times 3}$ are standard positive definite fourth order meso- and micro-elasticity tensors. For isotropic materials they are given by
	\begin{align}
		&\Ce = 2 \mue \, \J  +  \lame \one \otimes \one \, , &&
		\Cmic = 2 \mumi \, \J  + \lammi \one \otimes \one  \, ,
	\end{align}
    where $\one \in \R^{3\times 3}$ is the second order identity tensor and $\J\in \R^{3\times 3\times 3 \times 3}$ is the fourth order identity tensor. The macroscopic shear modulus is designated as $\muma$.
	The fourth order tensor $\Cc \in \R^{3 \times 3 \times 3 \times 3}$ is a positive semi-definite material tensor related to Cosserat micro-polar continua \cite{GhibaCosserat,Russo2023,AltenbachCosserat}, which accounts for infinitesimal rotations $\Cc: \so(3) \to \so(3)$, where $\so(3)$ is the space of skew-symmetric matrices. For isotropic materials there holds $\Cc = 2\muc \, \J$, where $\muc \geq 0$ is the Cosserat couple modulus. Lastly, $\mathbb{L} \in \R^{3 \times 3 \times 3 \times 3}$ is a positive definite fourth order tensor of weights for the characteristic length scale parameter $\Lc$, which is motivated by the geometry of a possible microstructure. For simplicity, we assume $\mathbb{L}=\J$ in the following.
	The forces and micro-moments are given by $\vb{f}$ and $\bm{M}$, respectively.
The balance equations are derived by variation and partial integration, reading
\begin{subequations}
		\begin{align}
			-\Di[\Ce \sym (\D \vb{u} - \bm{P}) + \Cc \skw (\D \vb{u} - \bm{P})] &= \vb{f} && \text{in} \quad \body \, , \label{eq:strong_u} \\
			-\Ce  \sym (\D \vb{u} - \Pm) - \Cc  \skw(\D \vb{u} - \Pm) + \Cmic \sym \Pm + \muma \, \Lc ^ 2  \Curl (\mathbb{L} \sym \Curl\Pm) &= \bm{M} && \text{in} \quad \body \, . \label{eq:strong_p}  
		\end{align}
	    \label[Problem]{eq:full_relaxed}
	\end{subequations}
We note that an alternative version of the model employs the $\HC{,\body}$-space \cite{Neff2014}, in which case the entire Curl of the microdistortion $\Curl \Pm$ appears in the energy functional, instead of just the symmetric part $\sym \Curl \Pm$.

In this work we are interested in constructing conforming finite elements for the complete $\HsC{,\body}$ space in the context of the relaxed micromorphic model. 
As shown in \cite{Lewintan2021}, the $\HsC{,\body}$-space is strictly larger than $\HC{,\body}$ and in general there holds
\begin{align}
    [\Hone(\body)]^{3\times 3} \subsetneq \HC{,\body} \subsetneq \HsC{,\body} \, .
\end{align}
Consequently, conforming finite elements for $[\Hone(\body)]^{3 \times 3}$ and $\HC{,\body}$, such as Lagrange and N\'ed\'elec \cite{Ned2,Nedelec1980}, are clearly conforming in $\HsC{,\body}$ and represent possible candidates for computations in the space. However, since the $\HsC{,\body}$-space is strictly larger, there exist solutions not belonging to $[\Hone(\body)]^{3 \times 3}$ or $\HC{,\body}$ but to $\HsC{,\body}$, such that for these, suboptimal convergence is observed, as demonstrated in \cite{SkyOn}.
We note that in \cite{Lewintan2021}, the space $\HdsC{,\body}$ has also been introduced and shown to be strictly larger than the $\HsC{,\body}$-space such that
\begin{align}
    [\Hone(\body)]^{3\times 3} \subsetneq \HC{,\body} \subsetneq \HsC{,\body} \subsetneq \HdsC{,\body} \, .
\end{align}
Nonetheless, the generalised Korn-inequality (needed for well-posedness in the case of $\Cc = 0$) continues to hold also with the very weak correction term $\dev\sym\Curl\Pm$, see \cite{Lewintan2021}. Consequently, requiring the control of only $\norm{\dev\sym\Curl P}_{\Le}$ allows to relax the considered variational problem even further, while consistently maintaining well-posedness. The latter demonstrates the decisive role that such coercivity estimates play in the deduction of sound mathematical models. Sharp criteria for the validity of such Korn-Maxwell-Sobolev type inequalities were given in the recent works \cite{Gmeineder0, Gmeineder1,Gmeineder2} and the references contained therein. Nevertheless, Observation 2.3 from \cite{Lewintan2021} states
\begin{align}
    &\dev \sym[\Pm \times \vb{v}] = 0 \qquad \iff \qquad \sym[\Pm \times \vb{v}] = 0 \qquad \forall \, \{\Pm, \vb{v} \} \in \R^{3 \times 3} \times \R^3 \, ,   
\end{align}
implying that for the corresponding traces there holds
\begin{align}
    &\tr_{\HdsC{}}\Pm = 0 \qquad \iff \qquad \tr_{\HsC{}}\Pm = 0 \, ,  
\end{align}
such that the spaces $\HsC{,\body}$ and $\HdsC{,\body}$ possess the same regularity for discretisations. Thus, finite elements for the $\HsC{,\body}$-space are automatically $\HdsC{,\body}$-conforming and vice versa, although the kernel of the $\HdsC{,\body}$-space is in fact larger \cite{Lewintan2021}. Consequently, for complete polynomial spaces $[\Po^p(\body)]^{3 \times 3}$ there cannot exist an $\HdsC{,\body}$-conforming finite element with a larger polynomial kernel than for a conforming $\HsC{,\body}$-discretisation and the elements coincide. 

To the authors knowledge, this work introduces the first conforming finite elements for $\HsC{,\body}$, which also respect the minimal tangential regularity of $\HC{,\body}$. Alternative constructions with higher regularity are found in \cite{SkyOn,Hu2021,Chen2022,HuDivDIV}. 
In these works the authors rely on vertex degrees of freedom to define the $\HsC{,\body}$-subspace. The latter imposes a higher regularity on the subspace than the one exhibited by N\'ed\'elec elements \cite{Nedelec1980,Ned2,Joachim2005,Zaglmayr2006,sky_polytopal_2022,Anjam2015,sky_higher_2023}, which are used to discretise $\HC{,\body}$. This is of particular importance to the relaxed micromorphic model, since the model is multi-scale in nature. In fact, as shown in \cite{sky_hybrid_2021,sky_higher_2023,SKY2022115298,Schroder2022,Sarhil2023,SkyPamm}, the characteristic length scalar parameter $\Lc$ allows the model to interpolate between macro and micro reactions, which are governed by linear elasticity models. 
In the limit of the characteristic length scale parameter $\Lc \to +\infty$ the microdistortion field $\Pm$ is expected to degenerate into a gradient field $\Pm = \D \vb{v}$ for some vector $\vb{v}$, such that it is compatible with the gradient of the displacement field $\D \vb{u}$. This is possible for 
\begin{align}
    \D \vb{u} \in \HC{,\body} \subset \HsC{,\body} \ni \Pm \, , 
\end{align}
due to $\HsC{,\body}$ being the larger space, but impossible for
\begin{align}
    \Pm \in [\Hone(\body) \otimes \Dev(3)] \oplus [\Le(\body) \otimes \one] \nsupseteq \HC{,\body} \, , 
\end{align}
which is the construction used in \cite{SkyOn}. The latter imposes deviatoric $\C^0(\body)$-continuity at the vertices and allows the identity tensor $\one$ to jump in the discrete subspace. As such, it is not conforming in $\HC{,\body}$. Consequently, such a construction may induce locking in the limit of $\Lc \to + \infty$, analogously to how shear-locking in the Reissner-Mindlin plate \cite{sky2023reissnermindlin} is due to incompatibility of the discrete spaces for a vanishing thickness $t \to 0$. Further, the $\HC{,\body}$-version of the relaxed micromorphic model introduces the consistent coupling condition \cite{dagostino2021consistent}
\begin{align}
    &\Pm \times \vb{n} = \D \widetilde{\vb{u}} \times \vb{n} && \text{on} && \surf_D \, ,
\end{align}
which controls the Dirichlet boundary $\surf_D^P = \surf_D^u = \surf_D$ of the microdistortion field. 
Again, the condition requires the consistency of the Sobolev trace spaces \cite{HIPTMAIR2023109905,DINEZZA2012521} $
    \tr_{\HC{}} \D \widetilde{\vb{u}} \in \tr_{\HC{}}[ \HC{,\body} ]  \, ,
$
which is not satisfied by a $\C^0(\body)$-construction. 
In the $\HsC{,\body}$-version of the model the consistent coupling condition is adjusted to $\sym(\Pm \times \vb{n}) = \sym(\D \widetilde{\vb{u}} \times \vb{n})$. There holds
\begin{align}
    \tr_{\HC{}}[ \HC{,\body} ] \subseteq \tr_{\HsC{}}[\HsC{,\body}] \, , 
\end{align}
implying a weaker imposition of Dirichlet boundary conditions which a consistent discretisation must satisfy.

This work is organized as follows: Firstly we introduce the relevant Hilbert spaces and operators used in this paper. Secondly, we discuss the $\di\Di$-sequence and its completion towards the relaxed micromorphic sequence. Next, we construct finite elements for $\HsC{,\body}$. The elements are then benchmarked via numerical examples showcasing convergence results. Lastly, we discuss our conclusions and outlook.

The low order elements are implemented in the open source finite element software NGSolve\footnote{https://ngsolve.org/} \cite{Sch2014,Sch1997} and are available as supplementary material to this paper\footnote{https://github.com/Askys/NGSolve\_HsymCurl}.
\subsection{Notation}
The following notation is used throughout this work.
Exceptions to these rules are made clear in the precise context.
\begin{itemize}
    \item vectors are defined as bold lower-case letters $\vb{v}, \, \bm{\xi}$
    \item matrices are bold capital letters $\bm{M}$
    \item fourth-order tensors are designated by the blackboard-bold format $\mathbb{A}$
    \item we designate the Cartesian basis as $\{\vb{e}_1, \, \vb{e}_2, \, \vb{e}_3\}$  
    \item the angle-brackets are used to define scalar products of arbitrary dimensions $\con{\vb{a}}{\vb{b}} = a_i b_i$, $\con{\bm{A}}{\bm{B}} = A_{ij}B_{ij}$
    \item the matrix product is used to indicate all partial-contractions between a higher-order and a lower-order tensor $\bm{A}\vb{v} = A_{ij} v_j \vb{e}_i$, $\mathbb{A}\bm{B} = A_{ijkl}B_{kl}\vb{e}_i \otimes \vb{e}_j$
    \item subsequently, we define various differential operators based on the Nabla-operator $\nabla = \partial_i \vb{e}_i$
    \item volumes and surfaces of the physical domain are identified via $\body$ and $\surf$, respectively. Their counterparts on the reference domain are $\Omega$ and $\Gamma$
    \item the volume of elements on the domain is denoted by $\elem$
    \item tangent, cotangent and normal vectors on the physical domain are designated by $\vb{t}$, $\vb{m}$ and $\vb{n}$, respectively. On the reference domain we use $\bm{\tau}$, $\bm{\mu}$ and $\bm{\nu}$ 
    \item in the following, the polytopes of a tetrahedral element are identified with multi-indices, e.g., edge $e_j$ with $j \in \mathcal{J} = \{(1,2),(1,3), \dots\}$  
\end{itemize}

\section{Hilbert spaces}
We introduce the classical Hilbert spaces and their respective norms 
\begin{subequations}
    \begin{align}
    \Le(\body) &=  \{ u: \body \to \mathbb{R} \, | \,  \|u\|_{\Le}^2 < \infty \} \, , & \|u \|_{\Le}^2 &= \int_{\body} \|u\|^2 \, \dd \body \,,\\[1ex]
	\Hone(\body) &= \{ u \in \Le(\body) \, | \,  \nabla u \in [\Le(\body)]^3 \} \, , & \|u \|_{\Hone}^2 &= \|u\|_{\Le}^2 + \| \nabla u \|_{\Le}^2 \,,\\[2ex]
	\Hc{,\, \body} &= \{ \vb{p} \in [\Le(\body)]^3 \, | \curl\vb{p} \in [\Le(\body)]^3   \} \, , & \| \vb{p} \|_{\Hc{}}^2 &= \| \vb{p} \|_{\Le}^2 + \| \curl\vb{p} \|_{\Le}^2 \, ,\\[2ex]
	\Hd{,\, \body} &= \{ \vb{p} \in [\Le(\body)]^3 \, | \di\vb{p} \in \Le(\body)   \} \, , & \| \vb{p} \|_{\Hd{}}^2 &= \| \vb{p} \|_{\Le}^2 + \| \di\vb{p} \|_{\Le}^2 \, .
\end{align}
\end{subequations}
Further, we introduce the matrix-valued Hilbert spaces 
\begin{subequations}
    \begin{align}
    \HC{,\body} &= [\Hc{,\body}]^3 \, ,  \\
    \HsC{,\body} &= \{ \Pm \in [\Le(\body)]^{3 \times 3} \; | \; \sym \Curl \Pm \in [\Le(\body)]^{3 \times 3} \} \, ,  \\
    \HdD{,\body} &= \{ \Pm \in  [\Le(\body)]^{3 \times 3} \; | \; \di\Di \Pm \in \Le(\body) \, , \quad \Pm = \Pm^{T} \} \, ,  \label{eq:HdD}
\end{align}
\end{subequations}
where the space $\HC{,\body}$ is to be understood as a row-wise matrix of the vectorial space $\Hc{,\body}$.
The respective norms of the matrix-valued spaces read
\begin{subequations}
    \begin{align}
    \| \Pm \|_{\HsC{}}^2 &= \| \Pm \|^2_{\Le} + \| \sym \Curl \Pm \|^2_{\Le} \, , \\
    \| \Pm \|_{\HdD{}}^2 &= \| \Pm \|^2_{\Le} + \| \di\Di \Pm \|^2_{\Le} \, .
\end{align}
\end{subequations} 
\begin{remark}[Alternative $\HdD{,\body}$-definition]
    Observe that other works may employ an alternative definition of the $\HdD{,\body}$-space \cite{pechstein_anisotropic_2012,pechstein_analysis_2018,NEUNTEUFEL2021113857}
    \begin{align}
        \HdD{,\body} &= \{ \Pm \in  [\Le(\body)]^{3 \times 3} \; | \; \di\Di \Pm \in \H^{-1}(\body)  \, , \quad \Pm = \Pm^{T}\} \, , 
    \end{align}
    where the regularity of $\di \Di \Pm$ is reduced to the $\H^{-1}(\body)$-setting.
    \label{re:alt_hdD}
\end{remark}
Specifically for the $\HsC{,\body}$-space we also introduce its deviatoric version
\begin{align}
    \HsCd{,\body} = \{ \Pm \in \HsC{,\body} \; | \; \tr \Pm = 0 \} \, ,
\end{align}
which is used in the context of the $\di\Di$-complex.
The matrix-valued differential operators are defined as 
\begin{subequations}
    \begin{align}
    \D \vb{v} &= \vb{v} \otimes \nabla \, , \\
    \Curl \Pm &= - \Pm \times \nabla \, , \\ 
    \sym \Curl \Pm &= -\sym(\Pm \times \nabla) \, , \\
    \di \Di \Pm &= \con{\nabla}{\Pm \cdot \nabla} \, .
\end{align}
\end{subequations}
The Sobolev traces of the $\HC{,\body}$- and $\HsC{,\body}$ spaces read
\begin{align}
	&\tr_{\HC{}} \bm{P} = \Pm  \Anti(\vb{n})^{T} \at_{\Xi} \, , &&  
	\tr_{\HsC{}} \bm{P} = \sym[\Pm  \Anti(\vb{n})^{T}] \at_{\Xi} \, , \label{eq:trsc}
\end{align}
where $\Xi$ represents an arbitrary interface in the domain with a corresponding surface normal vector $\vb{n}$ and 
\begin{align}
    \Anti(\vb{v}) = \begin{bmatrix}
    0 & -v_3 & v_2 \\
    v_3 & 0 & -v_1 \\
    -v_2 & v_1 & 0 
    \end{bmatrix} \, ,
\end{align}
is the anti-symmetric matrix for some vector $\vb{v} \in \mathbb{R}^3$, such that $\Anti(\vb{v})\vb{w} = \vb{v} \times \vb{w}$ with $\vb{w} \in \R^3$.
Lastly, we introduce the spaces 
\begin{align}
    &\Sym(3) \, , && \so(3) \, , && \spa\{\one\} \, , && \Dev(3) \, ,
\end{align}
being the vector-spaces of symmetric, skew-symmetric, volumetric and deviatoric tensors (trace-free), respectively.
The spaces are associated with the algebraic operators \begin{align}
    &\sym\Pm = \dfrac{1}{2} (\Pm + \Pm^{T}) \, , && \skw(\Pm) = \dfrac{1}{2}(\Pm - \Pm^{T}) \, , && \vol \Pm = \dfrac{1}{3} (\tr\Pm) \one \, , && \dev \Pm = \Pm -  \dfrac{1}{3} (\tr\Pm) \one \, ,
\end{align}
such that any second order tensor $\Pm \in \R^{3 \times 3}$ can be orthogonally decomposed into
\begin{align}
    \Pm = \dev \sym \Pm + \skw \Pm + \dfrac{1}{3} (\tr \Pm) \one  \, , 
\end{align}
where $\one$ is the identity tensor.
The operators map a tensor from $\mathbb{R}^{3 \times 3}$ to the respective space.

\section{The $\di\Di \,$- and relaxed micromorphic sequences}
In this section we recall the $\di\Di$- sequence and complete it to the relaxed micromorphic sequence. 
A thorough treatment of the exact $\di\Di$-sequence can be found in \cite{PaulyDiv,CRMECA_2023__351_S1_A8_0,Chen2022,Hu2021,arnold_complexes_2021}. 

The $\di\Di$-sequence is depicted in \cref{fig:diDicomplex} and conveys the identities 
\begin{subequations}
    \begin{align}
    \dev\D [\Hone(\body)]^3 &= \ker(\sym\Curl) \cap \HsCd{,\body} \, , \label{eq:devD} \\ \sym\Curl[\HsCd{,\body }] &= \ker(\di\Di) \cap \HdD{,\body} \, , \label{eq:symCdev} \\ 
    \di \Di \HdD{,\body} &= \Le(\body) \label{eq:surj} \, ,
\end{align}
\end{subequations}
which are exact on contractible domains. 
For the first identity we observe that any element of $\D [\Hone(\body)]^3$ is in $\ker(\sym\Curl)$ due to
    \begin{align}
        \Curl \D \vb{v} = 0 \qquad \forall \,\vb{v} \in  \D [\Hone(\body)]^3 \, ,
    \end{align}
    where the classical identity $\curl \nabla(\cdot)$ is applied row-wise. Due to $\ker(\sym) = \so(3)$ elements not in $\ker(\Curl)$ but in $\ker(\sym\Curl)$ must satisfy
    \begin{align}
        &\Curl \Pm = \bm{A} = \Anti\vb{a} \, , &&\bm{A} \in \so(3) \, ,
    \end{align}
    where $\bm{A} \in \so(3)$ can always be represented by the anti-symmetric matrix of some three-dimensional axial vector $\vb{a}$.
    Taking the row-wise divergence of both sides yields
    \begin{align}
        \Di \Anti \vb{a} = -\curl \vb{a} = 0 \, ,
    \end{align}
    which is satisfied by $\vb{a} = \nabla \lambda$ with $\lambda \in \Hone(\body)$, implying $\Curl \Pm = \Anti(\nabla \lambda)$. Consequently, the algebraic identity 
    \begin{align}
        \Curl(\lambda \one) = \lambda \one \Anti(\nabla)^{T} = \Anti(\nabla \lambda)^{T} \, ,
    \end{align}
    allows to set $\Pm = \lambda \one$ and to determine
    \begin{align}
        \ker(\sym\Curl) \cap \HsC{,\body} = \D[\Hone(\body)]^3 \cup [\Le(\body) \otimes \one ] \, .
    \end{align}
    Clearly, an orthogonal split of the kernel can be achieved by taking only the deviatoric part of the gradients 
    \begin{align}
        \ker(\sym\Curl) \cap \HsC{,\body} = \dev\D[\Hone(\body)]^3 \oplus [\Le(\body) \otimes \one] \, ,
        \label{eq:ker}
    \end{align}
    and since the space is restricted to trace-free tensors $\HsCd{,\body}$ in the $\di\Di$-sequence, only $\dev\D[\Hone(\body)]^3$ remains in the kernel, while $\Le(\body) \otimes \one$ is neglected.
    \begin{observation}[Dimension of {$\D [\Hone(\body)]^3 \cap [\Le(\body) \otimes \one]$}]
    Let $\D \vb{v} = \lambda \one$, then applying the Curl-operator yields 
    \begin{align}
        \Curl \D \vb{v} = \Curl (\lambda \one) = \Anti(\nabla \lambda)^{T} = 0 \, , 
    \end{align}
    implying a constant scalar field $\lambda = const$ and a linear vector field $\vb{v} \in [\Po^1(\body)]^3$.
    In fact, the vector field $\vb{v}$ is given by
    \begin{align}
        \D\vb{v} = c_0 \one \qquad \iff \qquad \vb{v} = c_0\begin{bmatrix}
            x + c_1 \\ y + c_2 \\ z + c_3
        \end{bmatrix}  \in [\Po^1(\body)]^3 \, , 
    \end{align}
    where $\{c_i\} \in \R$ are constants and $c_0 \in \R \setminus \{0\}$. 
    In other words, only one volumetric element is given by the gradients 
    \begin{align}
        \dim (\D [\Hone(\body)]^3 \cap [\Le(\body) \otimes \one]) = 1 \, ,
    \end{align}
    and is eliminated by the deviatoric operator.
    \label{ob:dim_grad}
    \end{observation}

    The next identity in the sequence follows via the classical identity $\di \curl(\cdot) = 0$ since there holds ${\di \Di \bm{D} = 0}$ if and only if $\Di \bm{D} = \curl \vb{a}$ for some vector $\vb{a}$. The algebraic identity
    \begin{align}
        \Di \Anti \vb{a} = -\curl\vb{a} \, ,
    \end{align} 
    along with the identity $\Di \Curl (\cdot) = 0$ yield
    \begin{align}
        \bm{D} = \Curl\Pm - \Anti\vb{a}  \, .
    \end{align}
    By restricting the kernel space $\ker(\di \Di)$ to its symmetric part one finds
    \begin{align}
        &\skw(\Anti\vb{a} - \Curl\Pm) = \Anti \vb{a} - \skw\Curl\Pm = 0 &&  \iff &&  \Anti \vb{a} = \skw\Curl\Pm  \, , 
    \end{align}
    leading to 
    \begin{align}
        \bm{D} = \Curl\Pm - \skw \Curl \Pm = \Curl\Pm - \dfrac{1}{2}(\Curl\Pm - \Curl\Pm^{T}) = \sym \Curl \Pm \, .
    \end{align}
    Further, we observe that 
    \begin{align}
        \range(\sym\Curl) = \range(\sym\Curl\dev) \, , \label{eq:symCurlrange}
    \end{align}
    as volumetric tensors are in $\ker(\sym\Curl)$ due to
    \begin{align}
        \sym \Curl (\lambda \one) = \sym[\Anti(\nabla \lambda)^{T}] = 0 \, ,
        \label{eq:vol_ker}
    \end{align}
    such that employing $\HsCd{,\body}$ instead of the full space does not influence the next space in the sequence.
    Finally, the last identity in the sequence is the surjection
    \begin{align}
        &\forall \, \lambda \in \Le(\body) \qquad \exists \bm{D} \in \HdD{,\body}: \qquad \di\Di \bm{D} = \lambda  \, . 
    \end{align}
    Note that 
    \begin{align}
        \di \Di \bm{D} = \di \Di (\sym\bm{D} + \skw \bm{D}) = \di \Di \sym\bm{D} \, ,
    \end{align}
    since 
    \begin{align}
        \di \Di \skw\bm{D} = \di \Di \Anti \vb{a} = -\di (\curl \vb{a}) = 0 \, .
    \end{align}
    In other words, restricting the $\HdD{,\body}$ space to symmetric tensors does not influence the range of the $\di\Di$-operator acting on the space.
Consequently, the $\di\Di$-sequence is exact on contractible domains.

\begin{figure}
	\centering
    \begin{tikzpicture}[scale = 0.6][line cap=round,line join=round,>=triangle 45,x=1.0cm,y=1.0cm]
		\clip(7,7.5) rectangle (31,10);
		\draw (9.7,9) node[anchor=north east] {$[\mathit{H}^1(\body)]^3$};
		\draw [->,line width=1.5pt] (9.7,8.5) -- (12.7,8.5);
        \draw (17.8,9) node[anchor=north east] {$\HsCd{,\body}$};
		\draw (10.2,9.5) node[anchor=north west] {$\dev \D$};
		\draw [->,line width=1.5pt] (17.9,8.5) -- (20.9,8.5);
		\draw (17.8,9.5) node[anchor=north west] {$\sym \Curl$};
		\draw (24.8,8.9) node[anchor=north east] {$\HdD{,\body}$};
		\draw [->,line width=1.5pt] (24.8,8.5) -- (27.8,8.5);
		\draw (25,9.5) node[anchor=north west] {$\di \Di$};
		\draw (27.8,9) node[anchor=north west] {$\Le(\body)$};
	\end{tikzpicture}
	\caption{The $\di\Di$ exact sequence where the last operator yields a surjection onto $\Le(\body)$. The range of each operator in the sequence is exactly the kernel of the next operator.}
	\label{fig:diDicomplex}
\end{figure} 

The exact relaxed micromorphic sequence is now given by the completion of the $\HsCd{,\body}$-space with respect to its non-deviatoric part.
Note that due to \cref{eq:symCurlrange}, only the kernel of the space is incomplete. Consequently, we complete the space by adding the space of volumetric tensors
\begin{align}
    \HsC{,\body} = \underbrace{[\dev\D[\Hone(\body)]^3 \oplus [\ker^\perp(\sym\Curl) \cap \HsC{,\body} ] ]}_{=\HsCd{,\body}} \oplus [\Le(\body) \otimes \one] \, ,
\end{align}
where $\ker^\perp(\sym\Curl)$ denotes the $\Le(\body)$-orthogonal complement of the kernel $\ker(\sym\Curl)$.
The completed sequence is depicted in \cref{fig:recomplex} and implies the identities 
\begin{subequations}
    \begin{align}
    \dev\D [\Hone(\body)]^3 \oplus [\Le(\body) \otimes \one] &= \ker(\sym\Curl) \cap \HsC{,\body} \, , \label{eq:devD} \\ \sym\Curl[\HsC{,\body }] &= \ker(\di\Di) \cap \HdD{,\body} \, , \label{eq:symCdev} \\ 
    \di \Di \HdD{,\body} &= \Le(\body) \label{eq:surj} \, ,
\end{align}
\end{subequations}
which are exact on contractible domains.

\begin{figure}
	\centering
	\begin{tikzpicture}[scale = 0.6][line cap=round,line join=round,>=triangle 45,x=1.0cm,y=1.0cm]
		\clip(7,5) rectangle (31,10);
		\draw (10.5,6.5) node[anchor=north east] {$\Le(\body) \otimes \one$};
		\draw [line width=1.5pt] (10.5,6) -- (15.8,6);
		\draw [->,line width=1.5pt] (15.8,6) -- (15.8,8);
		\draw (11.5,7) node[anchor=north west] {$\text{id}$};
		\draw (11.5,8) node[anchor=north west] {$\oplus$};

        \draw (10.5,9) node[anchor=north east] {$[\mathit{H}^1(\body)]^3$};
		\draw [->,line width=1.5pt] (10.5,8.5) -- (13.5,8.5);
        \draw (17.8,9) node[anchor=north east] {$\HsC{,\body}$};
		\draw (11,9.5) node[anchor=north west] {$\dev \D$};
		\draw [->,line width=1.5pt] (17.9,8.5) -- (20.9,8.5);
		\draw (17.8,9.5) node[anchor=north west] {$\sym \Curl$};
		\draw (24.8,8.9) node[anchor=north east] {$\HdD{,\body}$};
		\draw [->,line width=1.5pt] (24.8,8.5) -- (27.8,8.5);
		\draw (25,9.5) node[anchor=north west] {$\di \Di$};
		\draw (27.8,9) node[anchor=north west] {$\Le(\body)$};
	\end{tikzpicture}
	\caption{The relaxed micromorphic sequence as the completion of the $\di \Di$-sequence in the $\HsC{,\body}$-space via the kernel of the $\sym\Curl$-operator.}
	\label{fig:recomplex}
\end{figure} 

\section{Conforming finite elements}
In the following we present lowest order elements, linear elements, and a construction for elements of arbitrary order. The generalised construction makes use of the following polytopal definition \cite{sky_polytopal_2022,sky_higher_2023} of an $\Hone(\body)$-conforming polynomial space.
\begin{definition}[Tetrahedron $\U^p(\elem)$-polytopal spaces]
Each polytope of the tetrahedron $\elem$ is associated with a space of base functions.
\begin{itemize}
    \item each vertex $v_i$ is associated with the space of its respective base function $\ver^p_i(\elem)$. As such, there are four spaces in total $i \in \{1,2,3,4\}$ and each one is of dimension one, $\dim \ver^p_i(\elem) = 1 \quad \forall \, i \in \{1,2,3,4\}$. The base function of each respective vertex vanishes on all other vertices.
    \item for each edge $e_j$ there exists a space of edge functions $\edge^p_{j}(\elem)$ with  $j \in \mathcal{J} = \{(1,2),(1,3),(1,4),(2,3),(2,4),(3,4)\}$.
    The dimension of each edge space is given by $\dim \edge^p_j(\elem) = p-1$.
    The edge base functions of a respective edge vanish on all other edges.
    \item for each face $f_k$ there exists a space of face base functions $\face_{k}^p(\elem)$ with $k \in \mathcal{K} =  \{(1,2,3),(1,2,4),(1,3,4),(2,3,4)\}$, where the dimension of the spaces reads $\dim \face_{k}(\elem) = (p-2)(p-1)/2$.
    The base functions of a face are zero on all other faces. 
    \item lastly, the space of cell base function is given by $\cell^p(\elem)$ with the dimensionality $\dim \cell^p(\elem) = (p-3)(p-2)(p-1)/6$. Cell base functions are zero on the entire boundary of the element.
\end{itemize}
The polytope also implies the connectivity of the base functions. Vertex base functions are shared by all neighbouring elements with said vertex. Edge base functions are shared on interfacing element edges, face base functions on element faces, and cell base functions are unique to each element. 
\label{def:polyspace}
\end{definition}
The lowest order space $\U^1(T)$ is given by the barycentric coordinates
\begin{align}
    &\lambda_1(\xi,\eta,\zeta) = 1-\xi - \eta - \zeta \, , && \lambda_2(\xi,\eta,\zeta) =  \zeta \, , && \lambda_3(\xi,\eta,\zeta) = \eta \, , && \lambda_4(\xi,\eta,\zeta) = \xi \, ,
\end{align}
which we also subsequently use to construct the lowest order N\'ed\'elec elements \cite{Nedelec1980}. 
Possible definitions of the $\U^p(\elem)$-space are given for example by the Lagrange \cite{SKY2022115298,Schroder2022,Sarhil2023}, Bernstein \cite{sky_higher_2023,El-Amrani2022,ELAMRANI2023115758,AinsworthOpt} and Legendre \cite{Joachim2005,Zaglmayr2006} polynomials.

In the following, the reference tetrahedron is defined via
\begin{align}
    \Omega = \{ (\xi, \eta, \zeta) \in [0,1]^3 \; | \; \xi + \eta + \zeta \leq 1 \} \, , 
\end{align}
with the assumption that all physical elements are given by a non-degenerate mapping of it
\begin{align}
    &\vb{x}:\Omega \subset \R^3  \to T \subset \body \subset \R^3 \, , && \bm{J} = \D \vb{x} \, ,
\end{align}
where $\bm{J}$ is the associated Jacobi matrix. We note that the mapping does not have to be affine, such that curved physical elements are also viable.

\subsection{Preliminaries}\label{sec:Prelim}
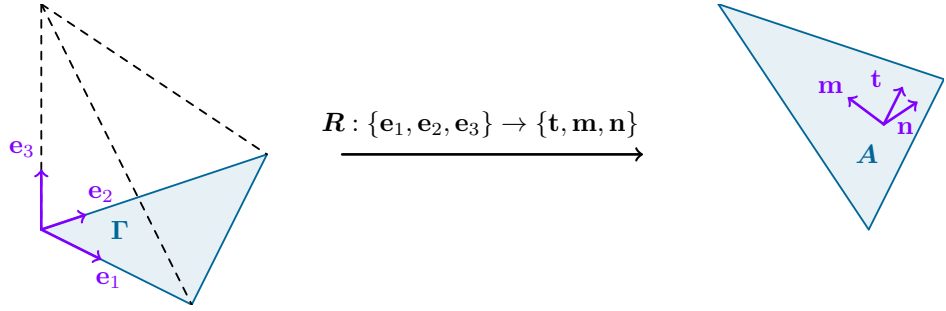
\begin{figure}
		\centering
		\definecolor{asl}{rgb}{0.4980392156862745,0.,1.}
\definecolor{zzttqq}{rgb}{0.,0.4,0.6}
\begin{tikzpicture}[line cap=round,line join=round,>=triangle 45,x=1cm,y=1cm]
\clip(-1,-1) rectangle (12,3);
\fill[line width=0.7pt,color=zzttqq,fill=zzttqq,fill opacity=0.1] (0,0) -- (3,1) -- (2,-1) -- cycle;
\fill[line width=0.7pt,color=zzttqq,fill=zzttqq,fill opacity=0.1] (9,3) -- (12,2) -- (11,0) -- cycle;
\draw [line width=0.7pt,color=zzttqq] (0,0)-- (3,1);
\draw [line width=0.7pt,color=zzttqq] (3,1)-- (2,-1);
\draw [line width=0.7pt,color=zzttqq] (2,-1)-- (0,0);
\draw [line width=0.7pt,dashed] (0,3)-- (0,0);
\draw [line width=0.7pt,dashed] (0,3)-- (2,-1);
\draw [line width=0.7pt,dashed] (0,3)-- (3,1);
\draw [line width=0.7pt,color=zzttqq] (9,3)-- (12,2);
\draw [line width=0.7pt,color=zzttqq] (12,2)-- (11,0);
\draw [line width=0.7pt,color=zzttqq] (11,0)-- (9,3);
\draw [-to,line width=1pt, color=asl] (11.2,1.4) -- (11.45,1.9);
\draw [-to,line width=1pt, color=asl] (11.2,1.4) -- (10.72,1.76);
\draw [-to,line width=1pt, color=asl] (11.2,1.4) -- (11.65,1.7);
\draw [-to,line width=1pt, color=asl] (0,0) -- (0.6,0.2);
\draw [-to,line width=1pt, color=asl] (0,0) -- (0.8,-0.4);
\draw [-to,line width=1pt, color=asl] (0,0) -- (0,0.8);
\draw [-to,line width=1pt] (4,1) -- (8,1);
\draw (0.6,-0.45) node[anchor=north west, color=asl] {$\mathbf{e}_1$};
\draw (0.5,0.7) node[anchor=north west, color=asl] {$\mathbf{e}_2$};
\draw (-0.55,1.3) node[anchor=north west, color=asl] {$\mathbf{e}_3$};
\draw (10.9,2.25) node[anchor=north west, color=asl] {$\mathbf{t}$};
\draw (10.2,2.1) node[anchor=north west, color=asl] {$\mathbf{m}$};
\draw (11.25,1.55) node[anchor=north west, color=asl] {$\mathbf{n}$};
\draw (3.6,1.7) node[anchor=north west] {$\boldsymbol{R}: \{\mathbf{e}_1, \mathbf{e}_2,  \mathbf{e}_3 \} \to \{\mathbf{t},\mathbf{m},\mathbf{n}\}$};
\draw (0.8,0) node[anchor=west, color=zzttqq] {$\boldsymbol{\Gamma}$};
\draw (10.7,1) node[anchor=west, color=zzttqq] {$\boldsymbol{\surf}$};
\end{tikzpicture}
		\caption{Mapping of the Cartesian basis from a face $\Gamma$ of the reference tetrahedron to an arbitrary face $\surf$ in 3D space via the rotation tensor $\bm{R} \in \SO(3)$.}
		\label{fig:mapr}
	\end{figure}

Conforming subspaces are built to satisfy a vanishing jump of the respective trace of the space.
As such, a discrete function, being element-wise in $\HsC{,\elem}$, belongs to a subspace of $\HsC{,\body}$ if and only if the jump of the trace vanishes 
\begin{align}
    \jump{\tr_{\HsC{}} \Pm}\at_\Xi = \jump{\sym(\Pm [\Anti\vb{n}]^{T})} \at_\Xi = 0 \, , 
\end{align}
for every arbitrarily defined interface $\Xi$ in the domain $\body$.
In order to find conforming finite elements, we reformulate the jump of the trace 
\begin{align}
    \jump{\sym(\Pm [\Anti\vb{n}]^{T})} \at_{\Xi_i} = 0 \quad \iff \quad \jump{\langle\Pm (\Anti\vb{n})^{T} ,\, \bm{S}_j \rangle}  \at_{\Xi_i} = \jump{\langle\Pm ,\, \bm{S}_j(\Anti\vb{n}) \rangle}  \at_{\Xi_i} = 0 \quad \forall \, \bm{S}_j \in \Sym(3) \, ,
\end{align}
where $\Xi_i$ now represents an element's interface.
A simple basis for $\Sym(3)$ is
\begin{align}
    \Sym(3) = \spa \{ \vb{e}_1 \otimes \vb{e}_1 , \; \vb{e}_2 \otimes \vb{e}_2 , \; \vb{e}_3 \otimes \vb{e}_3 , \; 2\sym(\vb{e}_1 \otimes \vb{e}_2) , \; 2\sym(\vb{e}_1 \otimes \vb{e}_3) , \; 2\sym(\vb{e}_2 \otimes \vb{e}_3) \} \, .
\end{align}
The relation to the surface of a finite element can now be derived by defining the rotation tensor 
\begin{align}
    &\bm{R}:\{\vb{e}_1,  \vb{e}_2, \vb{e}_3\} \to \{\vb{t}, \vb{m}, \vb{n}\} \, , && \bm{R} \in \SO(3) \, ,
\end{align}
where we consider a mapped triangular surface in 3D space with tangent unit vectors $\vb{t}$ and $\vb{m}$, and the normal unit vector $\vb{n}=\vb{t}\times\vb{m}$, see \cref{fig:mapr}. 
Consequently, the $\Sym(3)$ space can be rewritten as
\begin{align}
    \Sym(3) = \spa \{ \vb{t} \otimes \vb{t} , \; \vb{m} \otimes \vb{m} , \; \vb{n} \otimes \vb{n} , \; 2\sym(\vb{t} \otimes \vb{m}) , \; 2\sym(\vb{t} \otimes \vb{n}) , \; 2\sym(\vb{m} \otimes \vb{n}) \} \, .
\end{align}
The computation $\bm{S}_j(\Anti \vb{n})$ now yields
\begin{subequations}
   \begin{align}
    \vb{t} \otimes \vb{t} \times \vb{n} &= - \vb{t} \otimes \vb{m} \, ,  \label{eq:t1t2} \\
    \vb{m} \otimes \vb{m} \times \vb{n} &= \vb{m} \otimes \vb{t} \, , \label{eq:t2t1} \\
    \vb{n} \otimes \vb{n} \times \vb{n} &= 0 \, , \label{eq:nunu} \\
    2 \sym(\vb{t} \otimes \vb{m}) \times \vb{n} &= \vb{t} \otimes \vb{t} - \vb{m} \otimes \vb{m} \, , \label{eq:t1t1t2t2} \\
    2\sym(\vb{t} \otimes \vb{n}) \times \vb{n} &= -\vb{n} \otimes \vb{m} \, , \label{eq:nt2} \\
    2\sym(\vb{m} \otimes \vb{n}) \times \vb{n} &= \vb{n} \otimes \vb{t} \, . \label{eq:nt1}
\end{align}
\label{eq:trcond}
\end{subequations}
\begin{remark} [Normal-normal jump]
    Observe that, since \cref{eq:nunu} is always satisfied, five conditions instead of six suffice to satisfy the point-wise conformity of the function on an interface.
    \label{re:five-cond}
\end{remark}
Using \cref{eq:trcond}, we can interpret the trace conditions and relax the restrictive assumption of an orthonormal basis.
\begin{observation}[Regularity conditions] \label{ob:reg}
    Observe that \cref{eq:t1t2} and \cref{eq:t2t1} assert the continuity of the off-diagonal tangent-cotangent components of the tensor, whereas
    \cref{eq:nt1} and \cref{eq:nt2} assert the continuity of the normal-tangent components of the tensor. Due to \cref{eq:t1t1t2t2}, the continuity of each tangent-tangent component on its own is not independent. Rather, the constraint implies a jumping identity tensor $\one$ via a vanishing trace
    \begin{align}
        \tr(\vb{t}\otimes \vb{t} - \vb{m} \otimes \vb{m}) = \con{\vb{t}\otimes \vb{t} - \vb{m} \otimes \vb{m}}{\one} = 0 \, ,
    \end{align}
    such that $\vb{t}\otimes \vb{t} - \vb{m} \otimes \vb{m} \perp \one$.
\end{observation}
The conditions in \cref{eq:trcond} define five of the nine components of second-order tensors. The remaining four can be defined via
\begin{align}
    \spa \{\vb{t} \otimes \vb{n} ,\,\vb{m} \otimes \vb{n},\, \vb{n} \otimes \vb{n} , \, \one \} \, ,
    \label{eq:missing} 
\end{align}
which are allowed to jump between elements, since they are in the kernel of the trace operator $\tr_{\HsC{}}$.
\begin{observation}[Trace comparison]
    Observe that in comparison to the $\HC{}$-trace
    \begin{align}
        \tr_{\HC{}} \Pm = \Pm (\Anti \vb{n})^{T} \, ,
    \end{align}
    the $\HsC{}$-trace imposes one less constraint via \cref{eq:t1t1t2t2}, which binds the diagonal tangential-tangential components $\vb{t} \otimes \vb{t}$ and $\vb{m} \otimes \vb{m}$. In contrast, the $\HC{}$-trace constrains each of these components individually. 
\end{observation}
The main problem in defining base functions for the $\HsC{}$-space stems from \cref{eq:t1t1t2t2}. The condition demands the coupling of two tangential components of the tensor field throughout the element. Given a tetrahedralization of the domain $\cup_e \elem_e = \body$, the lowest dimensional polytopal interface is a vertex. While one can relate a vertex to a single tangential component via its neighbouring edges \cite{sky_polytopal_2022,sky_higher_2023}, it is not possible to relate it to two edges simultaneously, such that the relation is the same for all interfacing elements.  
Consequently, one is confronted with two options. Either reduce the regularity of the element further and construct a non-conforming element, or increase the regularity of the element to find a conforming finite element, but without the minimal regularity of the space. We note that the problem of constructing tetrahedral finite elements with coupled components is also related to the construction of $\HsD{,\body}$-conforming elements, where the off-diagonal components are coupled. In fact, in \cite{arnold_finite_2008} the authors prove that such a construction can only exist with full symmetric $\C^0(\body)$-continuity at the vertices. 
The latter approach to $\HsC{,\body}$-elements can be found in \cite{SkyOn}, where the authors simply require full-deviatoric $\C^0(\body)$-continuity at the vertices, while allowing the identity $\one$ to jump. As such, the regularity of the proposed elements is somewhere in between $[\Hone(\body)]^{3\times 3}$ and $\HsC{,\body}$, but the element is non-conforming in $\HC{,\body}$.   
This approach is inappropriate for the relaxed micromorphic model, where the gradient of the displacement field $\D \vb{u}$ interacts with the microdistortion field $\bm{P}$, thus demanding compatibility of the spaces $ \D \vb{u} \in \HC{,\body} \subset \HsC{,\body} \ni \Pm$ \cite{SKY2022115298,sky_hybrid_2021}. This is further emphasized by the consistent coupling condition \cite{dagostino2021consistent}, which requires consistency of the Sobolev trace spaces \cite{HIPTMAIR2023109905}.  
At this point we also emphasize that the minimal regularity may be critical in some applications in order to compute correct approximations, compare \cite{Ciarlet2014,Caorsi}. 

Unlike the case for $\HsD{,\body}$, we can fall back to the lower regularity of $\HC{,\body}$, before going all the way back to $[\Hone(\body)]^{3\times 3}$. This motivates a finite element construction based on the low order N\'ed\'elec elements \cite{Nedelec1980,Ned2}.

\subsection{Lowest order elements}
The lowest order N\'ed\'elec element $\Ned^0(T)$ \cite{Nedelec1980}, is equipped with the polynomial space 
\begin{align}
    \spa \Ned^0(T) = [\Po^0(T)]^3 \oplus \vb{x} \times [\widetilde{\Po}^0(T)]^3 = \spa \left \{ 
    \begin{bmatrix}
        1 \\ 0 \\ 0
    \end{bmatrix},
    \begin{bmatrix}
        0 \\ 1 \\ 0
    \end{bmatrix},
    \begin{bmatrix}
        0 \\ 0 \\ 1
    \end{bmatrix},
    \begin{bmatrix}
        0 \\ z \\ -y
    \end{bmatrix},
    \begin{bmatrix}
        -z \\ 0 \\ x
    \end{bmatrix},
    \begin{bmatrix}
        y \\ -x \\ 0
    \end{bmatrix}
    \right  \} \, ,
\end{align}
where $\widetilde{\Po}^p(\elem)$ is the space of homogeneous polynomials of degree $p$.
A lowest order discretisation of the matrix-valued $\HC{,\body}$-space can thus be given via $[\Ned^0(\body)]^3 \subset \HC{,\body}$, where the following polynomial space is spanned on the physical element
\begin{align}
    \spa [\Ned^0(\elem)]^3 = \spa \{\vb{e}_1, \vb{e}_2, \vb{e}_3\} \otimes ([\Po^0(\elem)]^3 \oplus \vb{x} \times [\widetilde{\Po}^0(\elem)]^3) = \R^3 \otimes ([\Po^0(T)]^3 \oplus \vb{x} \times [\widetilde{\Po}^0(\elem)]^3)  \, .
\end{align}
\begin{lemma}[Constant identity in N\'ed\'elec]
    The polynomial space of the matrix-valued N\'ed\'elec space contains the constant identity $\one \in  [\Ned^0(\elem)]^3$, but not the full linear identity $\Po^1(\elem) \otimes \one \nsubseteq [\Ned^0(\elem)]^3$.
    \label{lem:const}
\end{lemma}
\begin{proof}
    The constant identity is clearly in the first part of the N\'ed\'elec polynomial space
    \begin{align}
        \one \in \R^3 \otimes [\Po^0(\elem)]^3 = \R^{3 \times 3} \, .
    \end{align}
    The linear identity could therefore only appear in the second part of the space.
    However, the cross product in the definition ensures that $x$, $y$ and $z$ could never be in the positions $\vb{e}_1$, $\vb{e}_2$ or $\vb{e}_3$, respectively. Consequently, the following tensor fields are not in the set
    \begin{align}
        \{ x \one , y \one, z\one \} \nsubseteq  \R^3 \otimes ( \vb{x} \times [\widetilde{\Po}^0(\elem)]^3) \, , 
    \end{align}
    and as such, are not contained in the matrix-valued N\'ed\'elec element.
\end{proof}
With this in mind we define the polynomial space of the lowest order $\HsC{}$-conforming finite element 
\begin{align}
    &\Y^0(T) =  [\Ned^0(\elem)]^3 \oplus [\widetilde{\Po}^1(\elem) \otimes \one] \, , && \dim \Y^0(\elem) = 21 \, .
\end{align} 
The three additional identity fields are added as cell base functions on each element. Thus, our construction has a regularity between $\HC{,\body}$ and $\HsC{,\body}$, improving over previous constructions with higher regularity, compare with \cite{SkyOn}. The conformity of the construction is obvious, since the lowest order N\'ed\'elec basis is $\HC{}$-conforming and the identity fields vanish in the $\HsC{}$-trace.
Further, there holds
    \begin{align}
         \sym \Curl [\Y^0(\elem)] = \Po^0(\elem) \otimes \Sym(3) = \Sym(3)  \, ,
    \end{align}
since $ \curl[\spa \Ned^0(\elem)] = \R^3$ \cite{sky_polytopal_2022,Zaglmayr2006,Joachim2005}, such that
    \begin{align}
    \Curl[\Ned^0(\elem)]^3 = \Curl[\R^3 \otimes \Ned^0(\elem)] = \R^3 \otimes \curl[ \Ned^0(\elem)] = \R^{3 \times 3} \, , 
    \end{align}
    and $\sym\R^{3 \times 3} = \Sym(3)$ is a surjection.
As such, the space fits into the polynomial sequence in \cref{fig:lowest_seq}.
\begin{figure}
	\centering
	\begin{tikzpicture}[scale = 0.6][line cap=round,line join=round,>=triangle 45,x=1.0cm,y=1.0cm]
		\clip(7,5) rectangle (28,10);
		\draw (10.5,6.5) node[anchor=north east] {$\widetilde{\Po}^1(T) \otimes \one$};
		\draw [line width=1.5pt] (10.5,6) -- (14.5,6);
		\draw [->,line width=1.5pt] (14.5,6) -- (14.5,8);
		\draw (11.5,7) node[anchor=north west] {$\text{id}$};
		\draw (11.5,8) node[anchor=north west] {$\oplus$};

        \draw (10.5,9) node[anchor=north east] {$[\Po^1(T)]^3$};
		\draw [->,line width=1.5pt] (10.5,8.5) -- (13.5,8.5);
        \draw (15.5,9) node[anchor=north east] {$\Y^0(T)$};
		\draw (11.5,9.5) node[anchor=north west] {$\D$};
		\draw [->,line width=1.5pt] (17.9-2.4,8.5) -- (20.9-2.4,8.5);
		\draw (17.8-2.4,9.5) node[anchor=north west] {$\sym \Curl$};
		\draw (24.8-1.8,8.9) node[anchor=north east] {$\Po^0(T) \otimes \Sym(3)$};
		\draw [->,line width=1.5pt] (24.8-1.8,8.5) -- (27.8-1.8,8.5);
		\draw (25-1.8,9.5) node[anchor=north west] {$\di \Di$};
		\draw (27.8-1.8,8.9) node[anchor=north west] {$0$};
	\end{tikzpicture}
	\caption{A lowest order discrete polynomial sequence for approximations of the continuous relaxed micromorphic sequence. As shown in \cref{ob:dim_grad}, $\D [\Po^1(\elem)]^3$ contains only the constant identity.}
	\label{fig:lowest_seq}
\end{figure}
Now, using the N\'ed\'elec basis we directly present the base functions of our new finite element space.
\begin{definition}[Lowest order base functions]
    We give the base functions with their polytopal association.
    \begin{itemize}
        \item the base functions of the lowest $\Y^0(T)$ space on each edge $e_{ij}$ with $(i,j) \in \mathcal{J}$ are given by the N\'ed\'elec base functions
        \begin{align}
            &\bm{\rho}(\xi, \eta ,\zeta) = \vb{e} \otimes (\lambda_i \nabla_x \lambda_j - \lambda_j \nabla_x \lambda_i \, ) \, , && \vb{e} \in  \{ \vb{e}_1,\vb{e}_2,\vb{e}_3 \}  \, .
        \end{align}
        \item the cell base functions read
        \begin{align}
            \bm{\rho}(\xi, \eta ,\zeta) &= n \one \, , && n  \in \{2\xi - 1, \, 2 \eta - 1 , \, 2 \zeta - 1\}  \subseteq\Po^1(\Omega) \setminus \Po^0(\Omega) \, ,
        \end{align}
        which represent pure slopes and are orthogonal to the constant identity in the $\Le(\elem)$-sense.
    \end{itemize}
    \label{def:y0}
\end{definition}
The N\'ed\'elec base functions can also be defined on each edge $e_{ij}$ of the reference element via $\nabla_\xi$
    and mapped to the physical element using the covariant Piola transformation \cite{SKY2022115298,Zaglmayr2006} 
    \begin{align}
        \bm{\theta}_{ij} = \bm{J}^{-T} \bm{\vartheta}_{ij} = \bm{J}^{-T}(\lambda_i \nabla_\xi \lambda_j - \lambda_j \nabla_\xi \lambda_i)  \, .
        \label{eq:co}
    \end{align}
In order to compute the symmetric Curl of the base functions one can use the contravariant Piola transformation \cite{SKY2022115298} to first compute the curl of each N\'ed\'elec base function 
\begin{align}
		\mathrm{curl}_x \bm{\theta} = \dfrac{1}{\det\bm{J}} \bm{J} \mathrm{curl}_\xi \bm{\vartheta} \, , 
  \label{eq:curl_map}
\end{align}
such that the full Curl is given by $\Curl_x(\vb{e}_i \otimes \bm{\theta}) =\vb{e}_i \otimes \curl_x\bm{\theta}$, on which one applies the symmetry operator $\sym(\cdot)$. Observe that we do not need to compute the $\sym \Curl(\cdot)$ of the three identity fields, since they are in the kernel of the operator $n \one \in \ker(\sym \Curl)$.

Clearly, the space $\Y^0(\elem)$ contains at most linear functions. As such, one can further enrich the space with higher order identity fields. We define the enriched element 
\begin{align}
    \S^0(\elem) = \Y^0(\elem) \oplus [\Po^2(\elem) \setminus \Po^1(\elem)] \otimes \one \, .
\end{align}
The element is given by adding one quadratic identity field function on each edge. 
\begin{definition}[Enriched lowest order]
    The enriched lowest order element is given by the base functions of $\Y^0(\elem)$ plus additional quadratic identity fields
    \begin{align}
        \bm{\rho}(\xi, \eta ,\zeta) &= n \one \, , && n \in \edge_j^2(\elem) \, , && j \in \mathcal{J} \, ,
    \end{align}
    which are cell base functions.
    \label{def:s0}
\end{definition}
The dimension of the enriched element is $\dim \S^0(\elem) = \dim \Y^0(\elem) + 6 = 27$. The reasoning for such an enrichment is to further compensate the non-jumping constant identity field in $\Y^0(\elem)$. Note that additional base functions are discontinuous across element interfaces and can thus be statically condensated, such that the dimension of the element stiffness matrix remains unchanged.


\subsection{The linear element}
As shown in \cite{Zaglmayr2006}, it is possible to construct discrete spaces using exact polynomial sequences. In \cite{sky_higher_2023,sky_polytopal_2022,Joachim2005}, the authors construct N\'ed\'elec elements of arbitrary order by using the kernel of the previous space in the sequence to define the gradient base functions via
\begin{align}
    \R^3 \oplus \nabla [\U^{p+1}(\elem) \setminus \Po^1(\elem)] = \ker(\curl) \cap \Ned^p(\elem) = \ker(\curl) \cap \Nedtwo^p(\elem) \, ,
\end{align}
where $\U^{p+1}(\elem)$ is a discrete $\Hone(\body)$-conforming subspace on an element $\elem$. 
In practice, this approach translates to taking the gradients of all base functions aside from the vertex base functions, which underline the linear polynomial space, compare \cite{sky_polytopal_2022}.
Following the same approach for the relaxed micromorphic sequence, one could try to construct the kernel of the $\sym \Curl$-operator on the discrete $\HsC{}$ using
\begin{align}
    \dev \D [\U^{p+1}(\elem)\setminus \Po^1(\elem)]^3 \oplus [\Po^{p+1}(\elem)\otimes \one] = \ker(\sym \Curl) \cap \Y^p(\elem) \, ,
\end{align}
where $\Y^p(\body)$ is designated as the general $\HsC{,\body}$-conforming subspace, to be defined later.
However, for our linear element, using the lowest order N\'ed\'elec element basis results in linear dependence.
\begin{lemma}[Linearly dependent sum]
The construction $\Y^0(\elem) + \dev \D[\U^{2}(\elem) \setminus \Po^1(\elem)]^3$ is linearly dependent, i.e., $\dim\Y^0(\elem) + \dim \dev \D[\U^{2}(\elem) \setminus \Po^1(\elem)]^3> \dim (\Y^0(\elem) + \dev \D[\U^{2}(\elem) \setminus \Po^1(\elem)]^3)$.
    \label{th:linear_id}
\end{lemma}
\begin{proof}
    The lowest order N\'ed\'elec element of the second type can be constructed via 
    \begin{align}
        \Nedtwo^1(\elem) =  \Ned^0(\elem) \oplus \nabla[\U^2(\elem) \setminus \Po^1(\elem)] \, , 
    \end{align}
    as per \cite{Joachim2005,Zaglmayr2006}. 
    The analogous operation for matrices reads
    \begin{align}
        [\Nedtwo^1(\elem)]^3 = \R^3 \otimes \Nedtwo^1(\elem) = [\Ned^0(\elem)]^3 \oplus \D [\U^2(\elem) \setminus \Po^1(\elem)]^3 \, .
    \end{align}
    Therefore, there holds
    $\dim [\R^3 \otimes \Nedtwo^1(\elem)] = \dim (\R^3 \otimes [\Po^1(\elem)]^3) = \dim [\Po^1(\elem)]^{3\times 3} = 36$, such that the space contains the full linear polynomial space over matrices, including linear identity fields $\Po^1(\elem) \otimes \one$. 
    Now due to \cref{ob:dim_grad} we know that the set $\D [\U^2(\elem) \setminus \Po^1(\elem)]^3$ does not contain any identity field. Consequently, the dimension of $\dev \D [\U^2(\elem) \setminus \Po^1(\elem)]^3$ does not decrease. Further, the trace of the lowest order N\'ed\'elec space already contains the linear polynomial space 
    \begin{align}
        \tr [\Ned^0(\elem)]^3  = \Po^1(\elem) \, .
    \end{align}
    Consequently, one can construct the linear matrix-valued polynomial space also as
    \begin{align}
        [\Po^1(\elem)]^{3\times 3} = [\Ned^0(\elem)]^3 \oplus \dev \D [\U^2(\elem) \setminus \Po^1(\elem)]^3 \, ,
    \end{align}
    without losing any base function. This concludes the proof.
\end{proof}
With \cref{th:linear_id} at hand we use the linear N\'ed\'elec element of the second type and compensate for the jumping identity tensor using the quadratic homogeneous polynomials over second order identity tensor fields  
\begin{align}
    &\Y^1(\elem) = [\Nedtwo^1(\elem)]^3 \oplus [\widetilde{\Po}^2(\elem) \otimes \one] \, , && \dim \Y^1(\elem) = 42 \, , 
\end{align}
resulting again in a construction with regularity between $\HC{,\body}$ and $\HsC{,\body}$, which respects the sequence in \cref{fig:arb} with $p = 2$.
\begin{definition}[Linear element] \label{def:nef_enr}
    We define the base functions on their respective polytopes.
    \begin{itemize}
        \item The linear base functions of on each edge $e_{ij}$ with $(i,j) \in \mathcal{J}$ are given by the N\'ed\'elec basis
        \begin{subequations}
            \begin{align}
            &\bm{\rho}(\xi, \eta ,\zeta) = \left \{ \begin{aligned}
                &\vb{e} \otimes (\lambda_i \nabla_x \lambda_j - \lambda_j \nabla_x \lambda_i) \\
                & \vb{e} \otimes \nabla_x n  
            \end{aligned}    \right . \, , && \vb{e} \in \{ \vb{e}_1, \,\vb{e}_2, \,\vb{e}_3 \} \, , &&  n \in \edge_{ij}^2(\elem) \, ,
        \end{align} 
        \end{subequations}
        \item The cell base functions read
        \begin{align}
            \bm{\rho}(\xi, \eta ,\zeta) &= n \one \, , && n \in \edge_j^2(\elem) \, , && j \in \mathcal{J} \, ,
        \end{align}
        on each edge $e_j$ of the tetrahedron.
    \end{itemize}
    \label{def:y1}
\end{definition}
Observe that the range of the $\sym \Curl$-operator on the element is still the constant symmetric space, since the added functions are in the kernel such that 
\begin{align}
    \sym \Curl [\Y^1(\elem)] = \sym \Curl [\Po^1(\elem)]^{3 \times 3} \oplus \underbrace{\sym \Curl [\widetilde{\Po}^2(\elem) \otimes \one]}_{=0} = \Sym(3) \, , 
\end{align}
due to the surjections $\Curl [\Po^1(\elem)]^{3 \times 3} = \R^{3 \times 3}$ and $\sym \R^{3 \times 3} = \Sym(3)$. The application of the operator on the base functions can be computed using \cref{eq:curl_map}, where the symmetric Curl of the identity fields simply yields zero. 
\begin{remark}[Using the dev-operator for base functions]
    Note that using the dev-operator for the construction of the base functions via 
    \begin{align}
        &\bm{\rho} = \dev (\vb{e} \otimes \nabla_x n) = \vb{e} \otimes \nabla_x n - \dfrac{1}{3} \tr(\vb{e} \otimes \nabla_x n) \one \, , && n \in \edge_j^2(\elem) \,, && j \in \mathcal{J} \, , 
    \end{align}
    disturbs the $\HC{,\body}$-regularity of the construction, since no linear \textbf{jumping} identify fields $n \one$ with $n \in \Po^1(\elem)$ are present.
    The fields are needed in order to compensate the emerging identity terms in this base function definition.
    Consequently, such a formulation cannot achieve the higher convergence rates of $[\Nedtwo^1(\body)]^3$ for $\HC{,\body}$-fields with $\Pm \in \HC{,\body} \setminus [\Hone(\body)]^{3 \times 3}$.
    \label{re:mone}
\end{remark}

Analogously to the lowest order element, we define a further enriched linear element by adding the cubic jumping identity fields.
\begin{definition}[Enriched linear element]
    The enriched linear element is given by $\S^1(\elem) = \Y^1(\elem) \oplus [\Po^3(\elem) \setminus \Po^2(\elem)] \otimes \one$, for which the additional base functions read
    \begin{subequations}
        \begin{align}
            \bm{\rho}(\xi,\eta,\zeta) &= n\one \, , && n \in\edge_j^3(\elem) \setminus \edge_j^2(\elem) \, , && j \in \mathcal{J}\, , \\
            \bm{\rho}(\xi,\eta,\zeta) &= n\one \, , && n \in \face_k^3(\elem) \, , && k \in \mathcal{K}\, , 
    \end{align}
    which are cell base functions constructed via edge and face scalar functions.
    \end{subequations}
    \label{def:s1}
\end{definition}
We again note that all cell base functions can be statically condensated and thus do not influence the dimension of the element stiffness matrix or the global system of equations. The dimension of the space reads $\dim \S^1(\elem) = 52$. 

\begin{remark}[Enriched $\Ned^1(\elem)$]
    By the same arguments as for $\Ned^0(\elem)$ in \cref{lem:const}, the linear N\'ed\'elec element of the first type $\Ned^1(\elem)$ does not contain quadratic or cubic identities. Consequently, it can be used as basis for the enrichment instead of the linear N\'ed\'elec element of the second type $\Nedtwo^1(\elem)$. Applying the $\sym \Curl$-operator on the latter would lead to a linear symmetric space in the sequence $\Po^1(\elem) \otimes \Sym(3)$ and improved results in the $\HsC{}$-semi-norm. 
\end{remark}

\begin{remark}[An $\HC{,\body}$-non-conforming alternative] \label{re:unconforming}
    If the jumping second order identity tensor $\one$ is paramount, then a basis can be given using the linear Lagrangian space over deviatoric tensors plus a jumping linear identity
\begin{align}
    &\Dp^1(\elem) = [\U^1(\elem) \otimes \Dev(3)]\oplus [\Po^1(\elem) \otimes \one]  \, ,
\end{align}
which is the lowest order basis used in \cite{SkyOn}. 
We note that the latter basis enforces deviatoric $\C^0(\body)$-continuity, which is higher than the tangential continuity of the $\HC{,\body}$-space and N\'ed\'elec elements $[\Nedtwo(\body)]^3 \subset \HC{,\body}$. Consequently, the regularity of this construction lies between $[\Hone(\body)]^{3 \times 3}$ and $\HsC{,\body}$, seeing as
\begin{align}
    [\U^1(\body) \otimes \Dev(3)] \oplus [\Po^1(\body) \otimes \one] \subset [\Hone(\body) \otimes \Dev(3)] \oplus [\Le(\body) \otimes \one] \left \{ \begin{aligned}
        &\subset \HsC{,\body} \\
        & \nsubseteq \HC{,\body}  
    \end{aligned} \right . \, ,
\end{align}
due to $\Le(\body) \otimes \one \supset \{ \Pm \in \HC{,\body} \; | \; \dev\Pm = 0 \}$.
Therefore, this option is not consistent with the relaxed micromorphic sequence. 
\end{remark}

\subsection{An arbitrary order construction} \label{sec:arb}
In the following we develop a generalised construction for order $p \geq 2$. The element is then given by the construction itself plus the linear N\'ed\'elec element of the second type $[\Nedtwo^1(\elem)]^3$.

\begin{figure}
		\centering
		\definecolor{asl}{rgb}{0.4980392156862745,0.,1.}
\definecolor{aqaqaq}{rgb}{0.6274509803921569,0.6274509803921569,0.6274509803921569}
\definecolor{wqwqwq}{rgb}{0.3764705882352941,0.3764705882352941,0.3764705882352941}
\definecolor{zzttqq}{rgb}{0.,0.4,0.6}
\begin{tikzpicture}[line cap=round,line join=round,>=triangle 45,x=1cm,y=1cm]
\clip(2,1.5) rectangle (12,5.5);
\fill[line width=0.7pt,color=zzttqq,fill=zzttqq,fill opacity=0.1] (4,2) -- (6,3) -- (5,5) -- cycle;
\fill[line width=0.7pt,dashed,color=wqwqwq,fill=wqwqwq,fill opacity=0.1] (2,4) -- (6,3) -- (4,2) -- cycle;
\fill[line width=0.7pt,color=zzttqq,fill=zzttqq,fill opacity=0.1] (9,5) -- (8,2) -- (10,3) -- cycle;
\fill[line width=0.7pt,dashed,color=aqaqaq,fill=aqaqaq,fill opacity=0.1] (11,5) -- (10,3) -- (8,2) -- cycle;
\fill[line width=0.7pt,dashed,color=aqaqaq,fill=aqaqaq,fill opacity=0.1] (2,4) -- (5,5) -- (4,2) -- cycle;
\fill[line width=0.7pt,dashed,color=aqaqaq,fill=aqaqaq,fill opacity=0.1] (9,5) -- (11,5) -- (8,2) -- cycle;
\draw [line width=0.7pt,color=zzttqq] (4,2)-- (6,3);
\draw [line width=0.7pt,color=zzttqq] (6,3)-- (5,5);
\draw [line width=0.7pt,color=zzttqq] (5,5)-- (4,2);
\draw [line width=0.7pt,dashed,color=wqwqwq] (2,4)-- (6,3);
\draw [line width=0.7pt,dashed,color=wqwqwq] (6,3)-- (4,2);
\draw [line width=0.7pt,dashed,color=wqwqwq] (4,2)-- (2,4);
\draw [line width=0.7pt,color=zzttqq] (9,5)-- (8,2);
\draw [line width=0.7pt,color=zzttqq] (8,2)-- (10,3);
\draw [line width=0.7pt,color=zzttqq] (10,3)-- (9,5);
\draw [line width=0.7pt,dashed,color=aqaqaq] (11,5)-- (10,3);
\draw [line width=0.7pt,dashed,color=aqaqaq] (10,3)-- (8,2);
\draw [line width=0.7pt,dashed,color=aqaqaq] (8,2)-- (11,5);
\draw [line width=0.7pt,dotted] (5,5)-- (9,5);
\draw [line width=0.7pt,dotted] (8,2)-- (4,2);
\draw [line width=0.7pt,dotted] (6,3)-- (10,3);
\draw [-to,line width=1pt, color=asl] (4.4,2.2) -- (5.6,2.8);
\draw [-to,line width=1pt, color=asl] (5,2.5) -- (5,3.5);
\draw [-to,line width=1pt, color=asl] (5,2.5) -- (5.5,3.5);
\draw [-to,line width=1pt, color=asl] (8.4,2.2) -- (9.6,2.8);
\draw [-to,line width=1pt, color=asl] (9,2.5) -- (9,3.5);
\draw [-to,line width=1pt, color=asl] (9,2.5) -- (8,3.5);
\draw [line width=0.7pt,dashed] (5,3.5)-- (5.5,3.5);
\draw [line width=0.7pt,dashed] (8,3.5)-- (9,3.5);
\draw [-to,line width=1pt, color=asl] (5,2.5) -- (6,2.5);
\draw [line width=0.7pt,dashed,color=aqaqaq] (2,4)-- (5,5);
\draw [line width=0.7pt,dashed,color=aqaqaq] (5,5)-- (4,2);
\draw [line width=0.7pt,dashed,color=aqaqaq] (4,2)-- (2,4);
\draw [line width=0.7pt,dashed,color=aqaqaq] (9,5)-- (11,5);
\draw [line width=0.7pt,dashed,color=aqaqaq] (11,5)-- (8,2);
\draw [line width=0.7pt,dashed,color=aqaqaq] (8,2)-- (9,5);
\draw (5.7,3.) node[anchor=north west, color=asl] {$\mathbf{t}$};
\draw (9.6,2.9) node[anchor=north west, color=asl] {$\mathbf{t}$};
\draw (6.1,2.7) node[anchor=north west, color=asl] {$\mathbf{n}$};
\draw [-to,line width=0.7pt, color=asl] (8,2.5) -- (9,2.5);
\draw (9.020590941218504,2.5) node[anchor=north west, color=asl] {$\mathbf{n}$};
\draw (4.7,4.) node[anchor=north west, color=asl] {$\mathbf{m}$};
\draw (8.7,4.) node[anchor=north west, color=asl] {$\mathbf{m}$};
\draw (5.45,3.5) node[anchor=north west, color=asl] {$\mathbf{k}$};
\draw (7.7,3.5) node[anchor=north west, color=asl] {$\mathbf{k}_*$};
\end{tikzpicture}
		\caption{Two interfacing elements. The vectors $\vb{t}$, $\vb{m}$ and $\vb{n}$ are the same for both interfacing elements. The vectors $\vb{k}$ and $\vb{k}_*$ produce the same tangential projection on $\vb{m}$ and are orthogonal to $\vb{t}$ and the other intersecting face of the edge on each tetrahedron.}
		\label{fig:int}
	\end{figure}
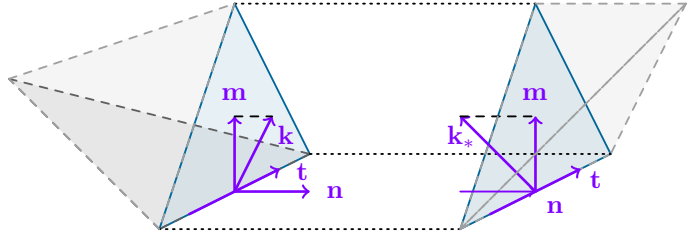

We approach the problem with the polytopal template methodology \cite{sky_polytopal_2022} using oriented polytopes \cite{SKY2022115298}. Templates on the reference element are given by $\tem$, whereas their counterparts on the physical element are denoted by $\phys$.  
We start by defining the oriented tangent vector $\bm{\tau}$ templates on the edges of the reference tetrahedron
\begin{figure}
		\centering
		\definecolor{asl}{rgb}{0.4980392156862745,0.,1.}
		\definecolor{asb}{rgb}{0.,0.4,0.6}
		\begin{tikzpicture}
			\begin{axis}
				[
				width=30cm,height=25cm,
				view={50}{15},
				enlargelimits=true,
				xmin=-1,xmax=2,
				ymin=-1,ymax=2,
				zmin=-1,zmax=2,
				domain=-10:10,
				axis equal,
				hide axis
				]
				\draw (-0.2, -0.2, -0.2) node[circle,fill=asb,inner sep=1.5pt] {};
				\draw (-0.2, -0.2, 1.2) node[circle,fill=asb,inner sep=1.5pt] {};
				\draw (1.2, -0.2, -0.2) node[circle,fill=asb,inner sep=1.5pt] {};
				\draw (-0.2, 1.2, -0.2) node[circle,fill=asb,inner sep=1.5pt] {};
				
				\addplot3[color=asb][line width=0.6pt,dotted]
				coordinates {(0,0,0)(0.5,0,0)(0,0.5,0)(0,0,0)};
				\addplot3[color=asb][line width=0.6pt,dotted]
				coordinates {(0,0,0)(0,0,0.5)};
				\addplot3[color=asb][line width=0.6pt,dotted]coordinates {(0.5,0,0)(0,0,0.5)};
				\addplot3[color=asb][line width=0.6pt,dotted]coordinates {(0,0.5,0)(0,0,0.5)};
				\fill[opacity=0.1, asb] (axis cs: 0,0,0) -- (axis cs: 0.5,0,0) -- (axis cs: 0,0.5,0) -- (axis cs: 0,0,0.5) -- cycle;
				
				\draw[color=asb] (-0.2, -0.2, -0.2) node[anchor=north east] {$_{v_{1}}$};
				\draw[color=asb] (-0.2, -0.2, 1.2) node[anchor=south east] {$_{v_{2}}$};
				\draw[color=asb] (-0.2, 1.2, -0.2) node[anchor=south west] {$_{v_{3}}$};
				\draw[color=asb] (1.2, -0.2, -0.2) node[anchor=north west] {$_{v_{4}}$};
				
				\draw[line width=.6pt, color=asb](-0.2, -0.2, 0)--(-0.2,-0.2,1);
				\draw[line width=.6pt, color=asb](0, -0.2, -0.2)--(1,-0.2,-0.2);
				\draw[line width=.6pt, color=asb](-0.2, 0, -0.2)--(-0.2,1,-0.2);
				\draw[line width=.6pt, color=asb](-0, 1.0, -0.2)--(1,0,-0.2);
				\draw[line width=.6pt, color=asb](0,-0.2,1)--(1,-0.2,0);
				\draw[line width=.6pt, color=asb](-0.2,0,1)--(-0.2,1,0);

				\draw[-to, line width=1.pt, color=asl](-0.2,-0.2,0.4)--(-0.2,-0.2,0.6);
				\draw[-to, line width=1.pt, color=asl](0.4,-0.2,-0.2)--(0.6,-0.2,-0.2);
				\draw[-to, line width=1.pt, color=asl](0.4,-0.2,0.6)--(0.6,-0.2,0.4);
				\draw[-to, line width=1.pt, color=asl](-0.2,0.4,0.6)--(-0.2,0.6,0.4);
				\draw[-to, line width=1.pt, color=asl](0.4,0.6,-0.2)--(0.6,0.4,-0.2);
				
				\draw[-to, line width=1.pt, color=asl,densely dashed](-0.2,-0.2,0.5)--(-0.4,-0.2,0.5);
				\draw[-to, line width=1.pt, color=asl,densely dashed](-0.2,-0.2,0.5)--(-0.2,-0.4,0.5);
				\draw[to-, line width=1.pt, color=asl,densely dashed](0.5,-0.2,-0.2)--(0.5,-0.4,-0.2);
				\draw[to-, line width=1.pt, color=asl,densely dashed](0.5,-0.2,-0.2)--(0.5,-0.2,-0.4);
				
				\draw[-to, line width=1.pt, color=asl,densely dashed](0.5,-0.2,0.5)--(0.5,0,0.5);
				\draw[-to, line width=1.pt, color=asl,densely dashed](0.5,-0.2,0.5)--(0.6,-0.1,0.6);
				
				\draw[-to, line width=1.pt, color=asl,densely dashed](-0.2,0.5,0.5)--(-0.1,0.6,0.6);
				\draw[to-, line width=1.pt, color=asl,densely dashed](-0.2,0.5,0.5)--(0,0.5,0.5);
				
				\draw[-to, line width=1.pt, color=asl,densely dashed](0.5,0.5,-0.2)--(0.6,0.6,-0.1);
				\draw[to-, line width=1.pt, color=asl,densely dashed](0.5,0.5,-0.2)--(0.5,0.5,0);

				\draw[-to,line width=1pt, color=asl](1, 1, 1-0.15)--(1.2, 1.2, 1.007-0.15);
				\draw[color=asl] (1.2, 1.2, 1.007-0.15) node[anchor=west] {edge template vectors};
				\draw[-to,line width=1pt, color=asl,densely dashed](1, 1, 1-0.3)--(1.2, 1.2, 1.007-0.3);
				\draw[color=asl] (1.2, 1.2, 1.007-0.3) node[anchor=west] {edge-face template vectors};

                \draw[color=asl] (-0.2, -0.2, 0.6) node[anchor=south west] {$\bm{\tau}$};
                \draw[color=asl] (-0.3, -0.3, 0.5) node[anchor=south east] {$\bm{\kappa}_2$};
                \draw[color=asl] (-0.31, -0.31, 0.5) node[anchor=north east] {$\bm{\kappa}_1$};
				
				\draw[color=asb] (-0.2, -0.2, 0.5) node[anchor=west] {$_{e_{12}}$};
				\draw[color=asb] (0.5, -0.2, -0.2) node[anchor=south] {$_{e_{14}}$};
				\draw[color=asb] (-0.2,0.5,0.5) node[anchor=east] {$_{e_{23}}$};
				\draw[color=asb] (0.5,0.5,-0.2) node[anchor=east] {$_{e_{34}}$};
				\draw[color=asb] (0.5,-0.2,0.6) node[anchor=south] {$_{e_{24}}$};
				\draw[color=asb] (0.1,0.05,0.05) node[anchor=south] {$_{f_{124}}$};
			\end{axis}
		\end{tikzpicture}
		\caption{Template vectors for the reference tetrahedron on their corresponding edges. Only vectors on the visible sides of the tetrahedron are depicted.}
		\label{fig:tet_nii}
	\end{figure}
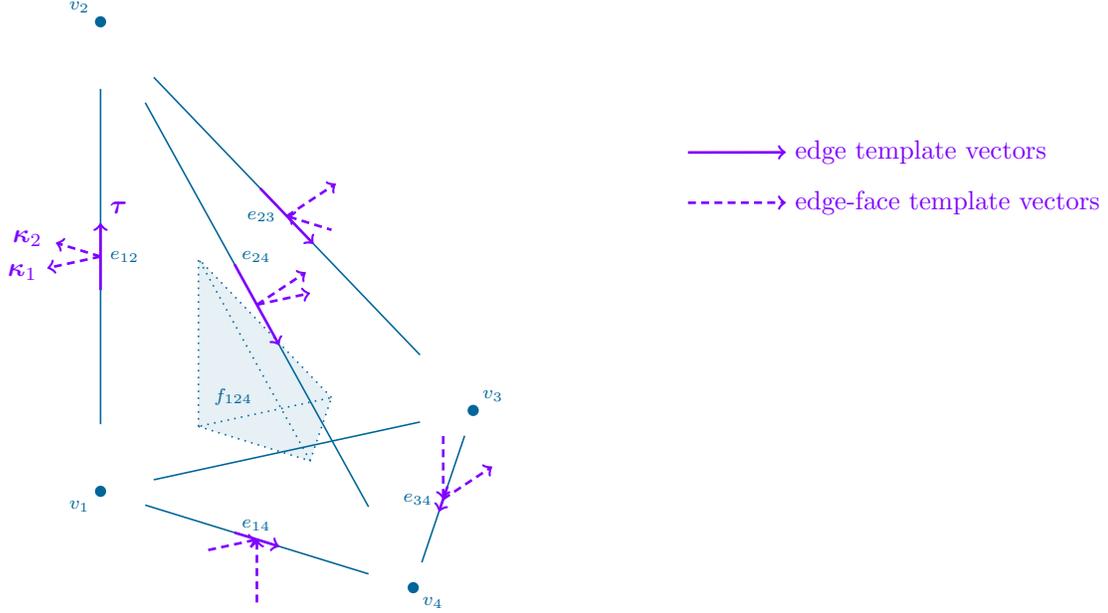
\begin{align}
    \tem_{12} &= \{ \vb{e}_3 \} \, , & \tem_{13} &= \{ \vb{e}_2 \} \, , & \tem_{14} &= \{ \vb{e}_1 \} \, , \notag \\
    	 \tem_{23} &= \{ \vb{e}_2 - \vb{e}_3  \} \, , & \tem_{24} &= \{ \vb{e}_1 - \vb{e}_3 \} \, , & 
      \tem_{34} &= \{ \vb{e}_1 - \vb{e}_2 \} \, .
      \label{eq:edge_tem}
\end{align}
The vectors are mapped to the physical element via 
\begin{align}
    &\vb{t} = \bm{J} \bm{\tau} \qquad \text{s.t.} \qquad \phys_{j} = \bm{J} (\tem_{j}) \, , && j \in \mathcal{J} \, .
\end{align}
Further, on each edge we define two oriented vectors $\bm{\kappa}_1 \perp \bm{\tau}$ and $\bm{\kappa}_2 \perp \bm{\tau}$, which are respectively orthogonal to the intersecting faces of the edge, see \cref{fig:tet_nii}.
The vectors are defined on one edge with orthogonal faces and mapped to all other edges via the covariant Piola transformation by permutations (\cref{eq:co}) of the vertex-ordering of the reference tetrahedron, compare with \cite{sky_polytopal_2022}
\begin{align}
    \tem_{12}^f &= \{ -\vb{e}_2,\, -\vb{e}_1 \} \, , & \tem_{13}^f &= \{ \vb{e}_3,\, -\vb{e}_1 \} \, , & \tem_{14}^f &= \{  \vb{e}_3,\, \vb{e}_2 \} \, , \notag \\
    \tem_{23}^f &= \{  \vb{e}_1 + \vb{e}_2 + \vb{e}_3,  \,-\vb{e}_1 \} \, , & \tem_{24}^f &= \{  \vb{e}_1 + \vb{e}_2 + \vb{e}_3,\, \vb{e}_2 \} \, , & \tem_{34}^f &= \{ \vb{e}_1 + \vb{e}_2 + \vb{e}_3, \,-\vb{e}_3 \} \, .
\end{align}
In essence, these vectors are permutations of the face normals $\bm{\nu}$. 
The covariant Piola transformation is also used to map these vectors to their physical counterparts
\begin{align}
    \vb{k}_i = \bm{J}^{-T} \bm{\kappa}_i \qquad \text{s.t.} \qquad \phys_{j}^f = \bm{J}^{-T}(\tem_{j}^f) \, , && j \in \mathcal{J} \, ,
\end{align}
thus conserving their tangential projection properties in the physical element.
On each reference face $f_{ijk}$ we define the templates of the oriented normal vectors $\bm{\nu}$
\begin{align}
    \tem_{123} &= \{ -\vb{e}_1 \} \, , &
    \tem_{124} &= \{ \vb{e}_2 \} \, , &
    \tem_{134} &= \{ -\vb{e}_3 \} \, , &
    \tem_{234} &= \{ -\vb{e}_1 - \vb{e}_2 - \vb{e}_3 \} \, .
\end{align}
The vectors are mapped to the normal vectors on the physical element via 
\begin{align}
    \vb{n} = (\cof \bm{J}) \bm{\nu} = (\det \bm{J}) \bm{J}^{-T} \bm{\nu} \qquad \text{s.t.} \qquad \phys_{k} = (\det \bm{J}) \bm{J}^{-T} (\tem_{k}) \, , && k \in \mathcal{K} \, .
\end{align}
With the normal vector at hand we also derive the cotangents on the edges of each face
\begin{align}
    &\vb{m} = \bm{Q} \vb{k} = \bm{Q} \bm{J}^{-T} \bm{\kappa} \, , && \bm{Q} = \norm{\vb{n}}^2 \one - \vb{n} \otimes \vb{n} \, , 
\end{align}
where $\bm{Q}$ is a scaled projection tensor.
The cotangent templates are therefore given by 
\begin{align}
    \phys_{j}^c = \bm{Q}(\phys_{j}^f) \, , && j \in \mathcal{J} \, .
\end{align}
The vectors are depicted in \cref{fig:int}.
Lastly, we construct two vectors orthogonal to the edge $\vb{d}_1 \perp \vb{t}$ and $\vb{d}_2 \perp \vb{t}$, which only depend on the tangent of the edge $\vb{t}$. An orthogonal vector $\vb{d}_2 \perp \vb{t}$ can be found without further information aside from the edge tangent via the algorithmic formula \cite{Ken}
\begin{align}
    \vb{d}_2 = \begin{bmatrix}
        (\sgn_* t_1) |t_3| \\
        (\sgn_* t_2) |t_3| \\
        - (\sgn_* t_3) |t_1| - (\sgn_* t_3) |t_2|
     \end{bmatrix} \perp \vb{t} \, , 
\end{align}
where the specialised signum function is defined as
\begin{align}
    \sgn_*(x) = \left \{ \begin{matrix}
        1 & \text{for} & x \geq 0 \\
        -1 & \text{for} & x < 0 
    \end{matrix} \right . \, .
\end{align}
\begin{remark}[Stability of $\vb{d}_2$]
The algorithmic approach in \cite{Ken} allows for a quasi-continuous function definition $\vb{d}_2 = \vb{d}_2(\vb{t})$, in the sense that the orthogonal vector $\vb{d}_2$ is directly dependent on the data of $\vb{t}$. Note that in general, such a function is not possible due to the hairy ball theorem \cite{Eugene}. 
The implementation of the algorithm is branch-free and square-root-function-free thanks to the use of the $copysign(\cdot)$-instruction. 
Further, there holds
\begin{align}
    \norm{\vb{t}} \leq \norm{\vb{d}_2} \leq \sqrt{2}  \norm{\vb{t}} \, ,
\end{align}
which asserts a stable evaluation of the orthogonal vector $\vb{d}_2$, as long as $\vb{t}$ does not vanish. 
This is important since an alternative algorithm based on testing the cross product with two different arbitrary vectors $\vb{a},\vb{b}$ can yield a very small orthogonal vector $\norm{\vb{a} \times \vb{t}} = \norm{\vb{a}}\norm{\vb{t}}|\sin(\phi)| \ll 1$, if the angle $\phi$ between the vectors is small, which may lead to numerical instability. The same problem can occur for a curved edge $\vb{s}$, where the normal is defined by the derivative $\vb{n} = \dd \vb{t} / \dd s$ with $\vb{t} = \vb{t}(\vb{s})$, if the change of the curve in normal direction is marginal. 
\end{remark}
The $\vb{d}_1$ vector is then given by the cross product
\begin{align}
    \vb{d}_1 = \vb{d}_2 \times \vb{t} \, .
\end{align}
The corresponding templates are given on each edge via
\begin{align}
    \phys_{j}^d &= \{ \vb{d}_1 , \, \vb{d}_2 \} \, , && j \in \mathcal{J} \, .
\end{align}
With all templates defined, we can now construct the base functions.
\begin{definition}[Base functions for $p \geq 2$]
    The construction is defined per polytope of the tetrahedron $\elem$. 
    \begin{itemize}
        \item on each edge $e_{j}$ with $j \in \mathcal{J}$ and polynomial power $p \geq 2$ we define the base functions
        \begin{align}
            &\bm{\rho}(\xi,\eta,\zeta) = n \,\vb{d} \otimes \vb{t}  \, , &&  n \in \edge_{j}^p(\elem) \, , && \vb{d} \in \phys_{j}^d \, , && \vb{t} \in \phys_{j} \, . 
        \end{align}
        Consequently, each scalar base function $n$ defines two tensorial base functions.         
        \item on each face $f_{k}$ with $k \in \mathcal{K}$ the base functions read
        \begin{subequations}
            \begin{align}
            \bm{\rho}(\xi,\eta, \zeta) &=  \left \{\begin{aligned}
                &n\, \vb{t} \otimes \vb{k}  \\
                &n \, \vb{m} \otimes \vb{k}  \\
                &n \, \vb{n} \otimes \vb{k} 
            \end{aligned} \right . \, ,
                 \qquad \begin{aligned}
                     &n \in \edge^p_{j}(\elem)  \, , \qquad \vb{t} \in \phys_{j}    \, , \qquad \vb{k} \in \{\vb{k}\in \phys_{j}^f \; | \; \vb{n} \times \vb{k} \neq 0 \} \, ,   \\
                 &\vb{m} \in \{\vb{m}\in \phys_{j}^c \; | \; \vb{n} \times \vb{m} \neq 0 \}  \, ,
                 \qquad  j \in \mathcal{J}_k \, , 
                  \qquad \vb{n} \in \phys_{k} \, ,
                 \end{aligned} 
                 \label{eq:edge-face}
                 \\ 
            \bm{\rho}(\xi,\eta, \zeta) &= \left \{\begin{aligned}
                &n\, \vb{t} \otimes \vb{m}  \\
                &n \, \vb{m} \otimes \vb{t}  \\ &n (\vb{t} \otimes \vb{t} - \vb{m} \otimes \vb{m})   \\
                &n \, \vb{n} \otimes \vb{t}  \\
                &n \, \vb{n} \otimes \vb{m} 
            \end{aligned} \right. \, , \qquad \begin{aligned}
                     &n \in \face^p_{k}(\elem)  \, , \qquad \vb{t} \in \phys_{j}   \, , \qquad \vb{n} \in \phys_{k} \, , \\
                 &\vb{m} \in \{\vb{m}\in \phys_{j}^c \; | \; \con{\vb{n}}{\vb{m}} = 0 \}  \, , \qquad j \in \min \mathcal{J}_k.
                 \end{aligned} 
                 \label{eq:facefns}
            \end{align} 
        \end{subequations}
        The first definition is for the respective face functions given by scalar functions on edges with $p \geq 2$, such that there are three face base functions for each scalar function on an edge of the face. On each face, there are three edges such that $\mathcal{J}_k \subset \mathcal{J}$ defines their index tuples. 
        In addition, each scalar base function on the face with $p \geq 3$ defines five tensorial face base functions, where $\min \mathcal{J}_k$ extracts the minimal index-tuple of an edge on the face, allowing for a consistent definition of the tangent $\vb{t}$ and cotangent $\vb{m}$ vectors on interfacing faces.
        \item on the cell we define the following base functions
        \begin{subequations}
            \begin{align}
                \bm{\rho}(\xi,\eta,\zeta) &= n \, \one \, , && n \in \bigoplus_{j \in \mathcal{J}} \edge_{j}^{p+1}(\elem)   \oplus \bigoplus_{k \in \mathcal{K}} \face^{p+1}_{k}(\elem) \oplus \cell^{p+1}(\elem)  \,,  \label{eq:ext-cell} \\
                \bm{\rho}(\xi,\eta, \zeta) &= \left \{\begin{aligned}
                &n\, \vb{t} \otimes \vb{n}  \\
                &n \, \vb{m} \otimes \vb{n}  \\ &n \, \vb{n} \otimes \vb{n}   
            \end{aligned} \right. \, ,  && \begin{aligned}
                     &n \in \face^p_{k}(\elem)  \, , \qquad \vb{t} \in \phys_{j}   \, , \qquad \vb{n} \in \phys_{k} \, , \\
                 &\vb{m} \in \{\vb{m}\in \phys_{j}^c \; | \; \vb{n} \times \vb{m} \neq 0 \}  \, , \qquad k \in \mathcal{K} \, , \qquad j \in \min \mathcal{J}_k \, ,
                 \end{aligned} \label{eq:face_basis} \\
            \bm{\rho}(\xi,\eta,\zeta) &= n \, \bm{T} \, , && n \in \cell^p(\elem)  \,, \qquad \bm{T} \in \Dev(3) \, ,
            \end{align}
        \end{subequations}
    where the first definition is of the identities via edge-cell with $p \geq 2$, face-cell with $p \geq 3$ and pure cell with $p \geq 4$, base functions. The second formula defines the face-cell base functions with $p \geq 3$, and the last definition is of the pure cell deviatoric base functions with $p \geq 4$. 
    \end{itemize}
    Recall that the definition per polytope implies connectivity analogously to \cref{def:polyspace}, such that edge base functions are shared only on interfacing edges, face base functions only on interfacing faces and the cell base functions are unique per element.
    \label{def:higher}
\end{definition}
\begin{remark}[Deviatoric tangent-tangent functions]
In \cref{eq:facefns}, the third definition is not orthogonal to identity fields since, in contrast to \cref{sec:Prelim}, the tangent and cotangent vectors do not constitute part of an orthonormal frame, such that the traces might not coincide
    \begin{align}
        \tr (\vb{t} \otimes \vb{t}) \neq \tr (\vb{m} \otimes \vb{m}) \, .
    \end{align} 
    This can be circumvented by redefining the cotangent as $\vb{m} = \norm{\vb{n}}^{-1} \vb{n} \times \vb{t}$.
    However, this does not influence the regularity of the construction.
\end{remark}
\begin{remark}[Simplified face and cell base functions]
    The face base functions in \cref{eq:edge-face} can be simplified to 
    \begin{align}
        \bm{\rho}(\xi,\eta, \zeta) &= n (\vb{e}_i \otimes \vb{k}) \, , && n \in \edge_j^p(\elem) 
        \, , && \vb{e}_i \in \{\vb{e}_1,\,\vb{e}_2,\,\vb{e}_3\} \, , && \vb{k} \in \{\vb{k}\in \phys_{j}^f \; | \; \vb{n} \times \vb{k} \neq 0 \} \, , && j \in \mathcal{J}_k \, ,
    \end{align}
    and the cell base functions in \cref{eq:face_basis} to 
    \begin{align}
        \bm{\rho}(\xi,\eta, \zeta) &= n (\vb{e}_i \otimes \vb{n}) \, , && n \in \face^p_{k}(\elem) 
        \, , && \vb{e}_i \in \{\vb{e}_1,\,\vb{e}_2,\,\vb{e}_3\} \, , && \vb{n} \in \phys_{k} \, , && k \in \mathcal{K} \, ,
    \end{align}
    using the Cartesian basis, since $\R^3 = \spa \{\vb{e}_1,\,\vb{e}_2,\,\vb{e}_3\} = \spa \{\vb{t} , \,\vb{m} , \,\vb{n}\}$. 
\end{remark}
A full finite element of order $p \geq 2$ can now be given using the construction plus $[\Nedtwo^1(T)]^3$. The symmetric Curl of the extended construction is given by
\begin{align}
    &\sym \Curl (n \bm{T}) =  \sym [\bm{T} \Anti(\bm{J}^{-T} \nabla_\xi n)^{T}]  \, , && n \in [\U^p(\Omega) \setminus \U^1(\Omega)] \, , && \bm{T} \in \phys(\elem) \, , 
\end{align}
for affinely mapped elements, where $\phys(T)$ represents the collection of all mapped second order template tensors. 
If the geometry is non-affine, then the full chain-rule is required
\begin{align}
    &\sym \Curl (n \bm{T}) =  \sym [(n \bm{T})_{,i} \Anti(\bm{J}^{-T} \vb{e}_i)^{T}]  \, , && n \in [\U^p(\Omega) \setminus \U^1(\Omega)] \, , && \bm{T} \in \phys(\elem) \, .
\end{align}
\begin{theorem}[Linear independence]
    The construction presented in \cref{def:higher} in combination with $[\Nedtwo^1(T)]^3$ is linearly independent. 
\end{theorem}
\begin{proof}
    The linear independence of $[\Nedtwo^1(\elem)]^3$ is obvious, and it spans $[\Po^1(\elem)]^{3 \times 3}$. Its linear independence from the extended construction is clear, as the construction relies on polynomials from $\Po^p(T) \setminus \Po^1(T)$. In order to prove the linear independence of the extended construction itself observe that on each edge exactly nine tensorial base functions are defined for each scalar function $n \in \edge^p(\elem)$. The tensors on each edge are given by
    \begin{align}
        \{ \vb{d}_1\otimes \vb{t}, \, \vb{d}_2\otimes \vb{t} , \, \vb{t}\otimes \vb{k}_1  , \, \vb{t}\otimes \vb{k}_2 , \,  \vb{m}_1\otimes \vb{k}_1 , \, \vb{m}_2\otimes \vb{k}_2 , \, \vb{n}_1\otimes \vb{k}_1 , \, \vb{n}_2\otimes \vb{k}_2 , \, \one\} \, .
    \end{align}
    It suffices to show that the latter represents a basis for $\R^{3 \times 3}$ in order to prove linear independence on the edges. Clearly, the pairs $\vb{d}_i$, $\vb{m}_i$, $\vb{k}_i$ and $\vb{n}_i$ are respectively linearly independent and orthogonal to $\vb{t}$, such that each pair spans the plane $\{\vb{t}\}^{\perp}=\{ \vb{v}\in\R^3 \; | \; \con{\vb{t}}{\vb{v}}=0\}$. Assuming a perpendicular intersection of two faces on an edge results in $\vb{m}_i \parallel \vb{k}_i$ and $\vb{k}_i \perp \vb{n}_i$. The $\vb{d}_i$ vectors can be chosen independently, such that we set $\vb{d}_i = \vb{k}_i$.
    Consequently, the set transforms to 
    \begin{align}
        \{ \vb{d}_1\otimes \vb{t}, \, \vb{d}_2\otimes \vb{t} , \, \vb{t}\otimes \vb{d}_1  , \, \vb{t}\otimes \vb{d}_2 , \,  \vb{d}_1\otimes \vb{d}_1 , \, \vb{d}_2\otimes \vb{d}_2 , \, \vb{d}_2\otimes \vb{d}_1 , \, \vb{d}_1\otimes \vb{d}_2 , \, \one\} \, \, , 
    \end{align}
    with some $\pm$ signs and scaling before each tensor, which are omitted for simplification of the presentation.
    Due to $\vb{t} \perp \vb{d}_1 \perp \vb{d}_2$ this is clearly a linearly independent set. The same holds true also for the general case and can be observed by the $\HC{,\elem}$- and $\HsC{,\elem}$-traces. 
    Observe that aside from the identity $\one$, the only tensors that do not vanish in the $\HC{,\elem}$-trace on the edge are defined with the edge tangent vector in the second position
    \begin{align}
        \tr_{\HC{}} (\vb{d}_i \otimes \vb{t}) \at_{\curv} =  (\vb{d}_i \otimes \vb{t}) \vb{t} \at_{\curv} \neq 0 \, ,
    \end{align}
    and are thus clearly independent of the other tensors. Their linear independence of each other is obvious due to $\vb{d}_1 \perp \vb{d}_2$.
    The edge-cell identity base function is linearly independent of the others since it vanishes in the $\HsC{,\elem}$-trace $n\one \in \ker(\tr_{\HsC{}})$.
    Lastly, the trace of edge-face base functions of one face vanish on the other intersecting face
    \begin{align}
        \tr_{\HC{}} (\vb{v} \otimes \vb{k}_1) \at_{\surf_2} =    (\vb{v}\otimes \vb{k}_1 )  (\Anti \vb{n}_2)^{T}   \at_{\surf_2} = 0  \qquad \forall \, \vb{v} \in \R^3 \, ,
    \end{align}
    due to $\vb{k}_1 \parallel \vb{n}_2$. By the same argument the functions do not vanish on their respective face due to $\vb{k}_1 \nparallel \vb{n}_1$.  
    Analogously, we observe that each face-scalar base function $n \in \face^p(\elem)$ is multiplied with the set
    \begin{align}
        \{\vb{t}\otimes \vb{m}, \, \vb{m}\otimes \vb{t}, \, \vb{t}\otimes \vb{t}-  \vb{m}\otimes \vb{m}, \, \vb{n}\otimes \vb{t}, \, \vb{n}\otimes \vb{m}, \, \vb{t}\otimes \vb{n}, \, \vb{m}\otimes \vb{n}, \, \vb{n}\otimes \vb{n}, \, \one \} \, ,
    \end{align}
    which is clearly a basis for $\R^{3\times 3}$.
    Lastly, each scalar cell base function $n \in \cell^p(\elem)$ is multiplied with the set $[\Dev(3) \oplus \{\one\}] = \R^{3 \times 3}$. This concludes the proof.
\end{proof}
\begin{theorem}[$\HsC{,\body}$-conformity] \label{th:conform}
    The construction from \cref{def:higher} is conforming in the space $\HsC{,\body}$.
\end{theorem}
\begin{proof}
    We start with the edge base functions and observe that on an edge 
    \begin{align}
        &\tr_{\HC{}} ( n \, \vb{d}_i \otimes \vb{t} ) \at_\curv = ( n \, \vb{d}_i \otimes \vb{t} ) \vb{t} \at_\curv = n \| \vb{t} \|^2 \vb{d}_i \at_\curv \, , && n \in \edge^p(\elem) \, ,
    \end{align}
    such that the component is the same for all interfacing elements and respects $\HC{,\body}$-conformity. We move on to the edge-face base functions. There holds, $\vb{k}$ is orthogonal to $\vb{t}$ and the plane of the other intersecting face of the element on the edge. Further, by design, its tangential projection on the cotangent $\vb{m}$ of the edge on the face is the same for two interfacing elements. Consequently, we have 
    \begin{align}
        &\jump{\tr_{\HsC{}} n\,\vb{v} \otimes \vb{k} } \at_\Xi = 0 \qquad \forall \, \vb{v} \in \{\vb{t}, \, \vb{m}, \,\vb{n}\} \, , && n \in \edge^p(\elem) \, ,
    \end{align}
    for the interface $\Xi$ between the two elements. The conformity follows via $\vb{k} \times \vb{n} = \vb{k}_* \times \vb{n}$, and by observing that $\vb{t}$, $\vb{m}$ and $\vb{n}$ are shared on interfacing faces.   
    The relations are visualised in \cref{fig:int}. Finally, due to the underlying scalar edge base function being $n \in \edge_j^p(\elem)$, the tensorial function vanishes also on all non-intersecting faces.
    The conformity of the pure-face and face-cell base functions is obvious as they are built directly according to the $\HsC{,\body}$-trace conditions, and vanish on all other faces due to the underlying scalar functions $n \in \face^p_k(\elem)$. Finally, the cell base functions do not influence the conformity, since their $\HsC{,\elem}$-trace vanishes on the entire boundary of the element. This is clear since they are either built as identity fields $n \one$, which vanish in the trace or as $\bm{\rho} \in \cell^p(\elem) \otimes \Dev(3)$, where $n \in \cell^p(\elem)$ vanish on the entire boundary of the element.  
\end{proof}
The main ingredient in the construction of this novel finite element are the face base functions in \cref{eq:edge-face} and cell base functions in \cref{eq:ext-cell}, whose significance we now clarify. For simplicity, assume the edge is given by the intersection of two orthogonal faces, see \cref{fig:intsect}. 
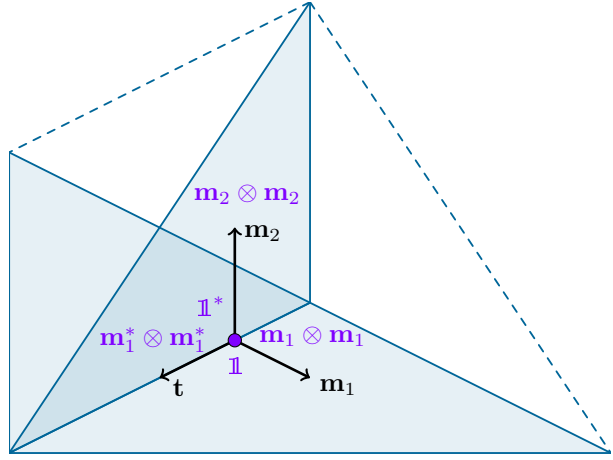
\begin{figure}
		\centering
		\definecolor{asl}{rgb}{0.4980392156862745,0.,1.}
\definecolor{xdxdff}{rgb}{0.4980392156862745,0.,1.}
\definecolor{zzttqq}{rgb}{0.,0.4,0.6}
\begin{tikzpicture}[line cap=round,line join=round,>=triangle 45,x=1cm,y=1cm]
\clip(6,5.5) rectangle (14,12);
\fill[line width=0.7pt,color=zzttqq,fill=zzttqq,fill opacity=0.1] (6,6) -- (10,12) -- (10,8) -- cycle;
\fill[line width=0.7pt,color=zzttqq,fill=zzttqq,fill opacity=0.1] (6,6) -- (14,6) -- (10,8) -- cycle;
\fill[line width=0.7pt,color=zzttqq,fill=zzttqq,fill opacity=0.1] (6,6) -- (6,10) -- (10,8) -- cycle;
\draw [line width=0.7pt,color=zzttqq] (6,6)-- (10,12);
\draw [line width=0.7pt,color=zzttqq] (10,12)-- (10,8);
\draw [line width=0.7pt,color=zzttqq] (10,8)-- (6,6);
\draw [line width=0.7pt,color=zzttqq] (6,6)-- (14,6);
\draw [line width=0.7pt,color=zzttqq] (14,6)-- (10,8);
\draw [line width=0.7pt,color=zzttqq] (10,8)-- (6,6);
\draw [line width=0.7pt,color=zzttqq] (6,6)-- (6,10);
\draw [line width=0.7pt,color=zzttqq] (10,8)-- (6,10);
\draw [line width=0.7pt,dashed,color=zzttqq] (10,12)-- (6,10);
\draw [line width=0.7pt,dashed,color=zzttqq] (10,12)-- (14,6);
\draw [-to,line width=1pt,color=black] (9,7.5) -- (10,7);
\draw [-to,line width=1pt,color=black] (9,7.5) -- (9,9);
\draw [-to,line width=1pt,color=black] (9,7.5) -- (8,7);
\draw (10.,7.1) node[anchor=north west,color=black] {$\mathbf{m}_1$};
\draw (9.2,7.8) node[anchor=north west,color=asl] {$\mathbf{m}_1 \otimes \mathbf{m}_1$};
\draw (8.75,7.8) node[anchor=north east,color=asl] {$\mathbf{m}_1^* \otimes \mathbf{m}_1^*$};
\draw (10.,9.7) node[anchor=north east,color=asl] {$\mathbf{m}_2 \otimes \mathbf{m}_2$};
\draw (9.,7.4) node[anchor=north,color=asl] {$\one$};
\draw (9.,7.7) node[anchor=south east,color=asl] {$\one^*$};
\draw (9.,9.15) node[anchor=north west,color=black] {$\mathbf{m}_2$};
\draw (8.05,7.1) node[anchor=north west,color=black] {$\mathbf{t}$};
\begin{scriptsize}
\draw [fill=xdxdff] (9,7.5) circle (2.5pt);
\end{scriptsize}
\end{tikzpicture}
		\caption{Edge on the perpendicular intersection of two faces, such that tangent $\vb{t}$ and cotangent $\vb{m}$ vectors form an orthonormal basis. On the node, the elements share the face degree of freedom for $\vb{m}_2 \otimes \vb{m}_2$. Further, each element defines its own face degrees of freedom for the component $\vb{m}_1 \otimes \vb{m}_1$, which may be constrained by further interfacing elements underneath. Finally, each element defines its own cell degree of freedom for the identity component $\one$, which implicitly determines $\vb{t} \otimes \vb{t}$. The un-shared components on the left element are marked with a star $*$ superscript.}
		\label{fig:intsect}
	\end{figure}
In that case, the tensorial basis of the base functions reads
\begin{align}
    \{ \vb{m}_1 \otimes \vb{m}_1 , \, \vb{m}_2 \otimes \vb{m}_2 , \, \one \} \, ,
\end{align}
and all other tensorial base functions on the edge are off-diagonal. Further, the cotangents are orthogonal to each other $\vb{m}_1 \perp \vb{m}_2$ due to the orthogonality of the intersecting faces. Now, let the tangent and cotangent vectors be unit vectors, then the identity $\one$ is simply 
\begin{align}
    \one = \vb{t} \otimes \vb{t} + \vb{m}_1 \otimes \vb{m}_1  + \vb{m}_2 \otimes \vb{m}_2  \, ,
\end{align}
which couples the tangent vectors as per \cref{ob:reg}.
Clearly, the identity, being a cell base function, controls the tangent-tangent component $\vb{t} \otimes \vb{t}$ on edges. This property allows the element to weakly capture also the regularity of $\HC{,\body}$, in which the tangent-tangent component is continuous on edges.
Lastly, observe that this characteristic is maintained also for non-perpendicular faces, as all other components aside from the identity are orthogonal to the tangent-tangent component $\vb{t} \otimes \vb{t}$, and are controlled by their own base functions.  

We designate the arbitrary order extension as $\widetilde{\Y}^p(\elem)$, such that the element construction reads
\begin{align}
        \Y^p(\elem) &= [\Nedtwo^1(\elem)]^3 \oplus \widetilde{\Y}^p(\elem) \, , && 
        p \geq 2 \, . 
\end{align}
Since the element is defined by the full polynomial space plus additional identities in the kernel, there holds
\begin{align}
    \dim \Y^p(\elem) = \dim [\Po^p(\elem)]^{3 \times 3} + \dim \widetilde{\Po}^{p+1}(\elem) = \dfrac{(p + 2)(p + 3)(3p + 4)}{2} \, . 
\end{align}
Due to $\widetilde{\Po}^{p+1}(\elem) \otimes \one \in \ker(\sym \Curl)$, and the surjections $\curl [\Po^p(\elem)]^3 = [\Po^{p-1}(\elem)]^3 \, \cap \,  \range(\curl)$ and $\sym \R^{3 \times 3} = \Sym(3)$, the $\Y^p(\body)$-element fits into the polynomial sequence depicted in \cref{fig:arb}.
\begin{figure}
	\centering
\begin{tikzpicture}[scale = 0.6][line cap=round,line join=round,>=triangle 45,x=1.0cm,y=1.0cm]
		\clip(4,5) rectangle (30,10);
		\draw (10.5,6.5) node[anchor=north east] {$[\Po^{p+1}(T) \setminus \Po^{0}(T)]  \otimes \one$};
		\draw [line width=1.5pt] (10.5,6) -- (14.5,6);
		\draw [->,line width=1.5pt] (14.5,6) -- (14.5,8);
		\draw (11.5,7) node[anchor=north west] {$\text{id}$};
		\draw (11.5,8) node[anchor=north west] {$\oplus$};

        \draw (10.5,9) node[anchor=north east] {$[\Po^{p+1}(T)]^3$};
		\draw [->,line width=1.5pt] (10.5,8.5) -- (13.5,8.5);
        \draw (15.5,9) node[anchor=north east] {$\Y^p(T)$};
		\draw (11.5,9.5) node[anchor=north west] {$\D$};
		\draw [->,line width=1.5pt] (17.9-2.4,8.5) -- (20.9-2.4,8.5);
		\draw (17.8-2.4,9.5) node[anchor=north west] {$\sym \Curl$};
		\draw (24.8-1.2,9) node[anchor=north east] {$\Po^{p-1}(T) \otimes \Sym(3)$};
		\draw [->,line width=1.5pt] (24.8-1.2,8.5) -- (27.8-1.2,8.5);
		\draw (25-1.2,9.5) node[anchor=north west] {$\di \Di$};
		\draw (27.8-1.2,8.9) node[anchor=north west] {$\Po^{p-2}(\elem)$};
	\end{tikzpicture}
	\caption{The general order discrete polynomial sequence for discretisations of the continuous relaxed micromorphic sequence. The constant identity is contained in $\D[\Po^{p+1}(\elem)]^3$, such that the deviatoric operator is not used.}
	\label{fig:arb}
\end{figure}
We note that the construction could be also be reduced to capture only full polynomial spaces $[\Po^p(\elem)]^{3 \times 3}$ by refining the identity cell base functions as
\begin{align}
    &\bm{\rho} = n \one \, , && n \in \bigoplus_{j \in \mathcal{J}} \edge_{j}^{p}(\elem)   \oplus \bigoplus_{k \in \mathcal{K}} \face^{p}_{k}(\elem) \oplus \cell^{p}(\elem)  \, ,
\end{align}
which would lead to the sequence in \cref{fig:arb2} with the element $\M^p(\elem) = \Y^p(\elem) \setminus [\widetilde{\Po}^{p+1}(\elem) \otimes \one]$, and the dimensionality $9(p+3)(p+2)(p+1)/6$ with $p \geq 2$.
\begin{figure}
	\centering
\begin{tikzpicture}[scale = 0.6][line cap=round,line join=round,>=triangle 45,x=1.0cm,y=1.0cm]
		\clip(5,5) rectangle (30,10);
		\draw (10.5,6.5) node[anchor=north east] {$[\Po^{p}(T) \setminus \Po^{0}(T)]  \otimes \one$};
		\draw [line width=1.5pt] (10.5,6) -- (14.5,6);
		\draw [->,line width=1.5pt] (14.5,6) -- (14.5,8);
		\draw (11.5,7) node[anchor=north west] {$\text{id}$};
		\draw (11.5,8) node[anchor=north west] {$\oplus$};

        \draw (10.5,9) node[anchor=north east] {$[\Po^{p+1}(T)]^3$};
		\draw [->,line width=1.5pt] (10.5,8.5) -- (13.5,8.5);
        \draw (15.6,9) node[anchor=north east] {$\M^p(T)$};
		\draw (11.5,9.5) node[anchor=north west] {$\D$};
		\draw [->,line width=1.5pt] (17.9-2.4,8.5) -- (20.9-2.4,8.5);
		\draw (17.8-2.4,9.5) node[anchor=north west] {$\sym \Curl$};
		\draw (24.8-1.2,9) node[anchor=north east] {$\Po^{p-1}(T) \otimes \Sym(3)$};
		\draw [->,line width=1.5pt] (24.8-1.2,8.5) -- (27.8-1.2,8.5);
		\draw (25-1.2,9.5) node[anchor=north west] {$\di \Di$};
		\draw (27.8-1.2,8.9) node[anchor=north west] {$\Po^{p-2}(\elem)$};
	\end{tikzpicture}
	\caption{A second general order discrete polynomial sequence for discretisations of the continuous relaxed micromorphic sequence with a reduced kernel yielding complete polynomial spaces.}
	\label{fig:arb2}
\end{figure}

\section{Numerical examples}\label{sec:num.examples}
In this section we demonstrate the behaviour of the finite elements with numerical examples. The following elements are tested (see \cref{def:y0,def:s0,def:y1,def:s1,def:higher})
\begin{align}
    &\Y^0(\body) \, , && \S^0(\body) \, , && \Y^1(\body) \, , && \S^1(\body) \, , && \Y^2(\body) \, , && \M^2(\body) \, ,  
\end{align}
where $\M^2(\body) = \Y^2(\body) \setminus [\widetilde{\Po}^3(\body)] \otimes \one$ drops cubic identity fields and thus is built using complete polynomial spaces as per \cref{fig:arb2}.

The first example verifies the regularity of the novel elements. In order to frame the results with respect to other discrete spaces, we compare with the  $\Hone(\body)$-conforming Lagrangian space $[\Lag^1(\body)]^{3\times 3}$, its deviatoric version with a jumping identity $\Dp^1(\body)$ \cite{SkyOn} as per \cref{re:unconforming}, and the two types of N\'ed\'elec elements, $[\Ned(\body)]^3$ and $[\Nedtwo(\body)]^3$ \cite{Nedelec1980,Ned2}.
We benchmark the elements for $\Hone(\body)$-, $\HC{,\body}$- and $\HsC{,\body}$-regularity. Specifically for $\HsC{,\body}$ we also examine the case of a jumping constant identity. Relative errors are measured in the $\Le(\body)$-norm
$
     \norm{\widetilde{\bm{P}} - \bm{P}^h}_{\Le} / \norm{\widetilde{\bm{P}}}_{\Le} \, , 
$
where $\widetilde{\bm{P}}$ represents the exact solution and $\bm{P}^h$ is the obtained discrete solution.

The second example demonstrates the significance of the new discrete spaces in the relaxed micromorphic model using a solution where the microdistortion field $\Pm$ is $\HsC{,\body}$-conforming but not $\HC{,\body}$-conforming.

\subsection{Regularity benchmark}
\begin{figure}
		\centering
		\definecolor{asl}{rgb}{0.4980392156862745,0.,1.}
\definecolor{zzttqq}{rgb}{0.,0.4,0.6}
\begin{tikzpicture}[line cap=round,line join=round,>=triangle 45,x=1cm,y=1cm,scale=0.8]
\clip(1.5,4.5) rectangle (20,14);
\fill[line width=0.7pt,color=zzttqq,fill=zzttqq,fill opacity=0.10000000149011612] (3,8) -- (15,5) -- (15,9) -- (3,12) -- cycle;
\fill[line width=0.7pt,color=zzttqq,fill=zzttqq,fill opacity=0.10000000149011612] (15,9) -- (18,10) -- (18,6) -- (15,5) -- cycle;
\fill[line width=0.7pt,color=zzttqq,fill=zzttqq,fill opacity=0.10000000149011612] (18,10) -- (6,13) -- (3,12) -- (15,9) -- cycle;
\draw [line width=0.7pt,color=zzttqq] (3,8)-- (15,5);
\draw [line width=0.7pt,color=zzttqq] (15,5)-- (15,9);
\draw [line width=0.7pt,color=zzttqq] (15,9)-- (3,12);
\draw [line width=0.7pt,color=zzttqq] (3,12)-- (3,8);
\draw [line width=0.7pt,color=zzttqq] (15,9)-- (18,10);
\draw [line width=0.7pt,color=zzttqq] (18,10)-- (18,6);
\draw [line width=0.7pt,color=zzttqq] (18,6)-- (15,5);
\draw [line width=0.7pt,color=zzttqq] (15,5)-- (15,9);
\draw [line width=0.7pt,color=zzttqq] (18,10)-- (6,13);
\draw [line width=0.7pt,color=zzttqq] (6,13)-- (3,12);
\draw [line width=0.7pt,color=zzttqq] (3,12)-- (15,9);
\draw [line width=0.7pt,color=zzttqq] (15,9)-- (18,10);
\draw [line width=0.7pt,dashed] (6,9)-- (18,6);
\draw [line width=0.7pt,dashed] (6,9)-- (3,8);
\draw [line width=0.7pt,dashed] (6,13)-- (6,9);
\draw [line width=0.7pt,dashed] (7,7)-- (7,11);
\draw [line width=0.7pt,dashed] (7,11)-- (10,12);
\draw [line width=0.7pt,dashed] (10,8)-- (7,7);
\draw [line width=0.7pt,dashed] (10,8)-- (10,12);
\draw [line width=0.7pt,dashed] (11,6)-- (14,7);
\draw [line width=0.7pt,dashed] (11,6)-- (11,10);
\draw [line width=0.7pt,dashed] (14,7)-- (14,11);
\draw [line width=0.7pt,dashed] (11,10)-- (14,11);
\draw [-to,line width=1pt,color=asl] (8.5,9.5) -- (9.3,9.3);
\draw [-to,line width=1pt,color=asl] (12.5,8.5) -- (13.3,8.3);
\draw [line width=0.7pt,dashed] (7,7)-- (10,12);
\draw [line width=0.7pt,dashed] (7,11)-- (10,8);
\draw [line width=0.7pt,dashed] (11,6)-- (14,11);
\draw [line width=0.7pt,dashed] (11,10)-- (14,7);
\draw [line width=0.7pt,dashed] (12,9.5)-- (15.428571428571429,10.642857142857142);
\draw [-to,line width=1pt] (15.428571428571429,10.642857142857142) -- (16.5,11);
\draw [-to,line width=1pt] (10.5,11) -- (10.5,13.5);
\draw [-to,line width=1pt] (16.5,7.5) -- (18.5,7);
\draw (18.65,7.15) node[anchor=north west] {$x$};
\draw (16.6,11.4) node[anchor=north west] {$y$};
\draw (10.25,14.1) node[anchor=north west] {$z$};
\draw (13.4,8.45) node[anchor=north west,color=asl] {$\mathbf{n}$};
\draw (9.35,9.5) node[anchor=north west,color=asl] {$\mathbf{n}$};
\draw (1.3,7.9) node[anchor=north west] {$_{(-3,-1,-1)}$};
\draw (5.4,6.95) node[anchor=north west] {$_{(-1,-1,-1)}$};
\draw (9.3,5.95) node[anchor=north west] {$_{(1,-1,-1)}$};
\draw (13.2,5.05) node[anchor=north west] {$_{(3,-1,-1)}$};
\begin{scriptsize}
\draw [fill=black] (3,8) circle (2.5pt);
\draw [fill=black] (7,7) circle (2.5pt);
\draw [fill=black] (11,6) circle (2.5pt);
\draw [fill=black] (15,5) circle (2.5pt);
\draw [color=black] (8.5,9.5)-- ++(-2.5pt,-2.5pt) -- ++(5pt,5pt) ++(-5pt,0) -- ++(5pt,-5pt);
\draw [color=black] (12.5,8.5)-- ++(-2.5pt,-2.5pt) -- ++(5pt,5pt) ++(-5pt,0) -- ++(5pt,-5pt);
\draw [color=black] (16.5,7.5)-- ++(-2.5pt,-2.5pt) -- ++(5pt,5pt) ++(-5pt,0) -- ++(5pt,-5pt);
\draw [color=black] (10.5,11)-- ++(-2.5pt,-2.5pt) -- ++(5pt,5pt) ++(-5pt,0) -- ++(5pt,-5pt);
\draw [color=black] (12,9.5)-- ++(-2.5pt,-2.5pt) -- ++(5pt,5pt) ++(-5pt,0) -- ++(5pt,-5pt);
\end{scriptsize}
\end{tikzpicture}
		\caption{The domain $\overline{\body} = [-3,3] \times [-1,1]^2$. At $x =-1$ and $x=1$, the vector $\vb{n}$ represents the corresponding normal vector on the plane, respectively.}
		\label{fig:exdom}
	\end{figure}
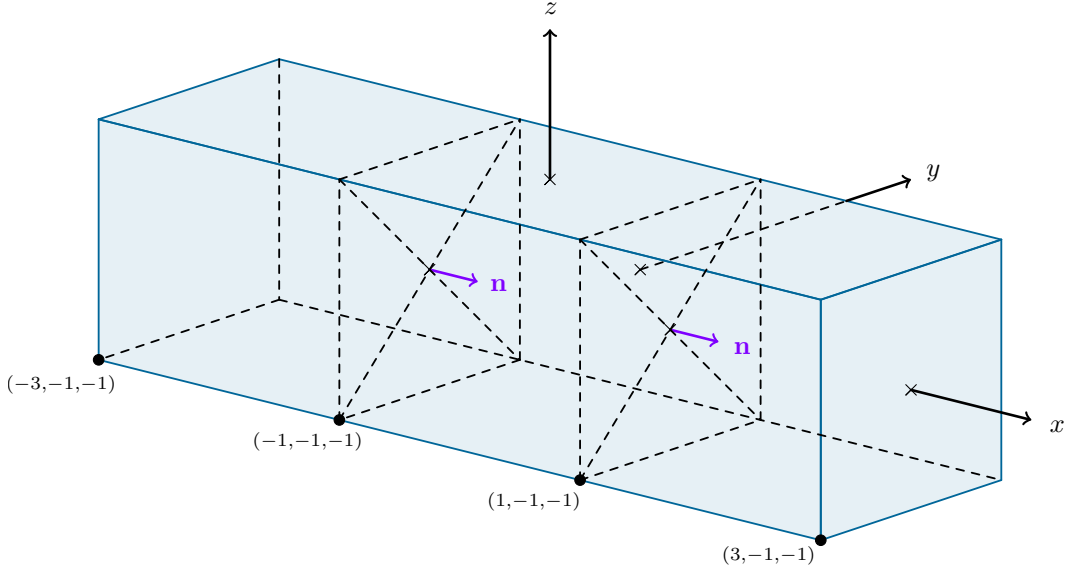

In order to demonstrate the regularity of the novel elements we consider four benchmarks with decreasing regularity. The domain $\overline{\body} = [-3,3] \times [-1,1]^2$ remains the same for all benchmarks and is depicted in \cref{fig:exdom}.
The problem is the quasi $\HsC{,\body}$-$\Le(\body)$-projection
\begin{align}
    &\int_\body \con{\delta \Pm}{\Pm} + \con{\sym\Curl \delta \Pm}{\sym \Curl\Pm} \, \dd \body  = \int_\body \con{\delta \Pm}{\widetilde{\Pm}} \, \dd \body  \qquad \forall \, \delta \Pm \in \HsC{,\body} \, .
\end{align}
Clearly, the problem is well-posed as the bilinear form represents the $\HsC{,\body}$-scalar product. Now, assuming the analytical solution is in the kernel $\widetilde{\Pm} \in \HsC{,\body} \cap \ker(\sym\Curl)$, then the problem is equivalent to the $\HsC{,\body}$-projection and the approximation $\Pm^h$ should converge to $\widetilde{\Pm}$.

\begin{figure}
    	\centering
    	\begin{subfigure}{0.48\linewidth}
    		\centering
    		\begin{tikzpicture}
    			\definecolor{asl}{rgb}{0.4980392156862745,0.,1.}
    			\definecolor{asb}{rgb}{0.,0.4,0.6}
    			\begin{loglogaxis}[
    				/pgf/number format/1000 sep={},
    				axis lines = left,
    				xlabel={degrees of freedom},
    				ylabel={$\| \widetilde{\Pm} - \Pm^h \|_{\Le} / \| \widetilde{\Pm} \|_{\Le} $ },
    				xmin=100, xmax=500000,
    				ymin=1e-4, ymax=1,
    				xtick={1e3,1e4,1e5,1e6},
    				ytick={1e-4,1e-2, 1},
    				legend style={at={(0.05,0.05)},anchor= south west},
    				ymajorgrids=true,
    				grid style=dotted,
    				]
    				\addplot[color=asl, mark=triangle] coordinates {
        ( 1218 , 0.20414945164197734 )
        ( 3771 , 0.14476878742112617 )
        ( 8556 , 0.11202820332625893 )
        ( 27630 , 0.07735160611103994 )
        ( 64056 , 0.059169543678884315 )
        ( 90531 , 0.05295534305332747 )
        ( 123450 , 0.04792173310449485 )
    				};
    				\addlegendentry{$\Y^0$}
        \addplot[color=blue, mark=diamond] coordinates {
        ( 2082 , 0.18997835840634797 )
        ( 6687 , 0.1348936863183479 )
        ( 15468 , 0.10459952132752627 )
        ( 50958 , 0.07246844547172501 )
        ( 80367 , 0.06288760126265695 )
        ( 119352 , 0.055558851951810574 )
    				};
    				\addlegendentry{$\S^0$}

    				\addplot[color=cyan, mark=square] coordinates {
        ( 2436 , 0.025330665768804564 )
        ( 7542 , 0.012116618994009374 )
        ( 17112 , 0.007040676459596872 )
        ( 55260 , 0.003222023755782561 )
        ( 86646 , 0.0023855648165126415 )
        ( 128112 , 0.0018367235925962199 )
    				};
    				\addlegendentry{$\Y^1$}
    				
    				\addplot[color=asb, mark=pentagon] coordinates {
        ( 3876 , 0.02511655264970166 )
        ( 12402 , 0.012049248175240071 )
        ( 28632 , 0.007012024521655207 )
        ( 55050 , 0.004573092860536905 )
        ( 94140 , 0.003213733904511236 )
        ( 148386 , 0.0023804256353025547 )
    				};
    				\addlegendentry{$\S^1$}
        
                    \addplot[color=olive, mark=o] coordinates {
( 1080 ,  0.02148669703299072 )
( 3224 ,  0.009471337821548546 )
( 7208 ,  0.005240809816547624 )
( 23000 ,  0.002294739958439281 )
( 53064 ,  0.0012816580151741052 )
( 102008 ,  0.0008171470438782495 )
    				};
    				\addlegendentry{$\Dp^1$}

                    \addplot[color=orange, mark=star] coordinates {
                    ( 2925 , 0.009325795934925968 )
                    ( 8379 , 0.004148536414721484 )
                    ( 18225 , 0.0023326803770812505 )
                    ( 33759 , 0.0014923810291763014 )
                    ( 56277 , 0.001036106373295894 )
                    ( 127449 , 0.0005826249240383437 )
    				};
    				\addlegendentry{$\Lag^1$}

                    \addplot[color=purple, mark=x] coordinates {
                    ( 2313 , 0.18275178552699906 )
                    ( 5100 , 0.14067288378718706 )
                    ( 15966 , 0.0960503192028097 )
                    ( 36408 , 0.07280555980207068 )
                    ( 69450 , 0.05857933821692631 )
                    ( 118116 , 0.04898719920320837 )
    				};
    				\addlegendentry{$\Ned^0$}

                    \addplot[color=violet, mark=|] coordinates {
                    ( 1572 , 0.03071849637359004 )
                    ( 4626 , 0.014720294689989227 )
                    ( 10200 , 0.008565446251486222 )
                    ( 31932 , 0.003926872968009862 )
                    ( 72816 , 0.002240923175267169 )
                    ( 138900 , 0.001446118371104673 )
    				};
    				\addlegendentry{$\Nedtwo^1$}
    				
    				\addplot[dashed,color=black, mark=none]
    				coordinates {
    					(5000, 0.5)
    					(1e+5, 0.18420157493201936)
    				};    				

                   \addplot[dashed,color=black, mark=none]
    				coordinates {
    					(5000, 5e-2)
    					(1e+5, 0.006786044041487267)
    				};    	
    			\end{loglogaxis}
    			\draw (4.5,4.9) 
    			node[anchor=south]{$\mathcal{O}(h^{1})$};
                \draw (4.7,3.15) 
    			node[anchor=south]{$\mathcal{O}(h^{2})$};
    		\end{tikzpicture}
    		\caption{}
    	\end{subfigure}
     \begin{subfigure}{0.48\linewidth}
    		\centering
    		\begin{tikzpicture}
    			\definecolor{asl}{rgb}{0.4980392156862745,0.,1.}
    			\definecolor{asb}{rgb}{0.,0.4,0.6}
    			\begin{loglogaxis}[
    				/pgf/number format/1000 sep={},
    				axis lines = left,
    				xlabel={degrees of freedom},
    				ylabel={$\| \widetilde{\Pm} - \Pm^h \|_{\Le} / \| \widetilde{\Pm} \|_{\Le} $ },
    				xmin=100, xmax=500000,
    				ymin=1e-4, ymax=1,
    				xtick={1e3,1e4,1e5,1e6},
    				ytick={1e-4,1e-2, 1},
    				legend style={at={(0.05,0.05)},anchor= south west},
    				ymajorgrids=true,
    				grid style=dotted,
    				]
    				\addplot[color=asl, mark=triangle] coordinates {
    					( 5894 ,  0.008253840553806588 )
( 14298 ,  0.0035339334037061835 )
( 26326 ,  0.0020409043400735857 )
( 41978 ,  0.001327867654118769 )
( 61254 ,  0.0008163461265295505 )
( 84154 ,  0.0006572155710791805 )
( 110678 ,  0.0004351619941417855 )
( 140826 ,  0.00033888047912723144 )
    				};
    				\addlegendentry{$\Y^2$}
        \addplot[color=blue, mark=diamond] coordinates {
                        ( 4814 ,  0.008317966284879411 )
( 11598 ,  0.0035626523233685033 )
( 21286 ,  0.0020531248469484757 )
( 33878 ,  0.0013391152871955604 )
( 49374 ,  0.0008249675174909552 )
( 67774 ,  0.0006628950781160799 )
( 89078 ,  0.0004398357118025796 )
( 113286 ,  0.0003422229818165495 )
    				};
    				\addlegendentry{$\M^2$}
    				
    				\addplot[dashed,color=black, mark=none]
    				coordinates {
    					(1e+4, 5e-2)
    					(1e+5, 0.005)
    				};    				
    			\end{loglogaxis}
    			\draw (4.9,3.15) 
    			node[anchor=south]{$\mathcal{O}(h^{3})$};
    		\end{tikzpicture}
    		\caption{}
    	\end{subfigure}
    	\caption{Relative error of the approximations of the $\Hone(\body)$-regular field.}
    	\label{fig:bench1}
    \end{figure}
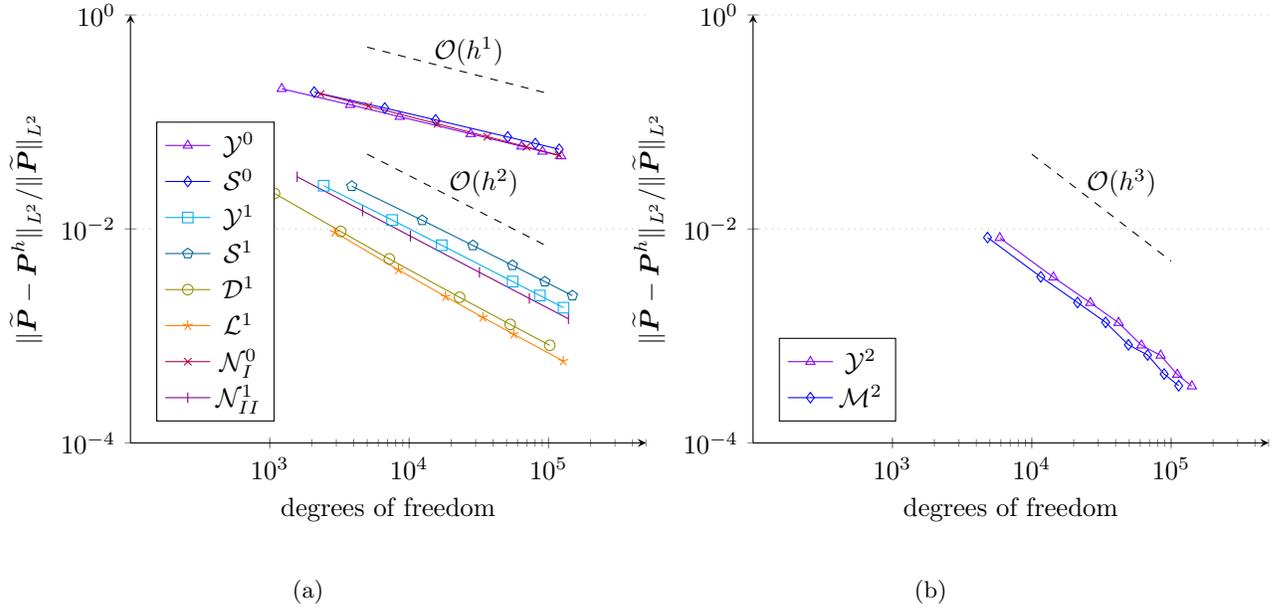

For completeness, we start with standard $\Hone(\body)$-regularity. The analytical solution is defined as
\begin{align}
    \widetilde{\Pm} = \dfrac{\sinh(x) }{10} \, \vb{e}_1 \otimes \vb{e}_1 \, ,
\end{align}
which is a smooth field $\widetilde{\Pm} \in [\Hone(\body)]^{3 \times 3} \cap \ker(\sym \Curl)$. The convergence rates for all formulations are depicted in \cref{fig:bench1}.
Clearly, all element formulations achieve optimal convergence rates. 
This is to be expected, since the Lagrangian elements satisfy $\C^0(\body)$-continuity and all other elements are of lower regularity. We note that the partial increase of the polynomial order of the second order tensor identity fields in our novel elements $\Y^p(\body)$, $\S^p(\body)$ does not yield higher convergence rates, which is the anticipated result, since the analytical solution is not a tensor identity field. 

\begin{figure}
    	\centering
    	\begin{subfigure}{0.48\linewidth}
    		\centering
    		\begin{tikzpicture}
    			\definecolor{asl}{rgb}{0.4980392156862745,0.,1.}
    			\definecolor{asb}{rgb}{0.,0.4,0.6}
    			\begin{loglogaxis}[
    				/pgf/number format/1000 sep={},
    				axis lines = left,
    				xlabel={degrees of freedom},
    				ylabel={$\| \widetilde{\Pm} - \Pm^h \|_{\Le} / \| \widetilde{\Pm} \|_{\Le} $ },
    				xmin=100, xmax=500000,
    				ymin=1e-4, ymax=1,
    				xtick={1e3,1e4,1e5,1e6},
    				ytick={1e-4,1e-3,1e-2,1e-1, 1},
    				legend style={at={(0.05,0.05)},anchor= south west},
    				ymajorgrids=true,
    				grid style=dotted,
    				]
    				\addplot[color=asl, mark=triangle] coordinates {
        ( 1218 , 0.14199241592226144 )
        ( 3771 , 0.1038856999673744 )
        ( 8556 , 0.08132818491840572 )
        ( 27630 , 0.05665818133672877 )
        ( 64056 , 0.043509820452789075 )
        ( 90531 , 0.038990133920417225 )
        ( 123450 , 0.03531971897588064 )
    				};
    				\addlegendentry{$\Y^0$}
        \addplot[color=blue, mark=diamond] coordinates {
        ( 2082 , 0.13925120525772294 )
        ( 6687 , 0.10198392065785902 )
        ( 15468 , 0.07986476110305496 )
        ( 50958 , 0.055676238738217275 )
        ( 80367 , 0.04837745224255066 )
        ( 119352 , 0.04277759159532118 )
    				};
    				\addlegendentry{$\S^0$}

    				\addplot[color=cyan, mark=square] coordinates {
        ( 2436 , 0.030549546684492262 )
        ( 7542 , 0.014229718557550439 )
        ( 17112 , 0.00819706057240546 )
        ( 55260 , 0.003729976306952729 )
        ( 86646 , 0.0027586854621050832 )
        ( 128112 , 0.002122628567664172 )
    				};
    				\addlegendentry{$\Y^1$}
    				
    				\addplot[color=asb, mark=pentagon] coordinates {
        ( 3876 , 0.030503377266451467 )
        ( 12402 , 0.014212332902973423 )
        ( 28632 , 0.008188935950770304 )
        ( 55050 , 0.0053167058826112485 )
        ( 94140 , 0.0037274386164237614 )
        ( 148386 , 0.0027570856558972274 )
    				};
    				\addlegendentry{$\S^1$}
        
                    \addplot[color=olive, mark=o] coordinates {
                    ( 1080 , 0.3949145582612925 )
                    ( 3224 , 0.3346215551879307 )
                    ( 7208 , 0.29140925381183347 )
                    ( 23000 , 0.23699661723590543 )
                    ( 53064 , 0.20385658744873858 )
                    ( 102008 , 0.1809839381278675 )
    				};
    				\addlegendentry{$\Dp^1$}

                    \addplot[color=orange, mark=star] coordinates {
                    ( 2925 , 0.1869560887921608 )
                    ( 8379 , 0.1522340721682401 )
                    ( 18225 , 0.13197926114799746 )
                    ( 33759 , 0.11820457312376786 )
                    ( 56277 , 0.10803226332758684 )
                    ( 127449 , 0.09372436728689254 )
    				};
    				\addlegendentry{$\Lag^1$}

                    \addplot[color=purple, mark=x] coordinates {
                    ( 2313 , 0.12121516792929196 )
                    ( 5100 , 0.09391650581070408 )
                    ( 15966 , 0.06443736970748729 )
                    ( 36408 , 0.048926943563323874 )
                    ( 69450 , 0.03939799501358051 )
                    ( 118116 , 0.0329611785669602 )
    				};
    				\addlegendentry{$\Ned^0$}

                    \addplot[color=violet, mark=|] coordinates {
                    ( 1572 , 0.033316618421783274 )
                    ( 4626 , 0.015328773993182264 )
                    ( 10200 , 0.008786354087433352 )
                    ( 31932 , 0.003981453643267778 )
                    ( 72816 , 0.0022616334685711553 )
                    ( 138900 , 0.0014560103024879704 )
    				};
    				\addlegendentry{$\Nedtwo^1$}
    				
    				\addplot[dashed,color=black, mark=none]
    				coordinates {
    					(5000, 0.5)
    					(1e+5, 0.3034811155014586)
    				};    
                   \addplot[dashed,color=black, mark=none]
    				coordinates {
    					(5000, 3e-2)
    					(1e+5, 0.01105209449592116)
    				};    
                    \addplot[dashed,color=black, mark=none]
    				coordinates {
    					(5000, 3e-3)
    					(1e+5, 0.00040716264248923605)
    				};    
    			\end{loglogaxis}
    			\draw (4.5,5.1) 
    			node[anchor=south]{$\mathcal{O}(h^{1/2})$};
                \draw (4.5,3.15) 
    			node[anchor=south]{$\mathcal{O}(h^{1})$};
                \draw (4.5,1.5) 
    			node[anchor=south]{$\mathcal{O}(h^{2})$};
    		\end{tikzpicture}
    		\caption{}
    	\end{subfigure}
     \begin{subfigure}{0.48\linewidth}
    		\centering
    		\begin{tikzpicture}
    			\definecolor{asl}{rgb}{0.4980392156862745,0.,1.}
    			\definecolor{asb}{rgb}{0.,0.4,0.6}
    			\begin{loglogaxis}[
    				/pgf/number format/1000 sep={},
    				axis lines = left,
    				xlabel={degrees of freedom},
    				ylabel={$\| \widetilde{\Pm} - \Pm^h \|_{\Le} / \| \widetilde{\Pm} \|_{\Le} $ },
    				xmin=100, xmax=500000,
    				ymin=1e-4, ymax=1,
    				xtick={1e3,1e4,1e5,1e6},
    				ytick={1e-4,1e-3,1e-2,1e-1, 1},
    				legend style={at={(0.05,0.05)},anchor= south west},
    				ymajorgrids=true,
    				grid style=dotted,
    				]
    				\addplot[color=asl, mark=triangle] coordinates {
    					( 5894 ,  0.013842836747087873 )
( 14298 ,  0.005112695854607393 )
( 26326 ,  0.002934581894435842 )
( 41978 ,  0.0022063192818508304 )
( 61254 ,  0.0016759064156563687 )
( 84154 ,  0.0010821953975087251 )
( 110678 ,  0.0009154500375910803 )
( 140826 ,  0.0006513949639583843 )
    				};
    				\addlegendentry{$\Y^2$}
        \addplot[color=blue, mark=diamond] coordinates {
                        ( 4814 ,  0.013910856996127301 )
( 11598 ,  0.0051450296296278245 )
( 21286 ,  0.0029550244019807353 )
( 33878 ,  0.0022169135479330974 )
( 49374 ,  0.0016815494063233787 )
( 67774 ,  0.0010875243859682327 )
( 89078 ,  0.0009191589740498611 )
( 113286 ,  0.0006539981840123787 )
    				};
    				\addlegendentry{$\M^2$}

    				\addplot[dashed,color=black, mark=none]
    				coordinates {
    					(1e+4, 5e-2)
    					(1e+5, 0.005)
    				};    	   				

    			\end{loglogaxis}
    			\draw (4.9,3.15) 
    			node[anchor=south]{$\mathcal{O}(h^{3})$};
    		\end{tikzpicture}
    		\caption{}
    	\end{subfigure}
    	\caption{Relative error of the field with a jumping normal component.}
    	\label{fig:bench2}
    \end{figure}
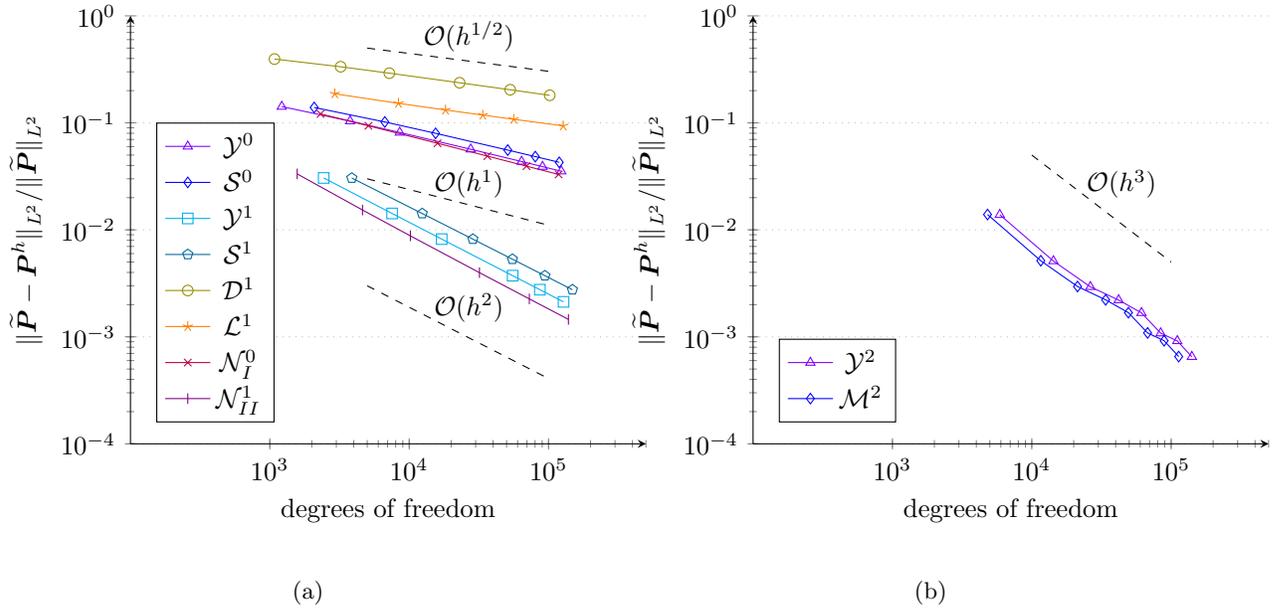

Next we consider a field with a jumping normal component.
The field reads 
\begin{align}
    \widetilde{\Pm} = \left \{ \begin{matrix}
        \cos(x) (\vb{e}_1 \otimes \vb{e}_1 + \vb{e}_2 \otimes \vb{e}_1 + \vb{e}_3 \otimes \vb{e}_1) & \text{for} &  -1< x <1 \\
        \sin(x) (\vb{e}_1 \otimes \vb{e}_1 + \vb{e}_2 \otimes \vb{e}_1 + \vb{e}_3 \otimes \vb{e}_1) & \text{otherwise} 
    \end{matrix} \right . \, .
\end{align}
The field lives in $\HC{,\body}$ but not in $[\Hone(\body)]^{3 \times 3}$, since $(\cdot) \otimes \vb{e}_1$ represent the normal components $(\cdot) \otimes \vb{n}$ on planes orthogonal to the $x$-axis, see \cref{fig:exdom}. 
Evidently, this implies that the normal components jump at $x = \pm 1$.
The results for all formulations are given in \cref{fig:bench2}.
As expected, the Lagrangian formulation $[\Lag^1(\body)]^{3 \times 3}$ is now incapable of achieving optimal convergence and in fact, its convergence rate is reduced to square-root. 
The same behaviour is observed for the $\Dp^1(\body)$-formulation, which confirms its non-conforming nature for $\HC{,\body}$-fields.
In comparison, the N\'ed\'elec elements continue to converge optimally, as they can account for the jump in the normal component. The same holds true for our novel elements $\Y^0(\body)$, $\S^0(\body)$, $\Y^1(\body)$ and $\S^1(\body)$, which enrich the N\'ed\'elec space with higher order discontinuous identity fields. More importantly, the $\Y^2(\body)$ and $\M^2(\body)$ formulations exhibit optimal cubic convergence. The latter proves that the extension in \cref{sec:arb} respects the regularity of $\HC{,\body}$. 
\begin{figure}
    	\centering
     \begin{subfigure}{0.48\linewidth}
    		\centering
    		\includegraphics[width=1\linewidth]{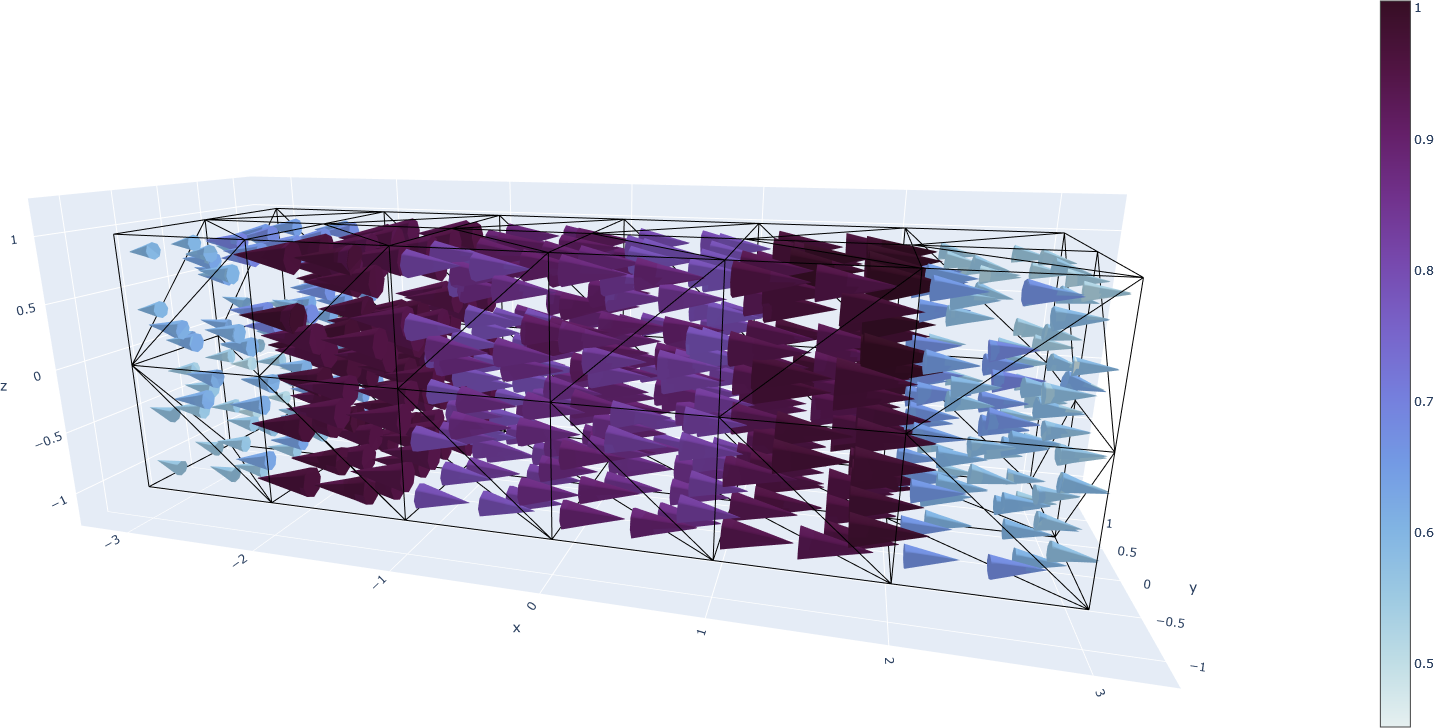}
    		\caption{}
    	\end{subfigure}
     \begin{subfigure}{0.48\linewidth}
    		\centering
    		\includegraphics[width=1\linewidth]{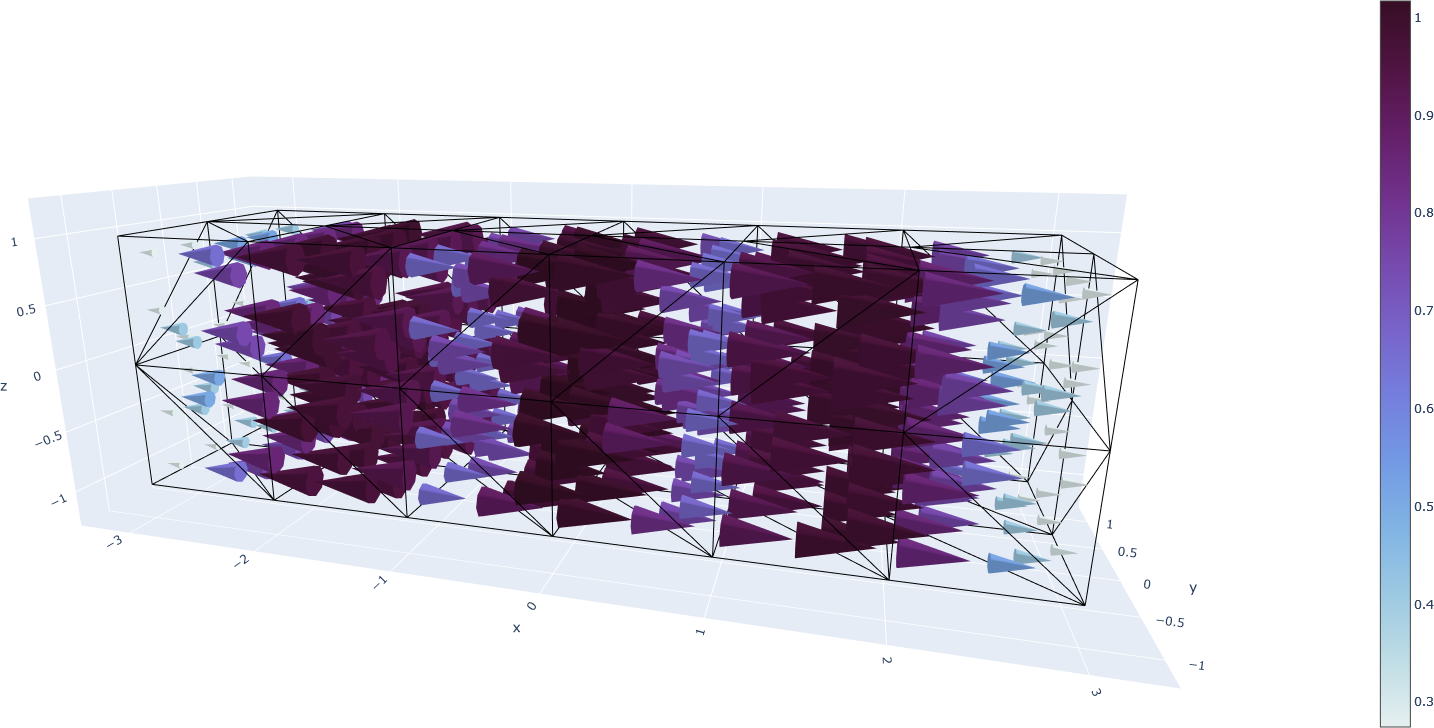}
    		\caption{}
    	\end{subfigure}
     \begin{subfigure}{0.48\linewidth}
    		\centering
    		\includegraphics[width=1\linewidth]{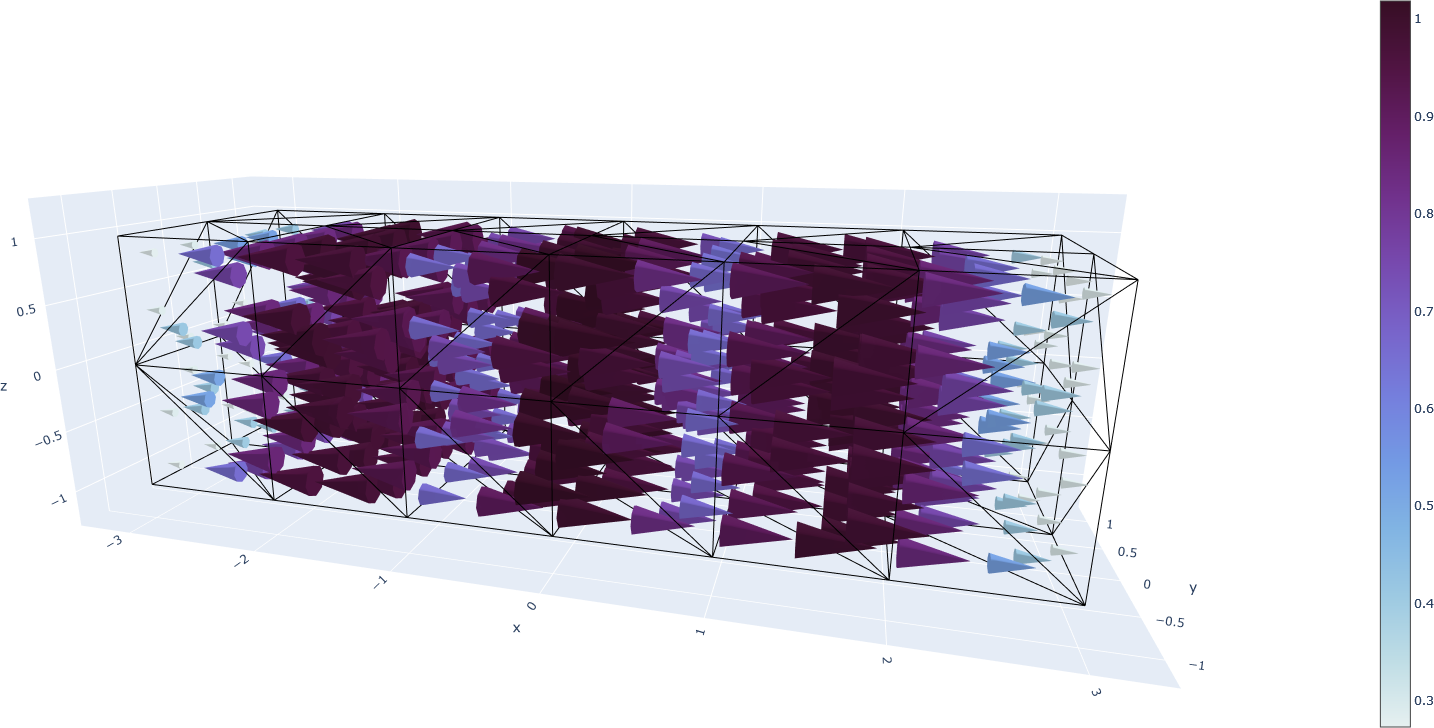}
    		\caption{}
    	\end{subfigure}
     \begin{subfigure}{0.48\linewidth}
    		\centering
    		\includegraphics[width=1\linewidth]{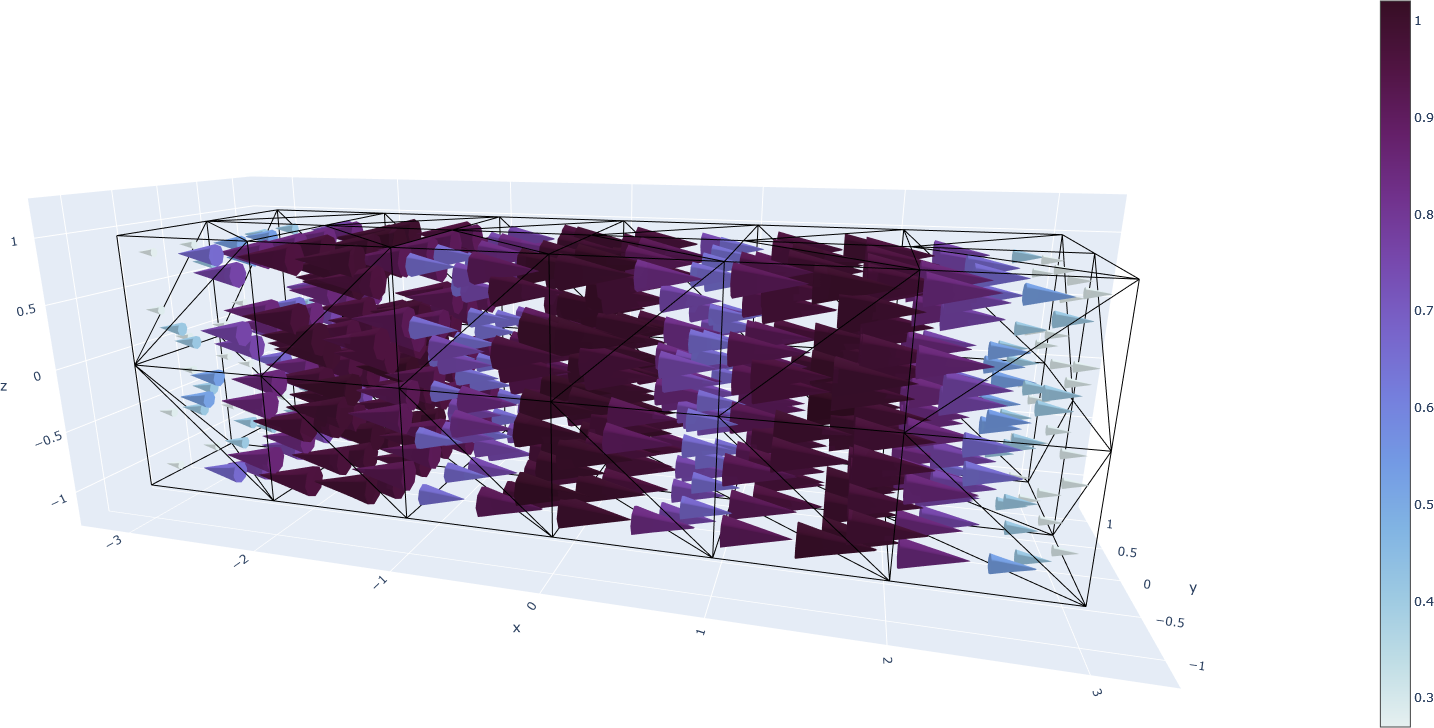}
    		\caption{}
    	\end{subfigure}
     \begin{subfigure}{0.48\linewidth}
    		\centering
    		\includegraphics[width=1\linewidth]{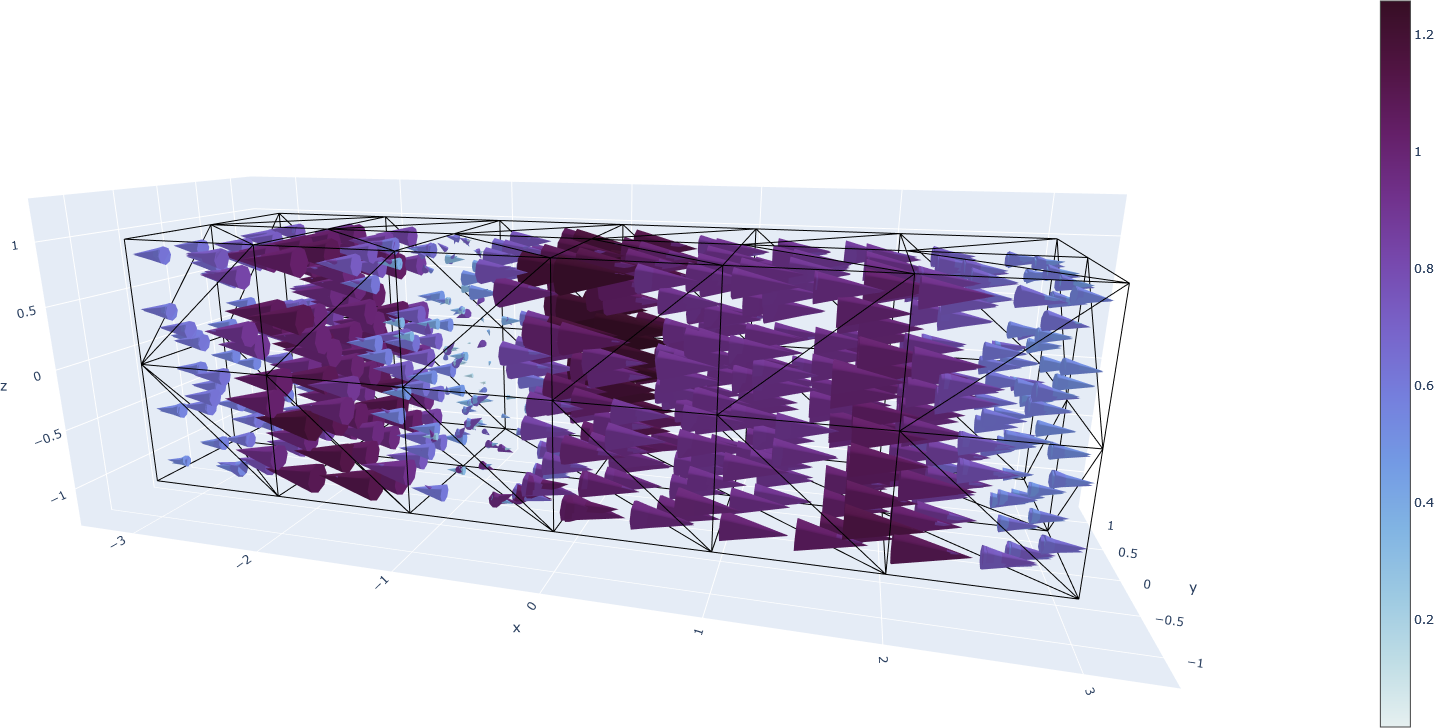}
    		\caption{}
    	\end{subfigure}
    	\begin{subfigure}{0.48\linewidth}
    		\centering
    		\includegraphics[width=1\linewidth]{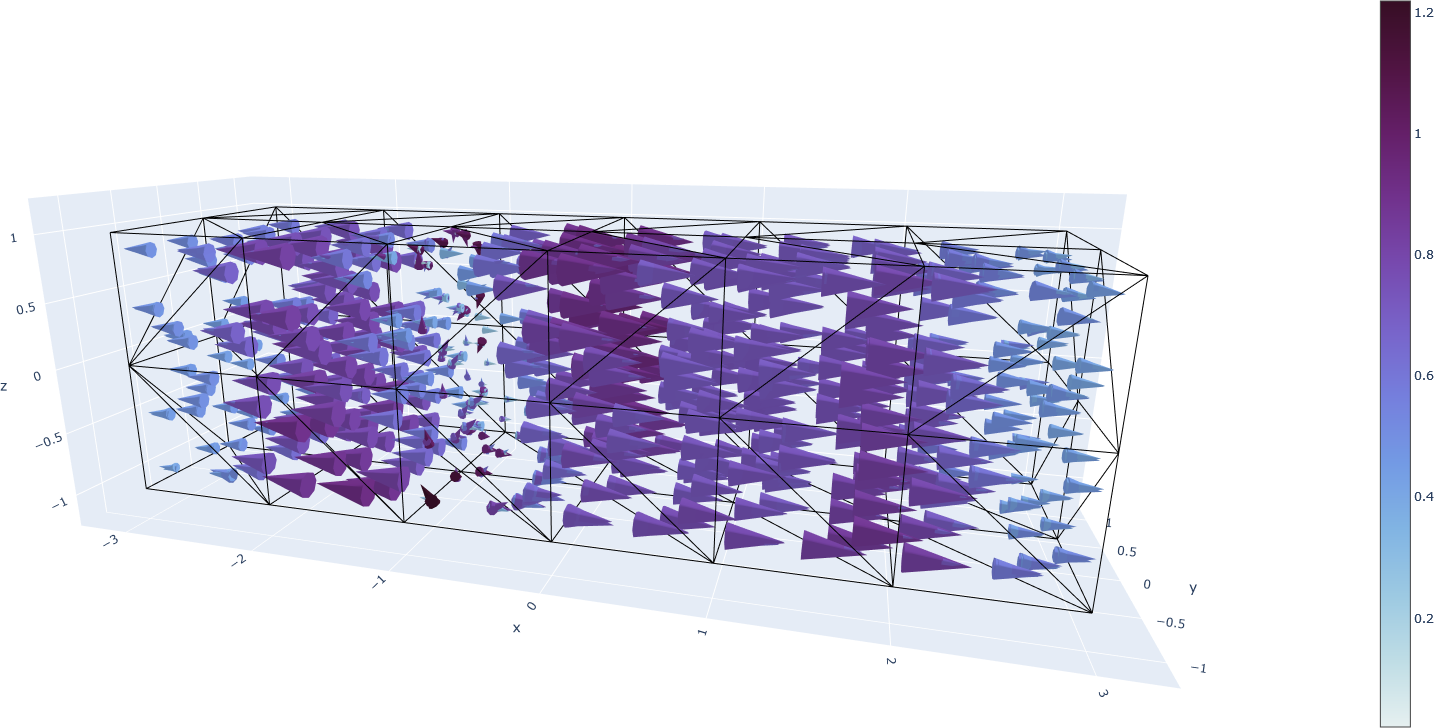}
    		\caption{}
    	\end{subfigure}
    	\caption{The jumping normal field on a coarse mesh of $144$ elements. Each cone represents one row in the matrix and they overlap at a point due to the similar orientation and magnitude. The depictions in (a) and (b) are the of the N\'ed\'elec elements $[\Ned^0(\body)]^3$ and $[\Nedtwo^1(\body)]^3$, which represent a decent approximation of the analytical solution. In (c) and (d) we depict the $\Y^1(\body)$ and $\Y^2(\body)$ approximations, which are clearly very similar to the solution by the linear N\'edelec element of the second type. In (e) and (f) one sees the solutions by $[\Lag^1(\body)]^{3 \times 3}$ and $\Dp^1(\body)$, which cannot correctly account for the jump in normal direction. }
    	\label{fig:ex1hc}
    \end{figure}
The results of the approximations on a coarse mesh of $144$ elements are depicted in \cref{fig:ex1hc}. Clearly, our novel elements match the performance of N\'ed\'elec elements. 

In the third benchmark we consider the semi-smooth jumping identity field
\begin{align}
    \widetilde{\Pm} = \left \{ \begin{matrix}
        2\sin(x+2y-3z)\one & \text{for} &  -1< x <1 \\
        \sin(x+2y-3z) \one & \text{otherwise} 
    \end{matrix} \right . \, .
\end{align}
The field is in $\HsC{,\body}$, but not in $\HC{,\body}$ or $[\Hone(\body)]^{3\times 3}$. The fact that $\widetilde{\Pm}$ is not in $[\Hone(\body)]^{3\times 3}$, is evident, since it jumps at $x = \pm 1$. The field is also not in $\HC{,\body}$ because the $\vb{e}_2 \otimes \vb{e}_2$- and $\vb{e}_3 \otimes \vb{e}_3$ components are not normal to the $x$-plane, on which the jump occurs. Finally, the field is in $\HsC{,\body}$ since identity tensor fields are in the kernel of its corresponding trace operator $\ker(\tr_{\HsC{}})$.
We examine the convergence rates in \cref{fig:bench3}.
\begin{figure}
    	\centering
    	\begin{subfigure}{0.48\linewidth}
    		\centering
    		\begin{tikzpicture}
    			\definecolor{asl}{rgb}{0.4980392156862745,0.,1.}
    			\definecolor{asb}{rgb}{0.,0.4,0.6}
    			\begin{loglogaxis}[
    				/pgf/number format/1000 sep={},
    				axis lines = left,
    				xlabel={degrees of freedom},
    				ylabel={$\| \widetilde{\Pm} - \Pm^h \|_{\Le} / \| \widetilde{\Pm} \|_{\Le} $ },
    				xmin=100, xmax=500000,
    				ymin=1e-3, ymax=1,
    				xtick={1e3,1e4,1e5,1e6},
    				ytick={1e-3,1e-2,1e-1, 1},
    				legend style={at={(0.05,0.05)},anchor= south west},
    				ymajorgrids=true,
    				grid style=dotted,
    				]
    				\addplot[color=asl, mark=triangle] coordinates {
        ( 1218 , 0.2963205179763696 )
        ( 3771 , 0.20301432380280757 )
        ( 8556 , 0.15135409308499143 )
        ( 27630 , 0.100384681985379 )
        ( 64056 , 0.07515648041333232 )
        ( 90531 , 0.06677567703514743 )
        ( 123450 , 0.06007922792317465 )
    				};
    				\addlegendentry{$\Y^0$}
        \addplot[color=blue, mark=diamond] coordinates {
        ( 2082 , 0.08099010609324685 )
        ( 6687 , 0.03249564391848271 )
        ( 15468 , 0.018442006005345603 )
        ( 50958 , 0.008227118123886662 )
        ( 80367 , 0.006049041267258999 )
        ( 119352 , 0.004633589731225922 )
    				};
    				\addlegendentry{$\S^0$}

    				\addplot[color=cyan, mark=square] coordinates {
        ( 2436 , 0.05862641205229836 )
        ( 7542 , 0.035208814693073125 )
        ( 17112 , 0.027997935014355622 )
        ( 55260 , 0.02195914317967153 )
        ( 86646 , 0.02022397662265587 )
        ( 128112 , 0.018868527548473193 )
    				};
    				\addlegendentry{$\Y^1$}
    				
    				\addplot[color=asb, mark=pentagon] coordinates {
        ( 3876 , 0.014207462062846472 )
        ( 12402 , 0.011593702591543001 )
        ( 28632 , 0.011172879351900549 )
        ( 55050 , 0.01052905220167135 )
        ( 94140 , 0.009891445556994873 )
        ( 148386 , 0.009319363858538314 )
    				};
    				\addlegendentry{$\S^1$}
        
                    \addplot[color=olive, mark=o] coordinates {
                    ( 1080 , 0.15825442660391614 )
                    ( 3224 , 0.07127090568088103 )
                    ( 7208 , 0.040585327424466536 )
                    ( 23000 , 0.018200338846796466 )
                    ( 53064 , 0.010270073694274746 )
                    ( 102008 , 0.006582479851065496 )
    				};
    				\addlegendentry{$\Dp^1$}

                    \addplot[color=orange, mark=star] coordinates {
                    ( 2925 , 0.13322845847104992 )
                    ( 8379 , 0.09432853212179958 )
                    ( 18225 , 0.07892303033954669 )
                    ( 33759 , 0.06979914780149965 )
                    ( 56277 , 0.06343743919222655 )
                    ( 127449 , 0.054790018513037136 )
    				};
    				\addlegendentry{$\Lag^1$}

                    \addplot[color=purple, mark=x] coordinates {
                    ( 2313 , 0.5182693313757134 )
                    ( 5100 , 0.40468063005725335 )
                    ( 15966 , 0.28097237355976473 )
                    ( 36408 , 0.21605144647118454 )
                    ( 69450 , 0.17626329838495208 )
                    ( 118116 , 0.14941793796907502 )
    				};
    				\addlegendentry{$\Ned^0$}

                    \addplot[color=violet, mark=|] coordinates {
                    ( 1572 , 0.3391890806760572 )
                    ( 4626 , 0.16807409230017423 )
                    ( 10200 , 0.11392133880780854 )
                    ( 31932 , 0.0783450293724329 )
                    ( 72816 , 0.06475213431157426 )
                    ( 138900 , 0.05699832488422488 )
    				};
    				\addlegendentry{$\Nedtwo^1$}
    				
    				\addplot[dashed,color=black, mark=none]
    				coordinates {
    					(5000, 0.6)
    					(1e+5, 0.22104188991842322)
    				};    		

                   \addplot[dashed,color=black, mark=none]
    				coordinates {
    					(20000, 3.5e-2)
    					(1.5e+5, 0.02501638202834432)
    				};    	 

                    \addplot[dashed,color=black, mark=none]
    				coordinates {
    					(5000, 1e-2)
    					(1e+5, 0.0013572088082974534)
    				};    	 
    			\end{loglogaxis}
    			\draw (4.5,5) 
    			node[anchor=south]{$\mathcal{O}(h^{1})$};
                \draw (5.1,2.75) 
    			node[anchor=south]{$\mathcal{O}(h^{1/2})$};
                \draw (4.6,1.1) 
    			node[anchor=south]{$\mathcal{O}(h^{2})$};
    		\end{tikzpicture}
    		\caption{}
    	\end{subfigure}
     \begin{subfigure}{0.48\linewidth}
    		\centering
    		\begin{tikzpicture}
    			\definecolor{asl}{rgb}{0.4980392156862745,0.,1.}
    			\definecolor{asb}{rgb}{0.,0.4,0.6}
    			\begin{loglogaxis}[
    				/pgf/number format/1000 sep={},
    				axis lines = left,
    				xlabel={degrees of freedom},
    				ylabel={$\| \widetilde{\Pm} - \Pm^h \|_{\Le} / \| \widetilde{\Pm} \|_{\Le} $ },
    				xmin=100, xmax=500000,
    				ymin=1e-3, ymax=1,
    				xtick={1e3,1e4,1e5,1e6},
    				ytick={1e-3,1e-2,1e-1, 1},
    				legend style={at={(0.05,0.05)},anchor= south west},
    				ymajorgrids=true,
    				grid style=dotted,
    				]
    				\addplot[color=asl, mark=triangle] coordinates {
    					( 7496 ,  0.026074001195626067 )
         ( 13904 ,  0.02015406937922714 )
( 23826 ,  0.015542267835413646 )
( 54784 ,  0.014089150763991526 )
( 64270 ,  0.013581945466269197 )
( 105050 ,  0.012808044632214173 )
    				};
    				\addlegendentry{$\Y^2$}
        \addplot[color=blue, mark=diamond] coordinates {
        ( 4814 ,  0.16548303953246898 )
                          ( 6056 ,  0.12080791016030785 )
( 18966 ,  0.044003484767024655 )
( 43264 ,  0.03261577970892015 )
( 50590 ,  0.02949420168388276 )
( 82550 ,  0.02556228949984371 )
( 122096 ,  0.02444722769881936 )
( 140424 ,  0.022280675511182017 )
    				};
    				\addlegendentry{$\M^2$}
    				
    				\addplot[dashed,color=black, mark=none]
    				coordinates {
    					(1e+4, 0.005)
    					(1e+5, 0.0034064603452898068)
    				};    	   				

    			\end{loglogaxis}
    			\draw (4.85,1.1) 
    			node[anchor=south]{$\mathcal{O}(h^{1/2})$};
    		\end{tikzpicture}
    		\caption{}
    	\end{subfigure}
    	\caption{Relative error of the jumping identity field.}
    	\label{fig:bench3}
    \end{figure}
Both the Lagrangian $[\Lag^1(\body)]^{3 \times 3}$ and the second type N\'ed\'elec $[\Nedtwo^1(\body)]^3$ formulations exhibit suboptimal convergences. In contrast, the $[\Ned^0(\body)]^3$, $\Dp^1(\body)$, the $\Y^0(\body)$ and $\S^0(\body)$ formulations converge optimally. 
Latter confirms the optimality of the novel elements for $\HsC{,\body}$-fields.
The result is also interesting in the case of $[\Ned^0(\body)]^3$, since it does not represent a conforming discretisation of $\HsC{,\body}$ with a minimal regularity, seeing as it contains a non-jumping constant identity field. However, it appears that the jump of the tangential projection of the constant identity field at $x =\pm 1$ does not dominate the produced error. 
The $\Y^1(\body)$ and $\S^1(\body)$ formulations exhibit suboptimal convergence rates. 
However, their added jumping identity fields allow them to further reduce the relative error.  

\begin{figure}
    	\centering
    	\begin{subfigure}{0.48\linewidth}
    		\centering
    		\begin{tikzpicture}
    			\definecolor{asl}{rgb}{0.4980392156862745,0.,1.}
    			\definecolor{asb}{rgb}{0.,0.4,0.6}
    			\begin{loglogaxis}[
    				/pgf/number format/1000 sep={},
    				axis lines = left,
    				xlabel={degrees of freedom},
    				ylabel={$\| \widetilde{\Pm} - \Pm^h \|_{\Le} / \| \widetilde{\Pm} \|_{\Le} $ },
    				xmin=100, xmax=500000,
    				ymin=1e-4, ymax=1,
    				xtick={1e3,1e4,1e5,1e6},
    				ytick={1e-4,1e-3,1e-2,1e-1, 1},
    				legend style={at={(0.05,0.05)},anchor= south west},
    				ymajorgrids=true,
    				grid style=dotted,
    				]
    				\addplot[color=asl, mark=triangle] coordinates {
        ( 1218 , 0.09758354923376852 )
        ( 3771 , 0.06608121642193518 )
        ( 8556 , 0.05006137975135952 )
        ( 27630 , 0.03386861607336075 )
        ( 64056 , 0.025447457679863277 )
        ( 90531 , 0.022665475013317535 )
        ( 123450 , 0.020437559426730093 )
    				};
    				\addlegendentry{$\Y^0$}
        \addplot[color=blue, mark=diamond] coordinates {
        ( 2082 , 0.009060243660397482 )
        ( 6687 , 0.004061854230153476 )
        ( 15468 , 0.002288945315077708 )
        ( 50958 , 0.0010181632603557071 )
        ( 80367 , 0.0007481447707554953 )
        ( 119352 , 0.0005728479581509543 )
    				};
    				\addlegendentry{$\S^0$}

    				\addplot[color=cyan, mark=square] coordinates {
        ( 2436 , 0.08709304891038404 )
        ( 7542 , 0.0724693916390017 )
        ( 17112 , 0.0633681046871743 )
        ( 55260 , 0.052254155551395426 )
        ( 86646 , 0.048517035818722236 )
        ( 128112 , 0.045481885267897 )
    				};
    				\addlegendentry{$\Y^1$}
    				
    				\addplot[color=asb, mark=pentagon] coordinates {
        ( 3876 , 0.04267515447140537 )
        ( 12402 , 0.03558029948245373 )
        ( 28632 , 0.03105347586997763 )
        ( 55050 , 0.02791044460697814 )
        ( 94140 , 0.025561147592463637 )
        ( 148386 , 0.02371978114879829 )
    				};
    				\addlegendentry{$\S^1$}
        

                    \addplot[color=orange, mark=star] coordinates {
                    ( 2925 , 0.2594845868235602 )
                    ( 8379 , 0.21439911999271827 )
                    ( 18225 , 0.1867384091956289 )
                    ( 33759 , 0.16759474554414192 )
                    ( 56277 , 0.15334094481136654 )
                    ( 127449 , 0.13317586666782125 )
    				};
    				\addlegendentry{$\Lag^1$}

                    \addplot[color=purple, mark=x] coordinates {
                    ( 2313 , 0.28226012722603483 )
                    ( 5100 , 0.24940926059148186 )
                    ( 15966 , 0.20733986752643221 )
                    ( 36408 , 0.1809985369641866 )
                    ( 69450 , 0.16260429356704667 )
                    ( 118116 , 0.14884798385297368 )
    				};
    				\addlegendentry{$\Ned^0$}

                    \addplot[color=violet, mark=|] coordinates {
                    ( 1572 , 0.27056355859797304 )
                    ( 4626 , 0.23221198897554815 )
                    ( 10200 , 0.2054165539033437 )
                    ( 31932 , 0.17168251527246398 )
                    ( 72816 , 0.1503624315657565 )
                    ( 138900 , 0.13537450050181576 )
    				};
    				\addlegendentry{$\Nedtwo^1$}
    				
    				\addplot[dashed,color=black, mark=none]
    				coordinates {
    					(1e+4, 0.01)
    					(1e+5, 0.00464158883361278)
    				}; 
                    \addplot[dashed,color=black, mark=none]
    				coordinates {
    					(1e+4, 0.4)
    					(1e+5, 0.2725168276231845)
    				}; 

                   \addplot[dashed,color=black, mark=none]
    				coordinates {
    					(1e+4, 0.0005)
    					(1e+5, 0.00010772173450159421)
    				}; 

    			\end{loglogaxis}
    			\draw (4.7,2.6) 
    			node[anchor=south]{$\mathcal{O}(h^{1})$};
       \draw (4.7,5.) 
    			node[anchor=south]{$\mathcal{O}(h^{1/2})$};
       \draw (4.8,0.5) 
    			node[anchor=south]{$\mathcal{O}(h^{2})$};
    		\end{tikzpicture}
    		\caption{}
    	\end{subfigure}
     \begin{subfigure}{0.48\linewidth}
    		\centering
    		\begin{tikzpicture}
    			\definecolor{asl}{rgb}{0.4980392156862745,0.,1.}
    			\definecolor{asb}{rgb}{0.,0.4,0.6}
    			\begin{loglogaxis}[
    				/pgf/number format/1000 sep={},
    				axis lines = left,
    				xlabel={degrees of freedom},
    				ylabel={$\| \widetilde{\Pm} - \Pm^h \|_{\Le} / \| \widetilde{\Pm} \|_{\Le} $ },
    				xmin=100, xmax=500000,
    				ymin=1e-4, ymax=1,
    				xtick={1e3,1e4,1e5,1e6},
    				ytick={1e-4,1e-3,1e-2,1e-1, 1},
    				legend style={at={(0.05,0.05)},anchor= south west},
    				ymajorgrids=true,
    				grid style=dotted,
    				]
    				\addplot[color=asl, mark=triangle] coordinates {
    					( 7496 ,  0.06652320994694018 )
( 23826 ,  0.055469582630965776 )
( 54784 ,  0.04896877948070954 )
( 105050 ,  0.04413612723677691 )
    				};
    				\addlegendentry{$\Y^2$}
        \addplot[color=blue, mark=diamond] coordinates {
                        ( 6056 ,  0.1153154870854622 )
( 18966 ,  0.0965724196193788 )
( 43264 ,  0.08532428120129093 )
( 82550 ,  0.07714644532391683 )
    				};
    				\addlegendentry{$\M^2$}
    				
    				\addplot[dashed,color=black, mark=none]
    				coordinates {
    					(1e+4, 0.01)
    					(1e+5, 0.0068129206905796135)
    				};    	   				

    			\end{loglogaxis}
    			\draw (4.8,2.7) 
    			node[anchor=south]{$\mathcal{O}(h^{1/2})$};
    		\end{tikzpicture}
    		\caption{}
    	\end{subfigure}
    	\caption{Relative error of the jumping constant identity field.}
    	\label{fig:bench4}
    \end{figure}
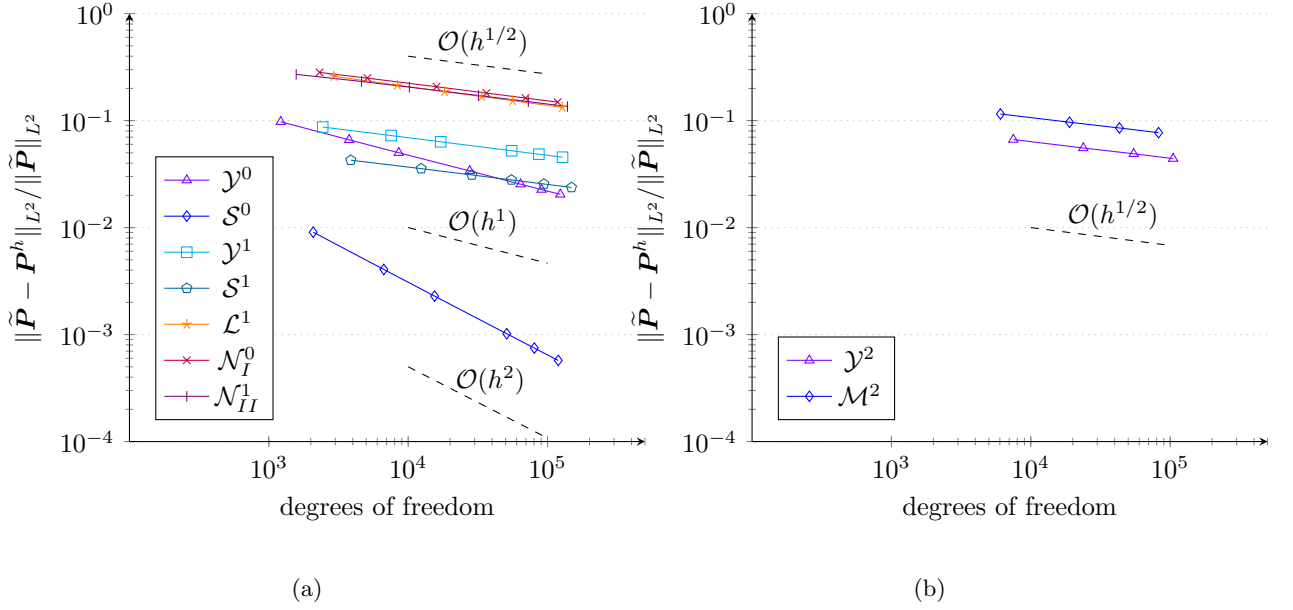

\begin{figure}
    	\centering
    	\begin{subfigure}{0.48\linewidth}
    		\centering
    		\includegraphics[width=1\linewidth]{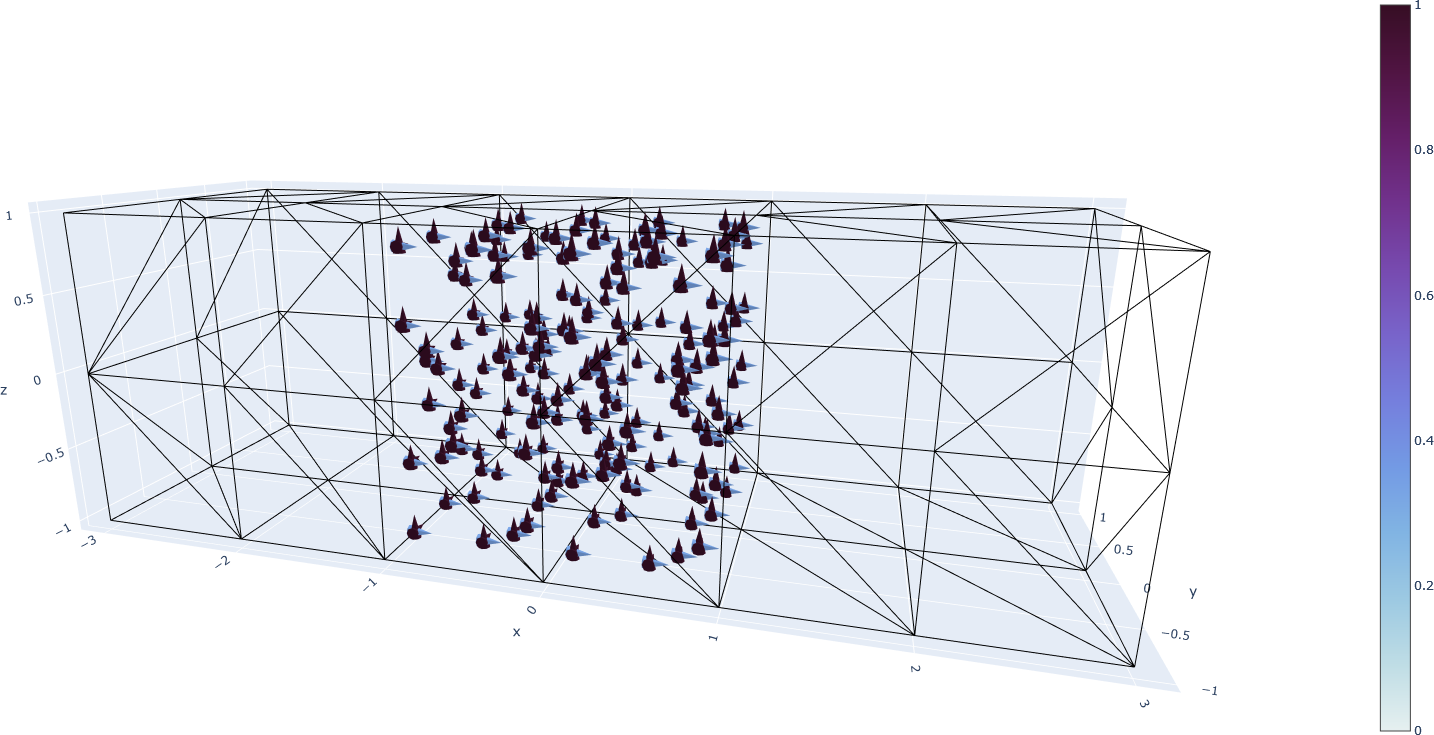}
    		\caption{}
      \label{fig:a}
    	\end{subfigure}
     \begin{subfigure}{0.48\linewidth}
    		\centering
    		\includegraphics[width=1\linewidth]{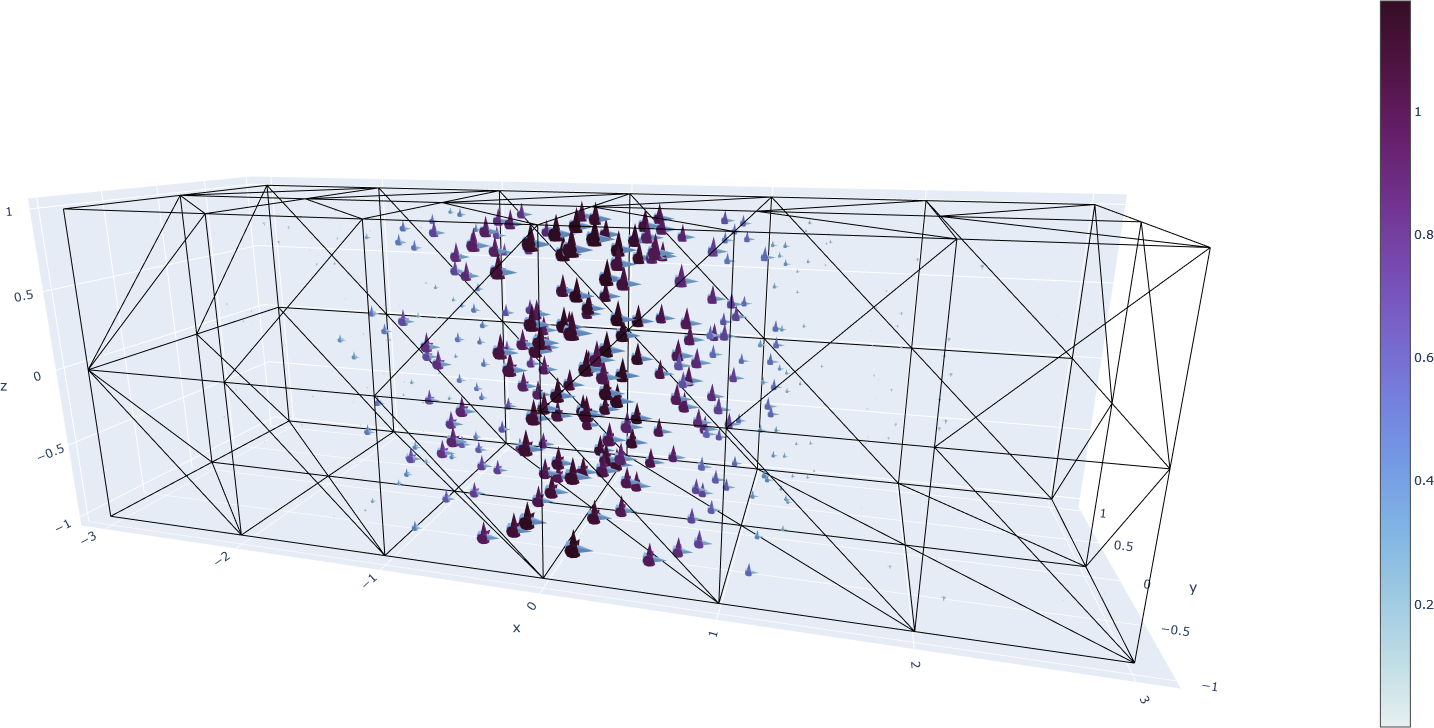}
    		\caption{}
    	\end{subfigure}
     \begin{subfigure}{0.48\linewidth}
    		\centering
    		\includegraphics[width=1\linewidth]{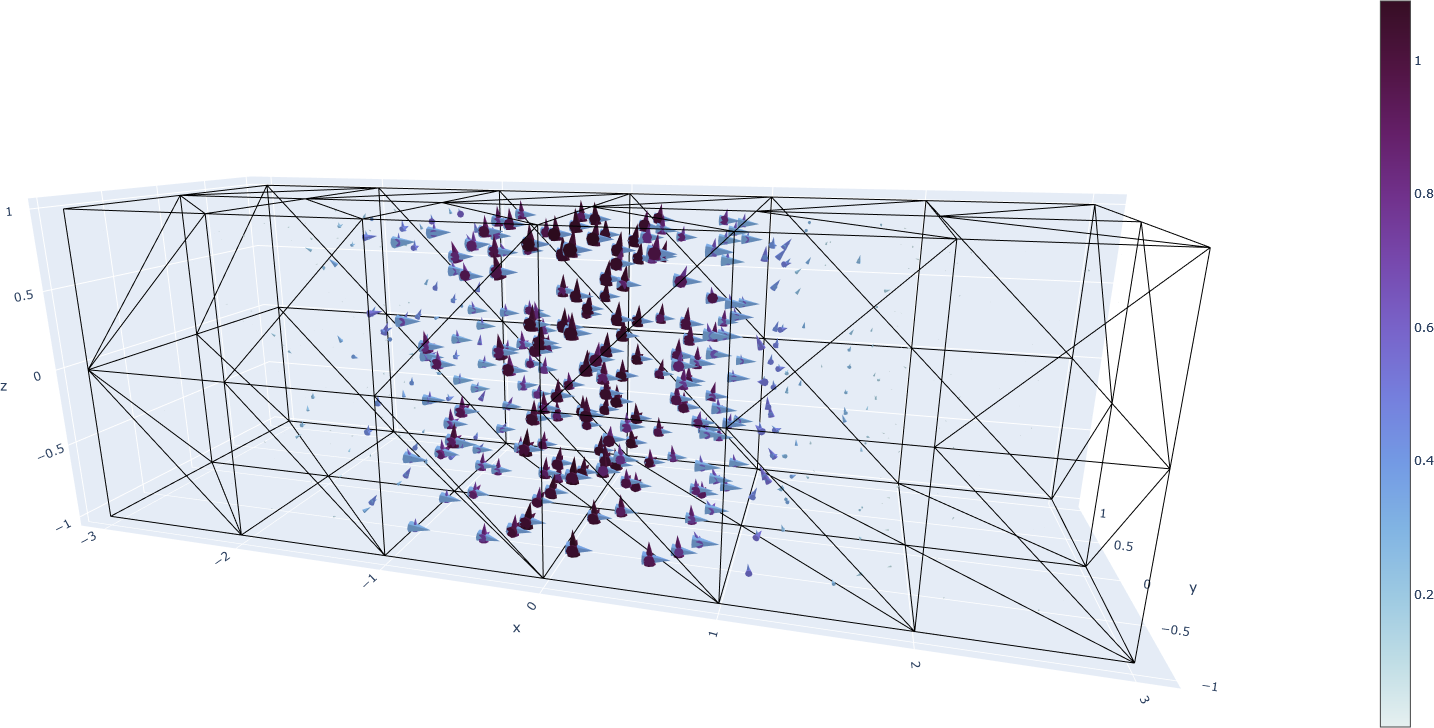}
    		\caption{}
    	\end{subfigure}
     \begin{subfigure}{0.48\linewidth}
    		\centering
    		\includegraphics[width=1\linewidth]{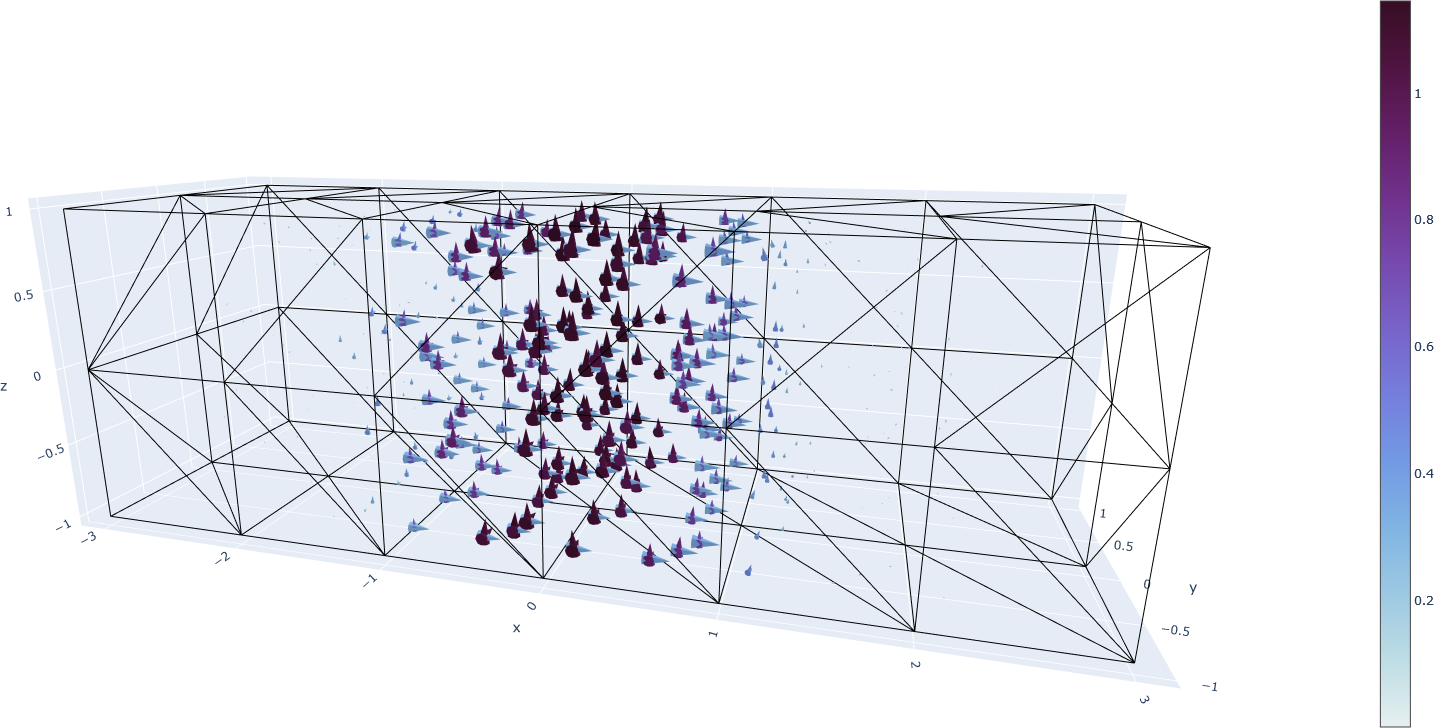}
    		\caption{}
    	\end{subfigure}
     \begin{subfigure}{0.48\linewidth}
    		\centering
    		\includegraphics[width=1\linewidth]{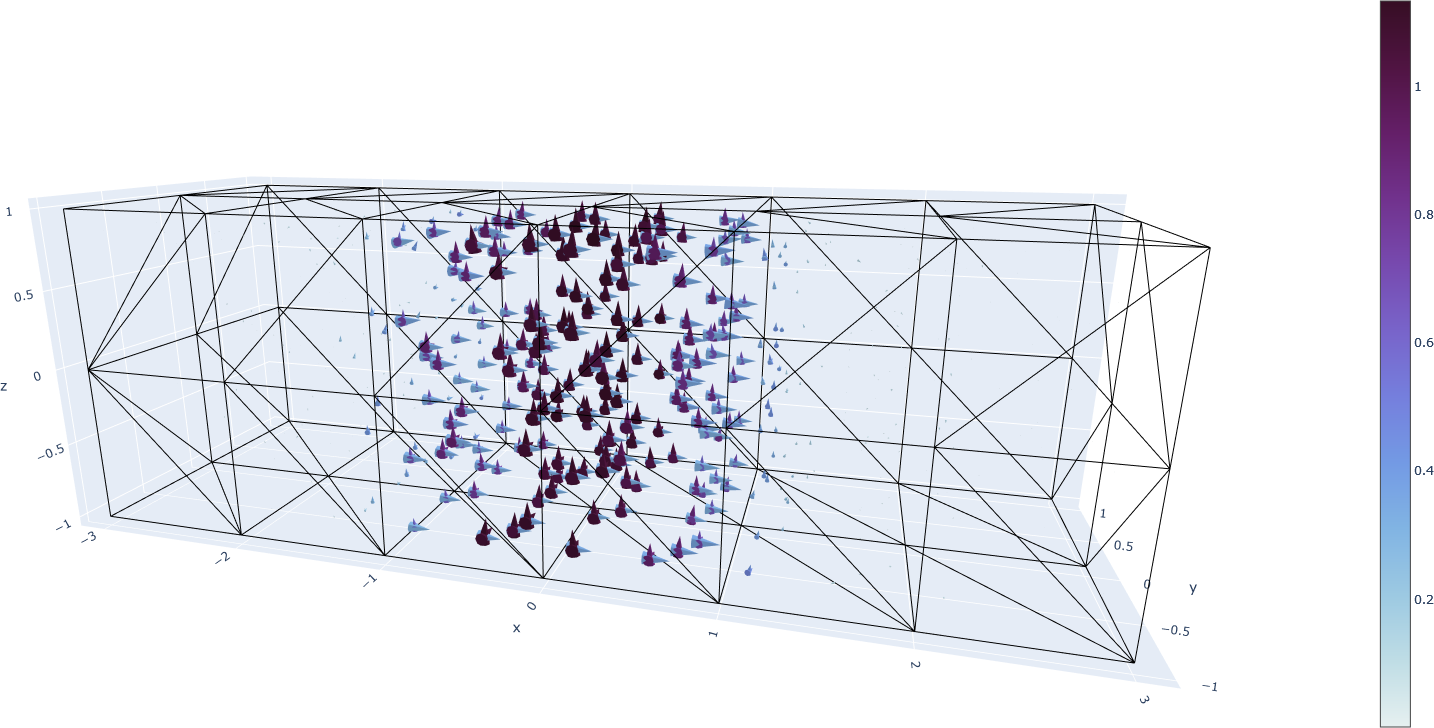}
    		\caption{}
    	\end{subfigure}
     \begin{subfigure}{0.48\linewidth}
    		\centering
    		\includegraphics[width=1\linewidth]{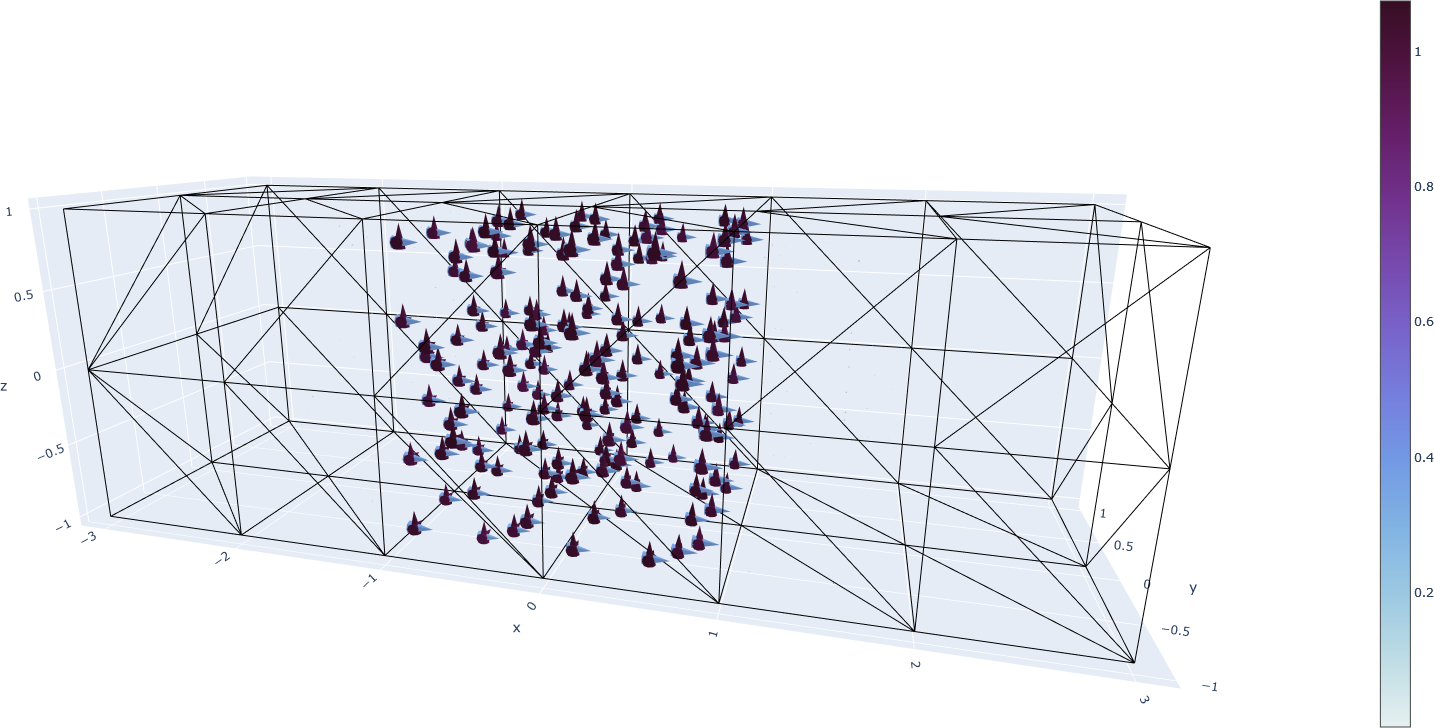}
    		\caption{}
      \label{fig:j}
    	\end{subfigure}
        \begin{subfigure}{0.48\linewidth}
    		\centering
    		\includegraphics[width=1\linewidth]{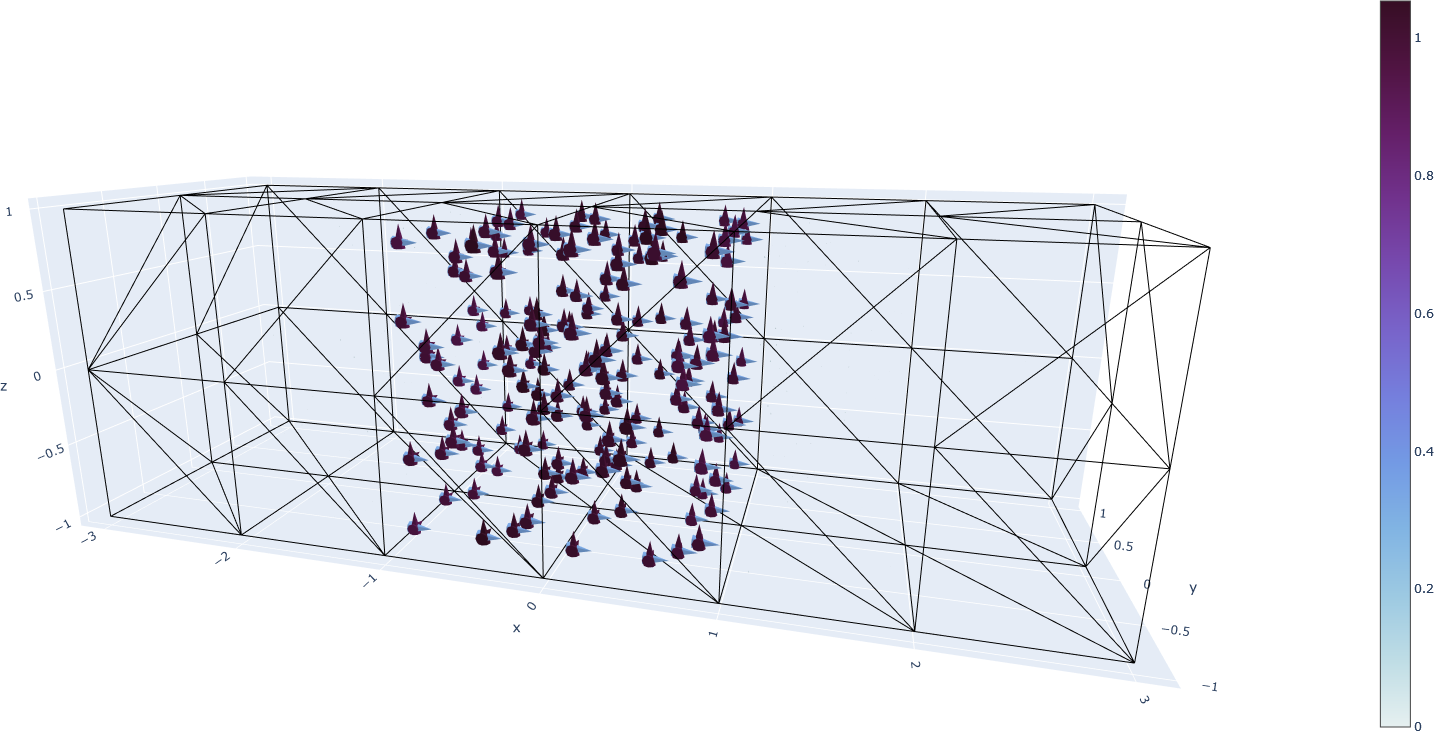}
    		\caption{}
    	\end{subfigure}
     \begin{subfigure}{0.48\linewidth}
    		\centering
    		\includegraphics[width=1\linewidth]{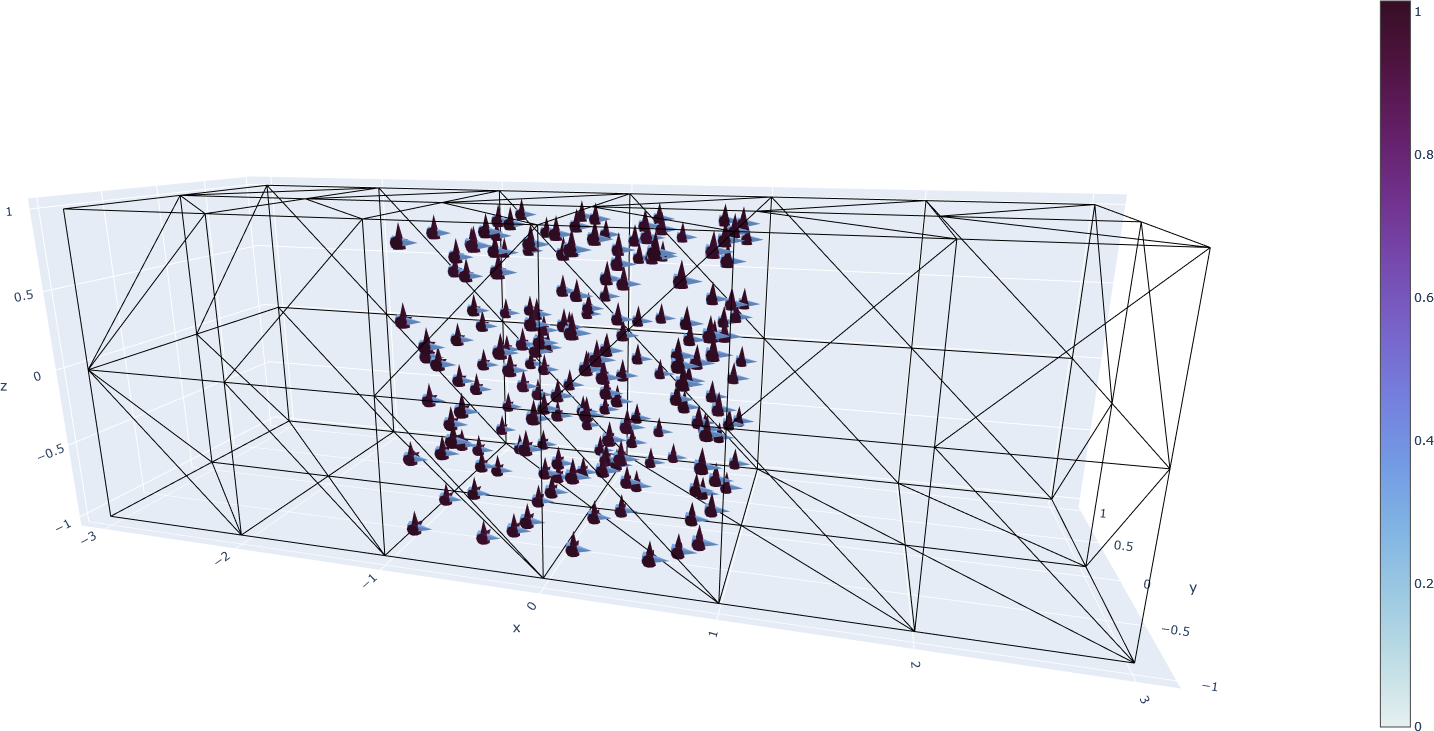}
    		\caption{}
    	\end{subfigure}
    	\caption{The jumping constant identity on a coarse mesh of $144$ elements. Each cone is a row in the matrix. The depiction in (a) is for $\Dp^1(\body)$ and represents the exact solution. In (b) one sees the solution by $[\Lag^1(\body)]^{3\times 3}$, which is followed by $[\Ned^0(\body)]^3$ and $[\Nedtwo^1(\body)]^3$ in (c) and (d). In (e)-(h) one finds the solutions for the pairs $\Y^1(\body)$, $\S^1(\body)$ and $\M^2(\body)$, $\Y^2(\body)$, respectively.}
    	\label{fig:ex1}
    \end{figure}

Lastly, we consider the constant-valued jumping identity field
\begin{align}
    \widetilde{\Pm} = \left \{ \begin{matrix}
        \one & \text{for} &  -1< x <1 \\
        0 & \text{otherwise} 
    \end{matrix} \right . \, ,
\end{align}
which lives in $\HsC{,\body}$. This field is especially relevant to our investigation, since our novel lowest order element is incapable of letting the constant identity jump. As such, it is of interest to examine the capacity of the novel elements to correct for the induced error.
The convergence estimates are depicted in \cref{fig:bench4}.
We observe that both the $\Y^0(\body)$ and $\S^0(\body)$ elements converge optimally, which implies that these formulations are in fact ideal for discretisations of $\HsC{,\body}$. 
The remaining formulations are all suboptimal in their convergence rates. The $\Dp^1(\body)$ formulation is not measured here since the field $\widetilde{\Pm}$ is contained in that discrete space, such that no error is induced.
The convergence rates represent one aspect of the error measurement. The second measure is given by the constant placed before the convergence slope.
While the $\Y^1(\body)$, $\S^1(\body)$, $\M^2(\body)$ and $\Y^2(\body)$ formulations indeed exhibit square-root convergence, it is clear that the relative error is being reduced through the addition of higher order jumping identity fields. The latter can be observed in \cref{fig:ex1}, where we depict the solution on a coarse mesh with $144$ elements.
The higher the order of the discontinuous identity field, the better the approximation compensates for the jump. 
This is particularly apparent for the $\S^1(\body)$-solution in \cref{fig:j}, where the magnitude of the field, as well as its orientation almost perfectly match that of \cref{fig:a}, which represents the exact solution via $\Dp^1(\body)$. 

\subsection{Dilatation in the relaxed micromorphic model}
The importance of $\HsC{,\body}$-conforming finite elements for the relaxed micromorphic model relates to its application with metamaterials \cite{ALBERDI2021104540,RizziMeta,dAgostino2020Meta,Madeo1,Madeo2,Demore}. For example, let a body $\body = \body_i \cup \body_o$ be composed of two domains with distinct complex materials, such that the micro-cell of each material is endowed with a different bulk modulus, then a dilatation field in the domain can jump between the two materials. In linear elasticity theory, this implies differing volumetric changes between material points due to the linearised relation $\dd \body = \det(\one + \D \vb{u}) \dd \body_0 \approx (1 + \tr \D \vb{u}) \dd \body_0$ between the current and reference configurations. Clearly, for this case, the capacity of the discrete subspace to correctly capture jumping identities becomes paramount. 

\begin{figure}
		\centering
		\definecolor{xfqqff}{rgb}{0.4980392156862745,0.,1.}
\definecolor{qqwwzz}{rgb}{0,0.4,0.6}
\definecolor{wqwqwq}{rgb}{0.3764705882352941,0.3764705882352941,0.3764705882352941}
\begin{tikzpicture}[line cap=round,line join=round,>=triangle 45,x=1cm,y=1cm,scale=1.2]
\clip(-4,-2.75) rectangle (4,2.5);
\fill[line width=0.7pt,color=wqwqwq,fill=wqwqwq,fill opacity=0.1] (-2,2) -- (2,2) -- (2,-2) -- (-2,-2) -- cycle;
\fill[line width=0.7pt,color=qqwwzz,fill=qqwwzz,fill opacity=0.1] (-1,1) -- (1,1) -- (1,-1) -- (-1,-1) -- cycle;
\draw [line width=0.7pt,color=wqwqwq] (-2,2)-- (2,2);
\draw [line width=0.7pt,color=wqwqwq] (2,2)-- (2,-2);
\draw [line width=0.7pt,color=wqwqwq] (2,-2)-- (-2,-2);
\draw [line width=0.7pt,color=wqwqwq] (-2,-2)-- (-2,2);
\draw [line width=0.7pt,color=qqwwzz] (-1,1)-- (1,1);
\draw [line width=0.7pt,color=qqwwzz] (1,1)-- (1,-1);
\draw [line width=0.7pt,color=qqwwzz] (1,-1)-- (-1,-1);
\draw [line width=0.7pt,color=qqwwzz] (-1,-1)-- (-1,1);
\draw [-to,line width=0.7pt,color=xfqqff] (0,2) -- (0,2.5);
\draw [-to,line width=0.7pt,color=xfqqff] (2,0) -- (2.5,0);
\draw [-to,line width=0.7pt,color=xfqqff] (0,-2) -- (0,-2.5);
\draw [-to,line width=0.7pt,color=xfqqff] (-2,0) -- (-2.5,0);
\draw [-to,line width=0.7pt,color=xfqqff] (-2,2) -- (-2.5,2.5);
\draw [-to,line width=0.7pt,color=xfqqff] (-2,-2) -- (-2.5,-2.5);
\draw [-to,line width=0.7pt,color=xfqqff] (2,-2) -- (2.5,-2.5);
\draw [-to,line width=0.7pt,color=xfqqff] (2,2) -- (2.5,2.5);
\draw [-to,line width=0.7pt] (0,0) -- (0.5,0);
\draw [-to,line width=0.7pt] (0,0) -- (0,0.5);
\draw [color=qqwwzz](-0.75,-0.75) node[anchor=south west] {$V_{i}$};
\draw [color=wqwqwq](-1.75,-1.75) node[anchor=south west] {$V_{o}$};
\draw [color=xfqqff](2,2) node[anchor=west] {$\mathbf{u}(\mathbf{x})$};
\draw [color=xfqqff](2,0) node[anchor=north west] {$\mathbf{u}(\mathbf{x})$};
\draw [color=xfqqff](1.8,-2) node[anchor=north] {$\mathbf{u}(\mathbf{x})$};
\draw [color=xfqqff](0,-2) node[anchor=north east] {$\mathbf{u}(\mathbf{x})$};
\draw [color=xfqqff](-2.05,-2) node[anchor=east] {$\mathbf{u}(\mathbf{x})$};
\draw [color=xfqqff](-2,0) node[anchor=south east] {$\mathbf{u}(\mathbf{x})$};
\draw [color=xfqqff](-1.8,2) node[anchor=south] {$\mathbf{u}(\mathbf{x})$};
\draw [color=xfqqff](0,2) node[anchor=south west] {$\mathbf{u}(\mathbf{x})$};
\draw (0.5,0) node[anchor=west] {$y$};
\draw (0,0.5) node[anchor=south] {$z$};
\begin{scriptsize}
\draw [fill=black] (0,0) circle (2.5pt);
\end{scriptsize}
\end{tikzpicture}
		\caption{Depiction of a total domain composed of an inner domain $\body_i$ and an outer domain $\body_o$, which are respectively equipped with different material parameters. A dilatation field $\vb{u}(\vb{x})$ is applied on the outer boundary $\partial \body_o$ .}
		\label{fig:dom2}
	\end{figure}
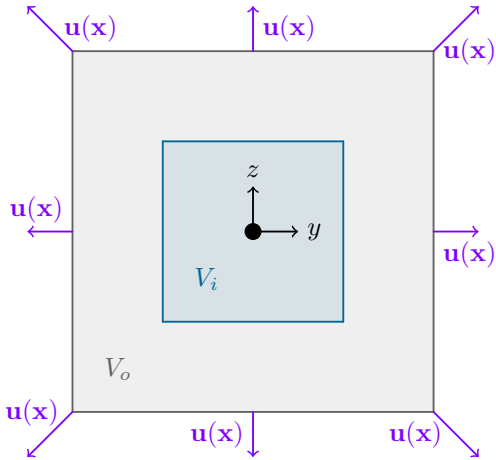

The following example demonstrates exactly this phenomenon, where we compare between formulations based on N\'ed\'elec elements $\{\vb{u},\Pm\} \in [\U^1(\body)]^3 \times [\Ned^0(\body)]^3$ and $\{\vb{u},\Pm\} \in [\U^2(\body)]^3 \times [\Nedtwo^1(\body)]^3$, and our newly proposed elements $\{\vb{u},\Pm\} \in [\U^1(\body)]^3 \times \S^0(\body)$ and $\{\vb{u},\Pm\} \in [\U^2(\body)]^3 \times \S^1(\body)$. We define the cubic domain $\overline{\body} = [-2,2]^3$, which is composed of the inner domain $\body_i = (-1,1)^3$ and the outer domain $\body_o = \body \setminus \body_i$, see \cref{fig:dom2}. 
We set the global material parameters $\muma = 76.9$, $\muc = 0$ and $\Lc = 1$, and 
equip each domain with its own local material parameters
\begin{align}
    &\lamma^i = 115.4 \, , && \lammi^i = 10 \lamma^i \, , && \mumi^i = 10 \muma \, , \notag \\
    &\lamma^o = \lamma^i / 10 \, , && \lammi^o = 100 \lamma^o \, , && \mumi^o = 10 \muma \, .
\end{align}
Using the homogenisation formula \cite{Neff2019}
\begin{align}
    \mue = \dfrac{\mumi\, \muma}{\mumi- \muma} \, , && 2\mue + 3 \lambda_\mathrm{e} = \dfrac{(2\mumi + 3\lambda_\mathrm{micro})(2\muma + 3\lambda_\mathrm{macro})}{(2\mumi + 3\lambda_\mathrm{micro})-(2\muma + 3\lambda_\mathrm{macro})} \, ,
\end{align}
we derive the corresponding meso material parameters
\begin{align}
    &\mue^i = \mue^o = 85.44 \, , && \lame^i = 128.22 \, , && \lame^o = 8.3 \, . 
\end{align}
Clearly, the materials differ in their Lam\'e constants $\lame$ and $\lammi$, which govern the trace of a strain measure. In contrast, their shear moduli $\mue$ and $\mumi$ are the same.
On the entire boundary of the outer domain we impose an expansion displacement field (see \cref{fig:dom2})
\begin{align}
    &\vb{u}(\vb{x}) = \dfrac{1}{10} \begin{bmatrix}
        x \\ y \\ z
    \end{bmatrix} \, , && \text{on} && \surf_D^u = \partial \body_o \, , 
\end{align}
and correspondingly, the consistent coupling condition \cite{dagostino2021consistent}
\begin{align}
    &\tr_{\HsC{}} \Pm = \tr_{\HsC{}}\D \vb{u} = \dfrac{1}{10}\tr_{\HsC{}}  \one = 0  && \text{on} && \surf_D^P = \partial \body_o \, .
\end{align}
Although $\one \in \ker(\tr_{\HsC{}})$, our element is not perfectly conforming in $\HsC{,\body}$ with respect to minimal regularity, such that the tangential projection of the constant identify field is imposed. Note that for the N\'ed\'elec formulation we employ the trace of the $\HC{,\body}$-space $\tr_{\HC{}} \Pm = (1/10) \tr_{\HC{}} \one$.  
We compare the elastic energy produced by the dilatation field. The energy $I(\vb{u},\Pm)$ is taken from \cref{eq:rmm} by removing the external work via $\vb{f}$ and $\bm{M}$, which is zero in this example. 

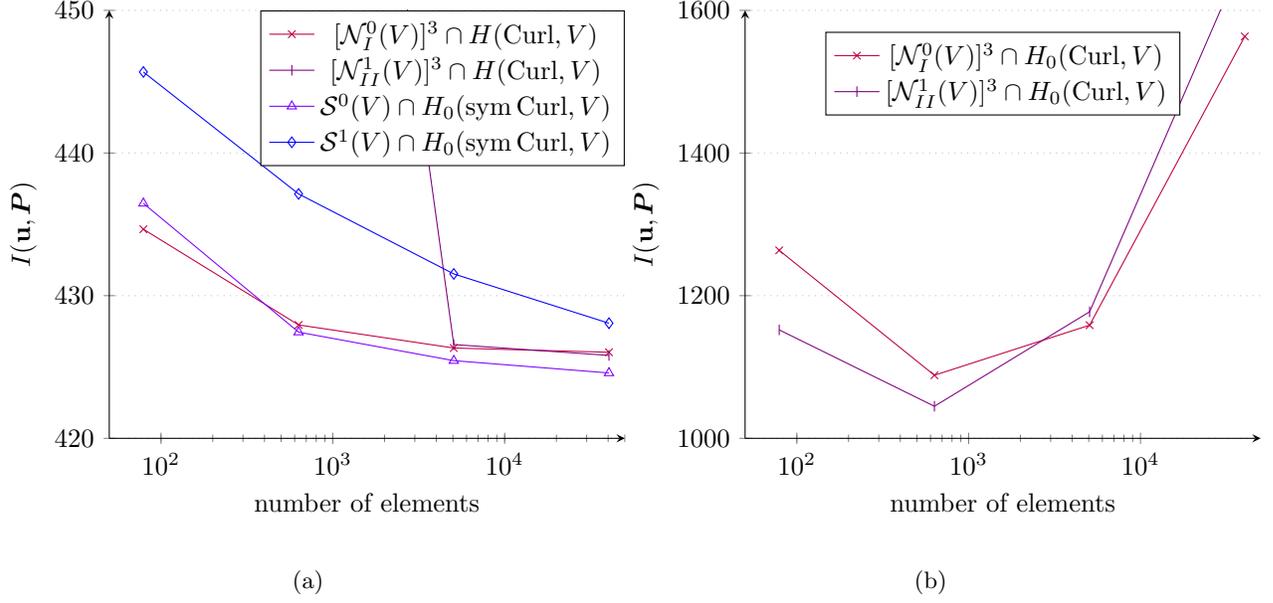
\begin{figure}
    	\centering
     \begin{subfigure}{0.48\linewidth}
    		\centering
    		\begin{tikzpicture}
    			\definecolor{asl}{rgb}{0.4980392156862745,0.,1.}
    			\definecolor{asb}{rgb}{0.,0.4,0.6}
    			\begin{semilogxaxis}[
    				/pgf/number format/1000 sep={},
    				axis lines = left,
    				xlabel={number of elements},
    				ylabel={$I(\vb{u},\Pm)$},
    				xmin=50, xmax=50000,
    				ymin=420, ymax=450,
    				xtick={1e2,1e3,1e4,1e5},
    				ytick={420, 430, 440, 450},
    				legend style={at={(1,1)},anchor= north east},
    				ymajorgrids=true,
    				grid style=dotted,
    				]
    				\addplot[color=purple, mark=x] coordinates {
    					( 79 , 434.65861576292957 )
                        ( 632 , 427.9492129667167 )
                        ( 5056 , 426.33095072985196 )
                        ( 40448 , 426.0324459273685 )
    				};
    				\addlegendentry{$[\Ned^0(\body)]^3 \cap H(\Curl,\body)$}

                    \addplot[color=violet, mark=|] coordinates {
                        ( 79 , 1387.343378339081 )
    					( 632 , 505.2248086867979 )
                        ( 5056 , 426.57200771073667 )
                        ( 40448 , 425.8126659822291 )
    				};
    				\addlegendentry{$[\Nedtwo^1(\body)]^3 \cap H(\Curl,\body)$}

            		\addplot[color=asl, mark=triangle] coordinates {
    					( 79 , 436.47237918528583 )
                        ( 632 , 427.44315348221573 )
                        ( 5056 , 425.4447055010866 )
                        ( 40448 , 424.582055891771 )
    				};
    				\addlegendentry{$\S^0(\body) \cap H_0(\sym\Curl,\body)$}

                    \addplot[color=blue, mark=diamond] coordinates {
                        ( 79 , 445.6807580202008 )
    					( 632 , 437.13179668562844 )
                        ( 5056 , 431.52392885882546 )
                        ( 40448 , 428.06900077583026 )
    				};
    				\addlegendentry{$\S^1(\body) \cap H_0(\sym\Curl,\body)$}
    				   				
    			\end{semilogxaxis}
    		\end{tikzpicture}
    		\caption{}
    	\end{subfigure}
     \begin{subfigure}{0.48\linewidth}
    		\centering
    		\begin{tikzpicture}
    			\definecolor{asl}{rgb}{0.4980392156862745,0.,1.}
    			\definecolor{asb}{rgb}{0.,0.4,0.6}
    			\begin{semilogxaxis}[
    				/pgf/number format/1000 sep={},
    				axis lines = left,
    				xlabel={number of elements},
    				ylabel={$I(\vb{u},\Pm)$},
    				xmin=50, xmax=50000,
    				ymin=1000, ymax=1600,
    				xtick={1e2,1e3,1e4},
    				ytick={1000, 1200, 1400, 1600},
    				legend style={at={(0.5,0.95)},anchor= north},
    				ymajorgrids=true,
    				grid style=dotted,
    				]
    				\addplot[color=purple, mark=x] coordinates {
    					( 79 , 1263.5674962115856 )
                        ( 632 , 1088.443006247015 )
                        ( 5056 , 1158.629869929022 )
                        ( 40448 , 1563.5013024952202 )
    				};
    				\addlegendentry{$[\Ned^0(\body)]^3 \cap H_0(\Curl,\body)$}
        
                    \addplot[color=violet, mark=|] coordinates {
                        ( 79 , 1152.0830412837704 )
                        ( 632 , 1044.9471767948949 )
                        ( 5056 , 1177.4300279042466 )
                        ( 40448 , 1678.8350746504289 )
    				};
    				\addlegendentry{$[\Nedtwo^1(\body)]^3 \cap H_0(\Curl,\body)$}    	
        
    			\end{semilogxaxis}
    		\end{tikzpicture}
    		\caption{}
    	\end{subfigure}
    	\caption{Elastic energy due to dilatation. The N\'ed\'elec based formulation is computed either with (a) a natural Neumann boundary   $H(\Curl, \body)$ or (b) complete Dirichlet boundary conditions $H_0(\Curl, \body)$.}
    	\label{fig:energy}
    \end{figure}

As seen in \cref{fig:energy}, for the N\'ed\'elec based formulations to produce similar energies to our novel elements, one must free the microdistortion on the boundary completely $\surf_N^P = \partial \body_o$, which is incompatible with the consistent coupling condition.
Even in latter case, the energies differ slightly $\Delta I \approx 0.4 \%$ due to the jump at the interface, such that the $\S^0(\body)$-formulation produces lower energies. The relative difference is given here with respect to the lower energy of the $\S^0(\body)$-formulation. 
In contrast, if a total Dirichlet boundary $\surf_D^P = \partial \body_o$ is enforced in $\Pm$, then the elastic energy in the medium increases significantly and we observe no convergence towards a specific energy value. We note that the same occurs even if a cubic formulation $\{\vb{u}, \Pm\} \in [\U^3(\body)]^3 \times [\Ned^2(\body)]^3$ is applied. The latter hints at the emergence of a boundary layer problem.      
    
\begin{figure}
    	\centering
     \begin{subfigure}{0.48\linewidth}
    		\centering
    		\includegraphics[width=0.8\linewidth]{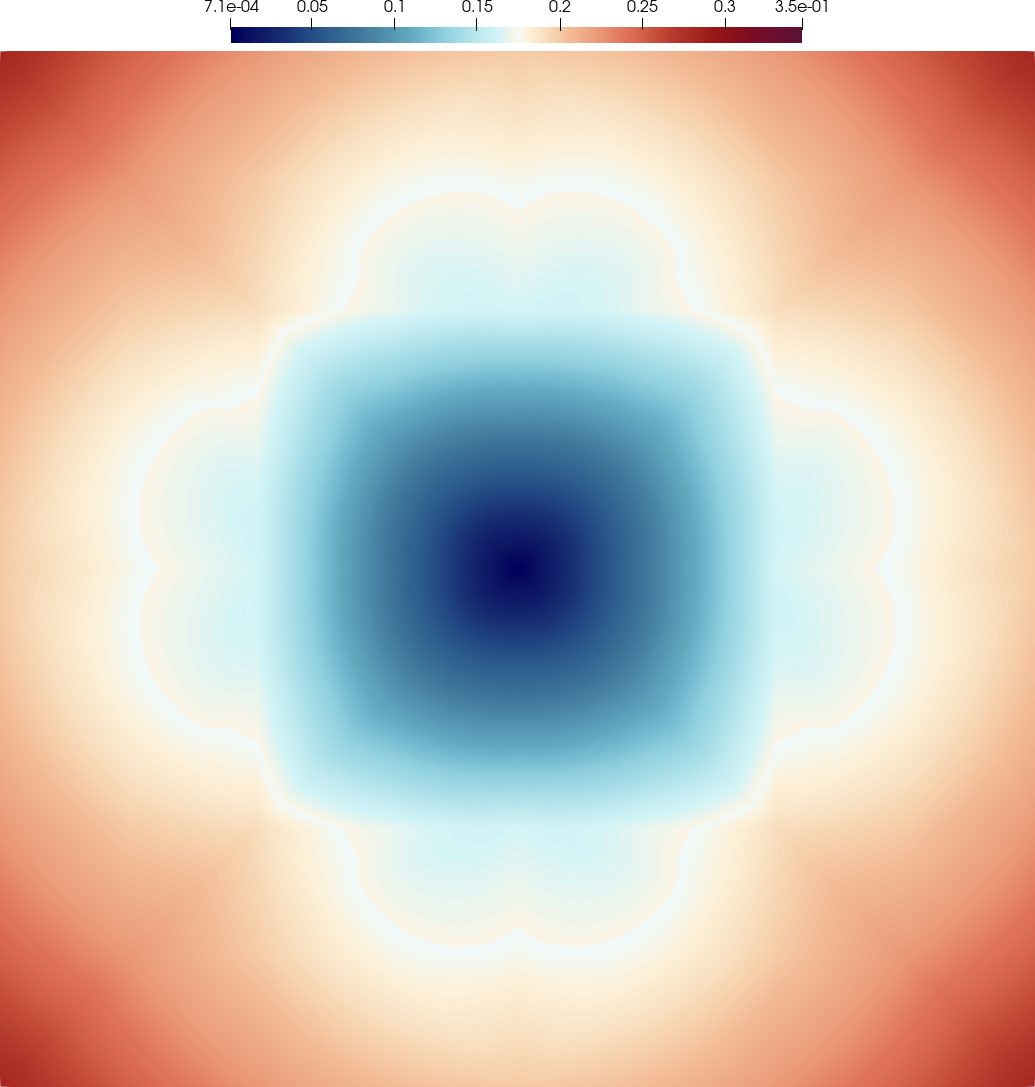}
    		\caption{}
    	\end{subfigure}
     \begin{subfigure}{0.48\linewidth}
    		\centering
    		\includegraphics[width=0.8\linewidth]{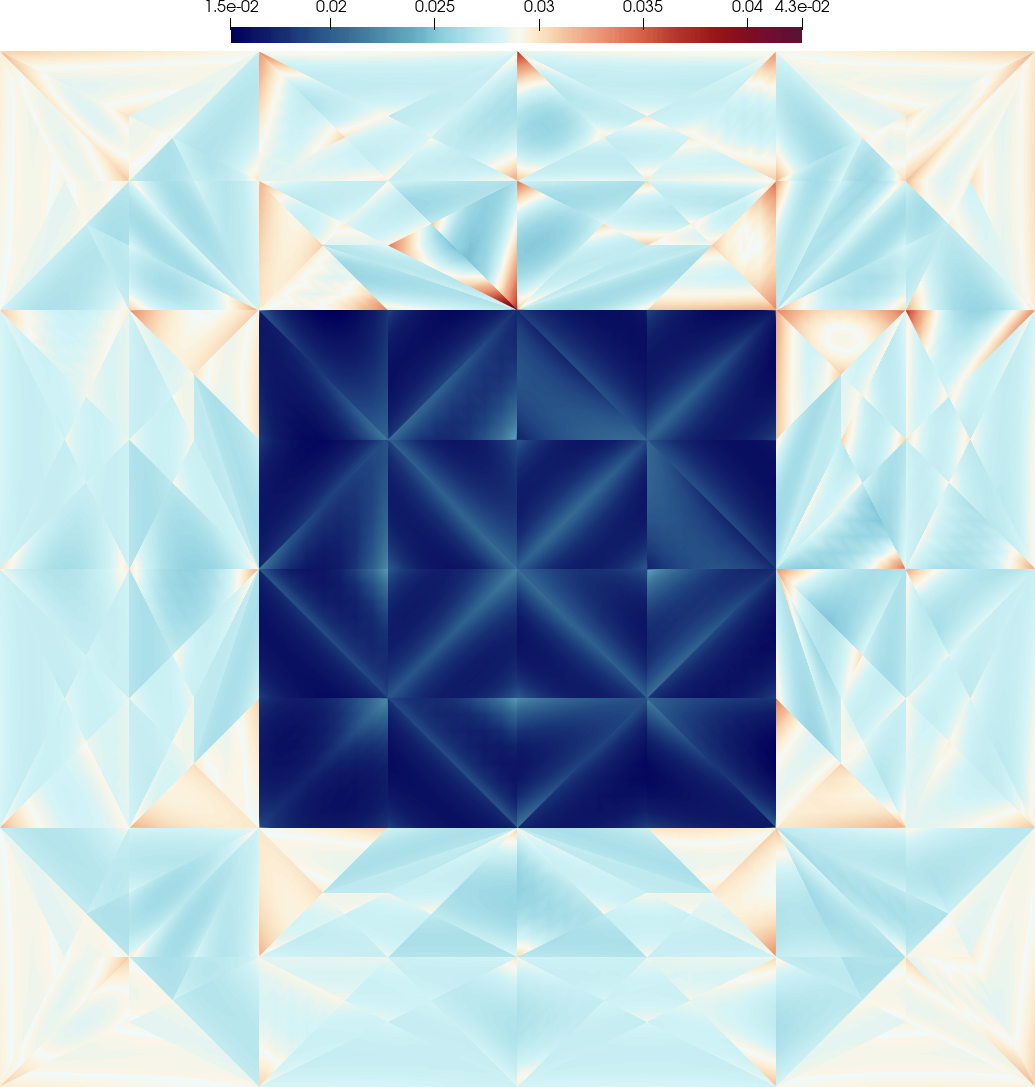}
    		\caption{}
    	\end{subfigure}
     \begin{subfigure}{0.48\linewidth}
    		\centering
    		\includegraphics[width=0.8\linewidth]{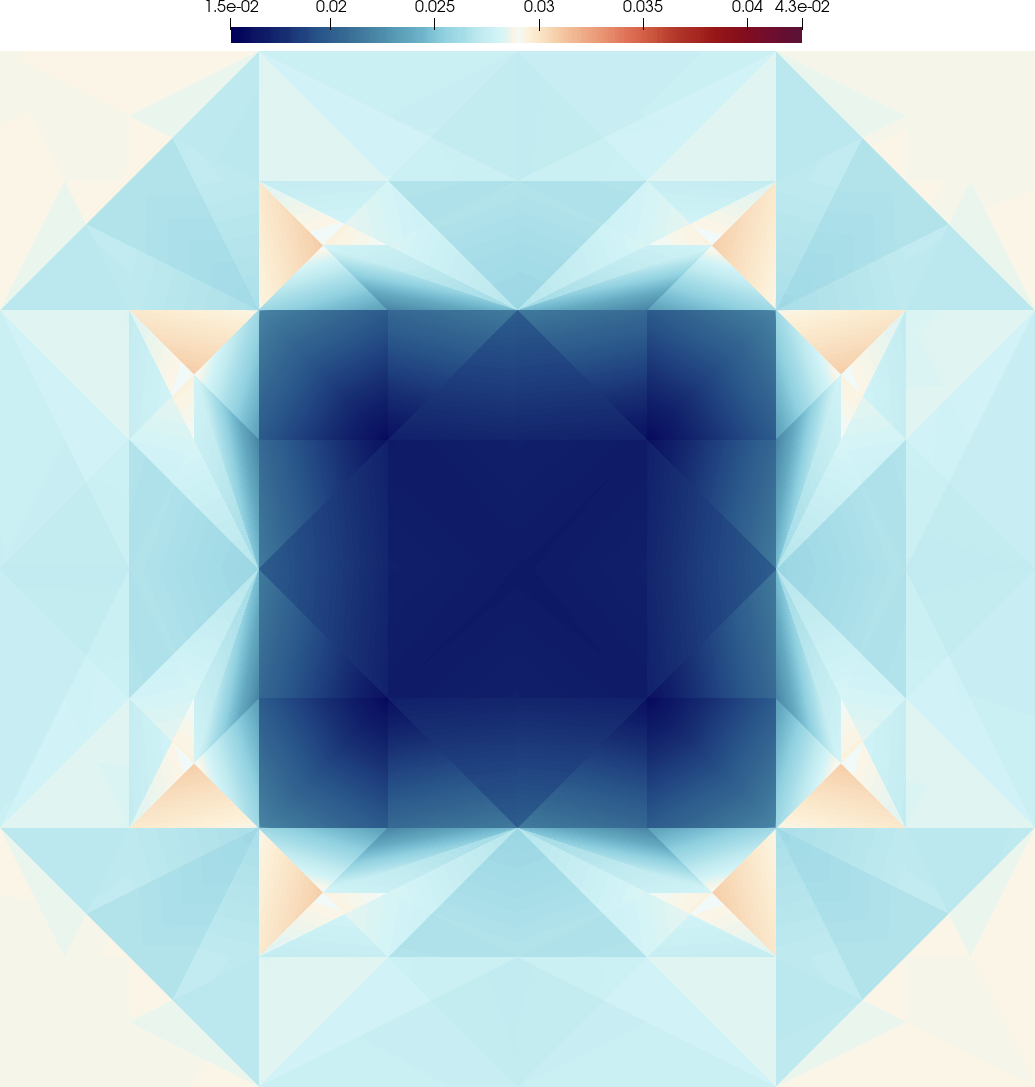}
    		\caption{}
    	\end{subfigure}
     \begin{subfigure}{0.48\linewidth}
    		\centering
    		\includegraphics[width=0.8\linewidth]{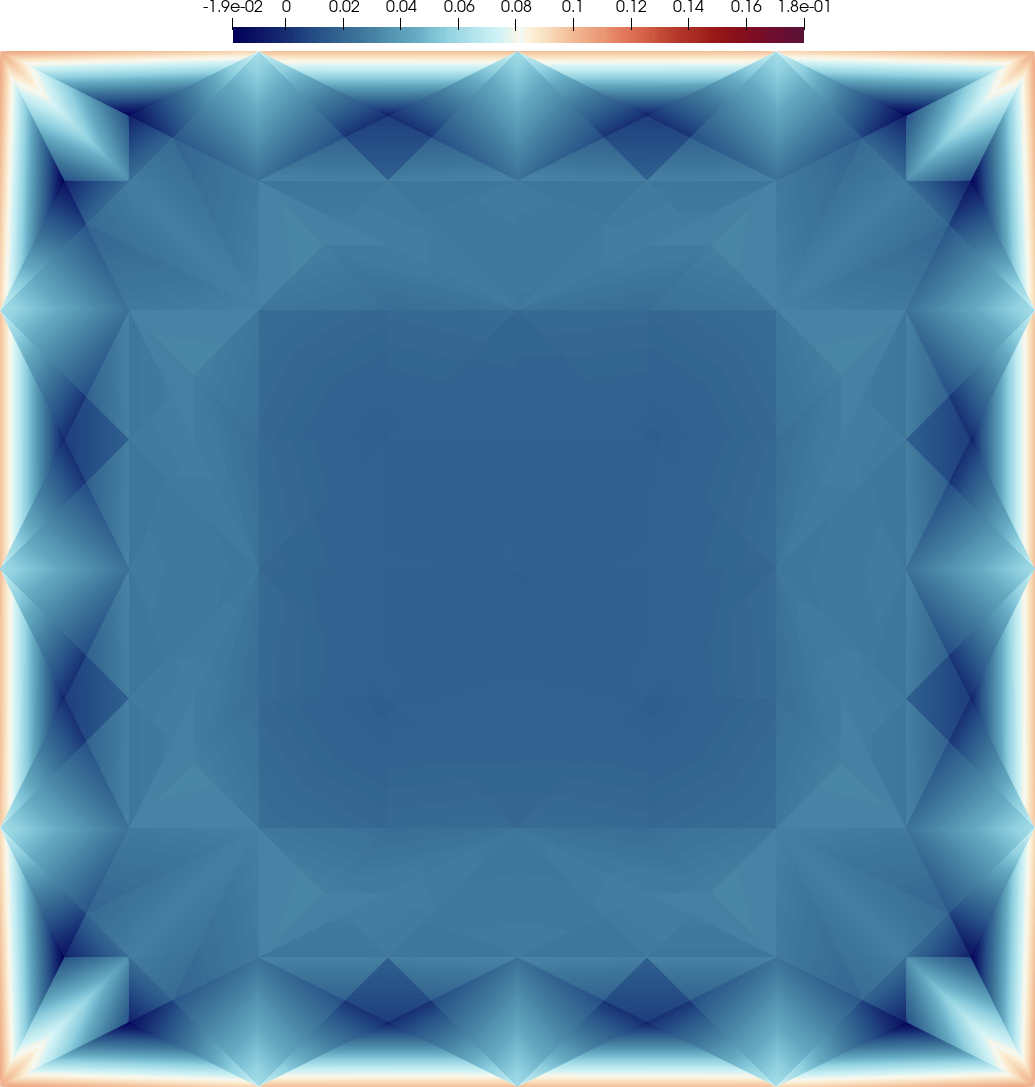}
    		\caption{}
    	\end{subfigure}
    	\caption{Resulting displacement field $\vb{u}$ due to dilatation (a). Trace of the microdistortion field $\tr\Pm$ using the novel $\S^0(\body)$-formulation (b), the N\'ed\'elec formulation with a natural Neumann boundary (c) and with a total Dirichlet boundary (d).}
    	\label{fig:ex2figs}
    \end{figure}

The resulting displacement field $\vb{u}$ as well as the trace of the microdistortion field $\tr\Pm$ are depicted in \cref{fig:ex2figs} for a mesh with $5056$ elements. 
Although the N\'ed\'elec based formulation with a natural Neumann boundary produces similar energies, it is clearly incapable of satisfactorily capturing the jumping trace of the microdistortion field. This is evident due to the visible transition in the intensity of the trace between the two materials.   
In contrast, our novel $\S^0(\body)$-formulation clearly distinguishes between the trace of the microdistortion in the outer and inner materials. Lastly, from the depiction it is obvious that imposing the consistent coupling condition on the N\'ed\'elec based formulation leads to a completely different state in the trace of the microdistortion. 

\section{Hexahedra on transfinite meshes}
\begin{figure}
    		\centering
    		\definecolor{asl}{rgb}{0.4980392156862745,0.,1.}
\definecolor{zzttqq}{rgb}{0.,0.4,0.6}
\definecolor{ududff}{rgb}{0.,0.4,0.6}
\begin{tikzpicture}[line cap=round,line join=round,>=triangle 45,x=1cm,y=1cm,scale = 0.5]
\clip(2.5,1.5) rectangle (29.5,12);
\fill[line width=0.7pt,color=zzttqq,fill=zzttqq,fill opacity=0.10000000149011612] (4,4) -- (7,3) -- (7,6) -- (4,7) -- cycle;
\fill[line width=0.7pt,color=zzttqq,fill=zzttqq,fill opacity=0.10000000149011612] (7,3) -- (9,5) -- (9,8) -- (7,6) -- cycle;
\fill[line width=0.7pt,color=zzttqq,fill=zzttqq,fill opacity=0.10000000149011612] (4,7) -- (6,9) -- (9,8) -- (7,6) -- cycle;
\fill[line width=0.7pt,dashed,color=zzttqq,fill=zzttqq,fill opacity=0.10000000149011612] (16,6) -- (25,3) -- (25,5) -- (16,8) -- cycle;
\fill[line width=0.7pt,dashed,color=zzttqq,fill=zzttqq,fill opacity=0.10000000149011612] (25,3) -- (29,7) -- (29,9) -- (25,5) -- cycle;
\fill[line width=0.7pt,dashed,color=zzttqq,fill=zzttqq,fill opacity=0.10000000149011612] (16,8) -- (20,12) -- (29,9) -- (25,5) -- cycle;
\draw [line width=0.7pt,color=zzttqq] (4,4)-- (7,3);
\draw [line width=0.7pt,color=zzttqq] (7,3)-- (7,6);
\draw [line width=0.7pt,color=zzttqq] (7,6)-- (4,7);
\draw [line width=0.7pt,color=zzttqq] (4,7)-- (4,4);
\draw [line width=0.7pt,color=zzttqq] (7,3)-- (9,5);
\draw [line width=0.7pt,color=zzttqq] (9,5)-- (9,8);
\draw [line width=0.7pt,color=zzttqq] (9,8)-- (7,6);
\draw [line width=0.7pt,color=zzttqq] (7,6)-- (7,3);
\draw [line width=0.7pt,color=zzttqq] (4,7)-- (6,9);
\draw [line width=0.7pt,color=zzttqq] (6,9)-- (9,8);
\draw [line width=0.7pt,color=zzttqq] (9,8)-- (7,6);
\draw [line width=0.7pt,color=zzttqq] (7,6)-- (4,7);
\draw [line width=0.7pt,dashed] (6,6)-- (4,4);
\draw [line width=0.7pt,dashed] (6,6)-- (9,5);
\draw [line width=0.7pt,dashed] (6,6)-- (6,9);
\draw [line width=0.7pt,dashed,color=zzttqq] (16,6)-- (25,3);
\draw [line width=0.7pt,dashed,color=zzttqq] (25,3)-- (25,5);
\draw [line width=0.7pt,dashed,color=zzttqq] (25,5)-- (16,8);
\draw [line width=0.7pt,dashed,color=zzttqq] (16,8)-- (16,6);
\draw [line width=0.7pt,dashed,color=zzttqq] (25,3)-- (29,7);
\draw [line width=0.7pt,dashed,color=zzttqq] (29,7)-- (29,9);
\draw [line width=0.7pt,dashed,color=zzttqq] (29,9)-- (25,5);
\draw [line width=0.7pt,dashed,color=zzttqq] (25,5)-- (25,3);
\draw [line width=0.7pt,dashed,color=zzttqq] (16,8)-- (20,12);
\draw [line width=0.7pt,dashed,color=zzttqq] (20,12)-- (29,9);
\draw [line width=0.7pt,dashed,color=zzttqq] (29,9)-- (25,5);
\draw [line width=0.7pt,dashed,color=zzttqq] (25,5)-- (16,8);
\draw [line width=0.7pt,dashed] (16,6)-- (20,10);
\draw [line width=0.7pt,dashed] (20,10)-- (20,12);
\draw [line width=0.7pt,dashed] (18,10)-- (18,8);
\draw [line width=0.7pt,dashed] (19,5)-- (23,9);
\draw [line width=0.7pt,dashed] (20,10)-- (29,7);
\draw [line width=0.7pt,dashed] (18,8)-- (27,5);
\draw [line width=0.7pt,dashed] (22,4)-- (26,8);
\draw [line width=0.7pt,dashed] (18,10)-- (27,7);
\draw [line width=0.7pt,dashed] (19,7)-- (23,11);
\draw [line width=0.7pt,dashed] (22,6)-- (26,10);
\draw [line width=0.7pt,dashed] (21,9)-- (21,7);
\draw [line width=0.7pt,dashed] (24,6)-- (24,8);
\draw [line width=0.7pt,dashed] (27,7)-- (27,5);
\draw [line width=0.7pt,dashed] (19,7)-- (19,5);
\draw [line width=0.7pt,dashed] (22,6)-- (22,4);
\draw [-to,line width=0.7pt] (7,3) -- (8.5,2.5);
\draw [-to,line width=0.7pt] (4,7) -- (4,9);
\draw [-to,line width=0.7pt] (6,6) -- (7,7);
\draw [line width=0.7pt,dashed] (23,11)-- (23,9);
\draw [line width=0.7pt,dashed] (26,10)-- (26,8);
\draw (6.1,3.1) node[anchor=north west] {$_{v_2}$};
\draw (9.1,5.391108562609959) node[anchor=north west] {$_{v_3}$};
\draw (7.25,6.45) node[anchor=north west] {$_{v_6}$};
\draw (3.,4.2) node[anchor=north west] {$_{v_1}$};
\draw (9.1,8.446614602664864) node[anchor=north west] {$_{v_7}$};
\draw (5.25,9.9) node[anchor=north west] {$_{v_8}$};
\draw (8.65,2.9) node[anchor=north west] {$\xi$};
\draw (3.5,10.15) node[anchor=north west] {$\zeta$};
\draw (3.,7.4) node[anchor=north west] {$_{v_5}$};
\draw [-to,line width=0.7pt] (11,6) -- (14,6);
\draw (10.6,7.2) node[anchor=north west] {$V_t = \cup_e \elem_e$};
\draw (5.55,5.3) node[anchor=north west,color=zzttqq] {$\Omega$};

\draw (17,7) node[anchor=north west,color=zzttqq] {$\elem_1$};
\draw (20,6) node[anchor=north west,color=zzttqq] {$\elem_2$};
\draw (23,5) node[anchor=north west,color=zzttqq] {$\elem_3$};
\draw (17+3.5,7+4) node[anchor=north west,color=zzttqq] {$\elem_4$};
\draw (20+3.5,6+4) node[anchor=north west,color=zzttqq] {$\elem_5$};
\draw (23+3.5,5+4) node[anchor=north west,color=zzttqq] {$\elem_6$};
\begin{scriptsize}
\draw [fill=ududff] (4,4) circle (2.5pt);
\draw [fill=ududff] (7,3) circle (2.5pt);
\draw [fill=ududff] (6,6) circle (2.5pt);
\draw [fill=ududff] (9,5) circle (2.5pt);
\draw [fill=ududff] (4,7) circle (2.5pt);
\draw [fill=ududff] (7,6) circle (2.5pt);
\draw [fill=ududff] (9,8) circle (2.5pt);
\draw [fill=ududff] (6,9) circle (2.5pt);
\end{scriptsize}
\end{tikzpicture}
	\caption{Transfinite mesh $\body_t$ of cuboids built by translating and resizing the reference hexahedron $\Omega$.}
	\label{fig:trans}
\end{figure}
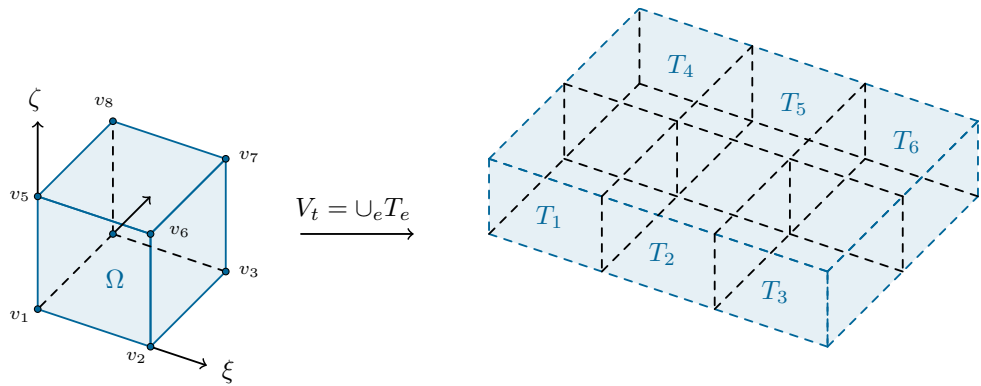

The elements we presented on tetrahedralizations do not achieve the minimal regularity of the continuous $\HsC{,\body}$-space in the lowest order. As clarified, the problem is related to the inability to relate a single coordinate system on a vertex to multiple edges on interfacing tetrahedra. This problem can be ameliorated if a transfinite mesh of cubes is used in the construction, see \cref{fig:trans}. 
Clearly, this limits the possible geometries of the domain significantly and is thus only applicable to certain cases. As such, we only provide the theory for the construction in this work, without examples.

\subsection{The lowest order element}
We start by investigating polynomial spaces on the unit hexahedron $\Omega = [0,1]^3$, see \cref{fig:trans}. For $\Hone(\body)$-conforming elements one employs the polynomial space
\begin{align}
    &\Q^{p.p.p}(\Omega) = \Po^p(\xi) \times \Po^p(\eta) \times \Po^p(\zeta) \, , && \dim \Q^{p.p.p}(\Omega) = (p+1)^3 \, . 
\end{align}
Its gradient space is given by 
\begin{align}
    \nabla\Q^{p.p.p}(\Omega) \subseteq  \Q^{p-1.p.p}(\Omega) \times \Q^{p.p-1.p}(\Omega) \times \Q^{p.p.p-1}(\Omega)  \, . 
\end{align}
Now, constructing its matrix-valued version and taking the deviator yields
\begin{align}
    \dev \D [\R^3 \otimes \Q^{p.p.p}(\Omega)] &\subseteq \dev [\R^3 \otimes [\Q^{p-1.p.p}(\Omega) \times \Q^{p.p-1.p}(\Omega) \times \Q^{p.p.p-1}(\Omega)] ] \notag \\
    &\subseteq \R^3 \otimes [\Q^{p-1.p.p}(\Omega) \times \Q^{p.p-1.p}(\Omega) \times \Q^{p.p.p-1}(\Omega)] \notag \\  &\quad- \dfrac{1}{3} [\Q^{p-1.p.p}(\Omega) \cup  \Q^{p.p-1.p}(\Omega) \cup  \Q^{p.p.p-1}(\Omega)] \one \notag \\
    &\subseteq \begin{bmatrix}
        \Q^{p.p.p}(\Omega) & \Q^{p.p-1.p}(\Omega) & \Q^{p.p.p-1}(\Omega) \\
        \Q^{p-1.p.p}(\Omega) & \Q^{p.p.p}(\Omega) & \Q^{p.p.p-1}(\Omega) \\
        \Q^{p-1.p.p}(\Omega) & \Q^{p.p-1.p}(\Omega) & \Q^{p.p.p}(\Omega) 
    \end{bmatrix} = \S^{p-1}(\Omega)  \, .
\end{align}
Clearly, the same result follows by the combination
\begin{align}
    &\S^{p-1}(\Omega) \supseteq \D [\R^3 \otimes \Q^{p.p.p}(\Omega)] \cup [\Q^{p.p.p}(\Omega) \otimes \one] \, , && \dim \S^{p-1}(\Omega) = 3(p+1)^3+ 6p(p+1)^2  \, .
+\end{align}
Now, since $\Q^{p.p.p}(\Omega) \otimes \one \in \ker(\sym \Curl)$, the next space in the sequence reads
\begin{align}
    \sym \Curl[\Q^{p-1.p.p}(\Omega) \times \Q^{p.p-1.p}(\Omega) \times \Q^{p.p.p-1}(\Omega)]^3 &\subseteq \sym [\Q^{p.p-1.p-1}(\Omega) \times \Q^{p-1.p.p-1}(\Omega) \times \Q^{p-1.p-1.p}(\Omega)]^3 \notag \\
    &\subseteq \begin{bmatrix}
        \Q^{p.p-1.p-1}(\Omega) & \Q^{p.p.p-1}(\Omega) & \Q^{p.p-1.p}(\Omega) \\
        \Q^{p.p.p-1}(\Omega) & \Q^{p-1.p.p-1}(\Omega) & \Q^{p-1.p.p}(\Omega) \\
        \Q^{p.p-1.p}(\Omega) & \Q^{p-1.p.p}(\Omega) & \Q^{p-1.p-1.p}(\Omega) 
    \end{bmatrix} 
    \notag \\
    &= \mathcal{W}^{p-1}(\Omega) \, .
\end{align}
Lastly, the last space in the sequence is given by
\begin{align}
    \di \Di \mathcal{W}^{p-1}(\Omega) = \Q^{p-1.p-1.p-1}(\Omega) \, .
\end{align}
Consequently, the polynomial sequence in \cref{fig:hexseq} is used to construct $\HsC{,\body}$-conforming finite elements on hexahedra. 
\begin{figure}
    		\centering
    		\begin{tikzpicture}[scale = 0.6][line cap=round,line join=round,>=triangle 45,x=1.0cm,y=1.0cm]
		\clip(6,5) rectangle (29,10);
		\draw (10.5,6.5) node[anchor=north east] {$\Q^{p.p.p}(\Omega) \otimes \one$};
		\draw [line width=1.5pt] (10.5,6) -- (14.5,6);
		\draw [->,line width=1.5pt] (14.5,6) -- (14.5,8);
		\draw (11.5,7) node[anchor=north west] {$\text{id}$};
		\draw (11.5,8) node[anchor=north west] {$\oplus$};

        \draw (10.5,9) node[anchor=north east] {$[\Q^{p.p.p}(\Omega)]^3$};
		\draw [->,line width=1.5pt] (10.5,8.5) -- (13.5,8.5);
        \draw (16,9) node[anchor=north east] {$\S^{p-1}(\Omega)$};
		\draw (11,9.5) node[anchor=north west] {$\dev \D$};
		\draw [->,line width=1.5pt] (17.9-2.4+0.5,8.5) -- (20.9-2.4+0.5,8.5);
		\draw (17.8-2.4+0.5,9.5) node[anchor=north west] {$\sym \Curl$};
		\draw (24.8-3.15,9) node[anchor=north east] {$\mathcal{W}^{p-1}(\Omega)$};
		\draw [->,line width=1.5pt] (24.8-3.2,8.5) -- (27.8-3.2,8.5);
		\draw (25-3.2,9.5) node[anchor=north west] {$\di \Di$};
		\draw (27.8-3.2,8.9) node[anchor=north west] {$\Q^{p-1.p-1.p-1}(\Omega)$};
	\end{tikzpicture}
	\caption{The polynomial relaxed micromorphic sequence on hexahedra.}
	\label{fig:hexseq}
\end{figure} 
With the polynomial space at hand, we can construct $\HsC{,\body}$-conforming finite elements using the polytopal template methodology \cite{sky_polytopal_2022,sky_higher_2023}. The vectorial template set on each edge reads
\begin{align}
    &\tem_{j} = \{\bm{\tau}, \, \bm{\mu}, \, \bm{\nu}\} \, , && j \in \mathcal{J} = \{ (0,1),(1,2),(2,3),(0,3),(0,4),(1,5),(2,6),(3,7),(4,5),(5,6),(6,7),(7,4) \} \, , 
\end{align}
which contains the tangent, cotangent and normal vectors, respectively. In order to maintain the same orientation across the mesh the vectors are simply defined as permutations of the Cartesian basis. For example, on the edge $e_{12}$ (see \cref{fig:trans}) the vectors read $ \{\bm{\tau}, \, \bm{\mu}, \, \bm{\nu}\} = \{\vb{e}_1, \,\vb{e}_2, \,\vb{e}_3\}$. 
Now, using the Lagrangian base functions
\begin{align}
    \lambda_1 &= (1 - \xi)(1-\eta)(1-\zeta) \, , & \lambda_2 &= \xi(1-\eta)(1-\zeta) \, , & \lambda_3 &= \xi\eta(1-\zeta) \, , & \lambda_4 &= (1 - \xi)\eta(1-\zeta) \, , \notag \\
    \lambda_5 &= (1 - \xi)(1-\eta)\zeta \, , & \lambda_6 &= \xi(1-\eta)\zeta \, , & \lambda_7 &= \xi\eta\zeta \, , & \lambda_8 &= (1 - \xi)\eta\zeta \, ,
\end{align}
we define the base functions of the lowest order $\S^0(\Omega)$-element directly.
\begin{definition}[Lowest order hexahedral element]
    The base functions of the lowest order element are defined polytope-wise.
    \begin{itemize}
        \item on every vertex $v_i$ we construct the base functions
        \begin{align}
            \bm{\varrho}(\xi,\eta,\zeta) &= \left \{ \begin{aligned}
                &\lambda_i (\vb{e}_1 \otimes \vb{e}_1-\vb{e}_2 \otimes \vb{e}_2)  \\
                &\lambda_i (\vb{e}_2 \otimes \vb{e}_2-\vb{e}_3 \otimes \vb{e}_3)  
            \end{aligned} \right . \, ,
        \end{align}
        such that each vertex defines two base functions.
        \item on every edge $e_{ij}$ with $(i,j) \in \mathcal{J}$ we construct the base functions
        \begin{align}
            \bm{\varrho}(\xi,\eta,\zeta) &= \left \{ \begin{aligned}
                &(\lambda_i + \lambda_j) \bm{\mu} \otimes \bm{\tau}
                \\
                &(\lambda_i + \lambda_j) \bm{\nu} \otimes \bm{\tau}
            \end{aligned} \right . \, ,  && \{\bm{\tau}, \, \bm{\mu}, \, \bm{\nu}\} \in \tem_{ij} \, ,
        \end{align}
        such that each edge defines four base functions.
        \item the cell base functions are given by
        \begin{align}
            \bm{\varrho}(\xi,\eta,\zeta) &= \lambda_i \one \, ,
        \end{align}
        for every vertex $v_{i}$. 
    \end{itemize}
\end{definition}
There are $16$ vertex base functions, $24$ edge base functions and $8$ cell base functions. Thus, the space is complete due to $\dim \S^0(\Omega) = 48$. The linear independence of the construction is obvious, since the vertex base functions plus the cell base functions span $\Q^{1.1.1}(\Omega) \otimes \Diag(3)$, and the edge base functions span off-diagonal terms. The conformity is also clear, since the vertex base functions impose the $\C^0(\body)$-continuity of $\Diag(3) \cap \Dev(3)$, the edge base functions are built using the tangent vector $(\cdot)\otimes \bm{\tau}$, such that they uphold $\Hc{,\body}$-conformity, and the cell base functions are in the kernel of the $\HsC{}$-trace operator.
\begin{remark}[Extension to higher orders]
    The extension of the element to higher orders follows analogously. 
    One must simply apply the polytopal templates to the corresponding scalar base functions while respecting the dimension of $\S^p(\Omega)$.
\end{remark}

\subsection{Consistent transformations}
The transfinite geometry of the reference grid of cuboids can be used to construct a curved domain with conforming elements. We shortly discuss this case in this section.
\begin{lemma}[Consistent transformation]
    Let the domain $\body$ be such that the transfinite grid $\cup_e \elem_e$ (see \cref{fig:trans}) can be mapped to it by a $\C^1(\body)$-continuous function
\begin{align}
    &F:\bigcup_e \elem_e \to \body \, , && F \in \C^1(\body) \, ,
\end{align}
then the transformation
\begin{align}
    \bm{\rho} = \bm{J} \bm{\varrho} \bm{J}^{-1} \, ,
\end{align}
maintains $\HsC{,\body}$-conformity.
\end{lemma}
\begin{proof}
    Each base function on the reference grid is defined by the dyadic product of two vectors $\vb{a} \otimes \vb{b}$. Clearly, the vectors match on interfacing elements of the reference grid. Now, due to the imposed $\C^1(\body)$-continuity of the mapping, vectors mapped as $\bm{J} \vb{a}$ continue to match on interfacing elements. Further, the covariant Piola transformation $\vb{b} \bm{J}^{-1}$ maintains tangential conformity. Lastly, deviatoric fields remain such under the transformation due to
    $
        \dev (\bm{J} \one \bm{J}^{-1} ) = \dev(\one) = 0 \, .
    $
    This completes the proof.
\end{proof}

\section{Conclusions and outlook}
In this work we presented the relaxed micromorphic sequence as the completion of the $\di \Di$-sequence in the $\HsC{,\body}$-space with respect to the kernel of the $\sym \Curl$-operator. The sequence is associated with the relaxed micromorphic model, in which the microdistortion field $\Pm$ is defined in the full $\HsC{,\body}$-space $\Pm \in \HsC{,\body}$. The consistency requirements of the microdistortion field in the relaxed micromorphic model motivate the introduction of novel $\HsC{,\body}$-conforming finite elements, which also respect tangential $\HC{,\body}$-conformity. To the authors knowledge, this work introduces the first of such elements. The novel elements presented here are based on extensions of the low order N\'ed\'elec elements. Our rigorous proofs along with the numerical examples demonstrate the correct conformity of the formulations in $\HC{,\body}$. Further, the investigation of $\HsC{,\body}$-conforming fields shows that the new lowest order elements are optimal, and that the higher order elements can greatly reduce the error induced by their low order non-jumping identity fields. The elements are given via closed formulas for the base functions, which also allow to directly map them to curved geometries. Latter is a topic for future work. Finally, we mention that the elements are also conforming in $\HdsC{,\body}$, since the spaces possess the same regularity. In fact, for complete polynomial spaces $[\Po^p(\body)]^{3 \times 3}$ the elements coincide. The question whether it is possible to build differing $\HdsC{,\body}$-conforming finite element using incomplete polynomial spaces remains open, to be addressed in future works.   

The application of the novel elements is demonstrated with an example in the relaxed micromorphic model, where the dilatation between two materials jumps at the interface due to differing materials. One observes the excellent behaviour of the new elements by the comparison to N\'ed\'elec based formulations and the capacity of the novel elements to correctly depict the jump of the trace of the microdistortion field at the interface. With the new elements at hand, further investigations of micro-structured materials in the relaxed micromorphic model are made possible.  

\section*{Acknowledgements}
Patrizio Neff acknowledges support in the framework of the DFG-Priority Programme
2256 “Variational Methods for Predicting Complex Phenomena in Engineering Structures and Materials”, Neff 902/10-1, Project-No. 440935806. Michael Neunteufel acknowledges support by the Austrian Science Fund (FWF) project F65.

\bibliographystyle{spmpsci}   

\footnotesize{
\bibliography{Ref}   
}


\end{document}